\documentclass[a4paper,10pt,twoside]{article}

\usepackage{graphicx}
\usepackage{caption}
\usepackage{subcaption}
\usepackage{multirow}
\usepackage[utf8]{inputenc}
\usepackage[T1]{fontenc}
\usepackage{ dsfont }
\usepackage[english]{babel}
\usepackage{charter}

\usepackage{indentfirst}
\usepackage{xcolor}
\usepackage{verbatim}
\usepackage{hyperref}
\usepackage{pifont}

\usepackage[top=3.cm, bottom=4.0cm, left=2.2cm, right=2.2cm]{geometry}
\usepackage{amsmath}
\usepackage{amsthm}	
\usepackage{amsfonts}	
\usepackage{amssymb	
            ,bbm
            ,units 
            ,stmaryrd
           }

\usepackage{enumerate}
\usepackage[square,numbers,sort&compress]{natbib}
\usepackage{
            ulem		% use \sout{...}
           ,soul		%  use \st{..}
} 

\numberwithin{equation}{section}
\numberwithin{figure}{section}
 \usepackage[nodayofweek]{datetime}

\renewcommand*{\thefootnote}{\fnsymbol{footnote}}

\title{ A flexible split-step scheme for solving McKean-Vlasov Stochastic Differential Equations }

\author{
\normalsize Xingyuan Chen\textit{$^{a}$} \\
        \small   X.Chen-176@sms.ed.ac.uk 
\and
 \normalsize Gon\c calo dos Reis\textit{$^{a,b,}$}\footnote{G.d.R. acknowledges support from the \emph{Funda{\c c}$\tilde{\text{a}}$o para a Ci$\hat{e}$ncia e a Tecnologia} (Portuguese Foundation for Science and Technology) through the project UIDB/00297/2020 and UIDP/00297/2020 (Center for Mathematics and Applications, CMA/FCT/UNL).} \\
        \small  G.dosReis@ed.ac.uk
}

\date{%
    \footnotesize 
    $^{a}$~School of Mathematics, University of Edinburgh, The King's Buildings, Edinburgh, UK
    \\
    $^{b}$~Centro de Matem\'atica e Aplica\c c$\tilde{\text{o}}$es (CMA), FCT, UNL, Portugal
    \\
    \longdate \today \ (\currenttime)
    \vspace{-1.0cm}
}

\theoremstyle{plain}
\newtheorem{Theorem}{Theorem}[section]
\newtheorem{lemma}[Theorem]{Lemma}
\newtheorem{Proposition}[Theorem]{Proposition}

\newtheorem{definition}[Theorem]{Definition}

\newtheorem{remark}[Theorem]{Remark}

\newtheorem{assumption}[Theorem]{Assumption}

\newcommand{\bE}{\mathbb{E}}
\newcommand{\bF}{\mathbb{F}}

\newcommand{\bN}{\mathbb{N}}
\newcommand{\bP}{\mathbb{P}}

\newcommand{\bR}{\mathbb{R}}
\newcommand{\bS}{\mathbb{S}}

\newcommand{\cF}{\mathcal{F}}

\newcommand{\cP}{\mathcal{P}}

\definecolor{darkgreen}{rgb}{0,0.35,0}

\newcommand{\1}{\mathbbm{1}}

\newcommand{\tn}{{t_n}}
\newcommand{\tnp}{{t_{n+1}}}

\newcommand{\dd}{\mathrm{d}}

\newcommand{\hx}{ \hat{X} }
\newcommand{\hz}{ \hat{Z} }
\newcommand{\kt}{ \kappa(t) }
\newcommand{\hm}{ \hat{\mu} }

\newcommand{\wm}{ \widetilde {\mu}}

\begin{document}

\selectlanguage{english}

\maketitle

%%%%%%%%%%%%%%%%%%%%%%%%%%%%%%%%%%%%%%%%%%%%%%%%%%%%%%%%%%%%%%%%%
%%%%%%%%%%%%%%%%%%%%%%%%%%%%%%%%%%%%%%%%%%%%%%%%%%%%%%%%%%%%%%%%%
%%%%%%%%%%%%%%%%%%%%%%%%%%%%%%%%%%%%%%%%%%%%%%%%%%%%%%%%%%%%%%%%%
%%%%%%%%%%%%%%%%%%%%%%%%%%%%%%%%%%%%%%%%%%%%%%%%%%%%%%%%%%%%%%%%%
%%%%%%%%%%%%%%%%%%%%%%%%%%%%%%%%%%%%%%%%%%%%%%%%%%%%%%%%%%%%%%%%%
\renewcommand*{\thefootnote}{\arabic{footnote}}

%%%%%%%%%%%%%%%%%%%%%%%%%%%%%%%%%
\begin{abstract}
We present an implicit Split-Step explicit Euler type Method (dubbed SSM) for the simulation of McKean–Vlasov Stochastic Differential Equations (MV-SDEs) with drifts of superlinear growth in space, Lipschitz in measure and non-constant Lipschitz diffusion coefficient. The scheme is designed to leverage the structure induced by the interacting particle approximation system, including parallel implementation and the solvability of the implicit equation.

The scheme attains the classical $1/2$ root mean square error (rMSE) convergence rate in stepsize and closes the gap left by \cite{reis2018simulation} regarding efficient implicit methods and their convergence rate for this class of McKean-Vlasov SDEs. A  sufficient condition for mean-square contractivity of the scheme is presented. Several numerical examples are presented, including a comparative analysis to other known algorithms for this class (Taming and Adaptive time-stepping) across parallel and non-parallel implementations. 
\end{abstract}
%%%%%%%%%%%%%%%%%%%%%%%%%%%%%%%%%
{\bf Keywords:} 
McKean-Vlasov equations, split-step methods, interacting particle systems, superlinear growth 2000 MSC: 65C05 , 65C30 , 65C35

%
%
%
%%%%%%%%%%%%%%%%%%%%%%%%%%%%%%%%%%%%%%%%%%%%%%%%%%%%%%%%%%%%%%
%%%%%%%%%%%%%%%%%%%%%%%%%%%%%%%%%%%%%%%%%%%%%%%%%%%%%%%%%%%%%%
%%%%%%%%%%%%%%%%%%%%%%%%%%%%%%%%%%%%%%%%%%%%%%%%%%%%%%%%%%%%%%
%%%% \BEGIN SECTION
%%%%%%%%%%%%%%%%%%%%%%%%%%%%%%%%%%%%%%%%%%%%%%%%%%%%%%%%%%%%%%
% \newpage 
\footnotesize
\setcounter{tocdepth}{2}
\tableofcontents
\normalsize
% \footnotesize

%%%%%%%%%%%%%%%%%%%%%%%%%%%%%%%%%%%%%%%%%%%%%%%%%%%%%%%%%%%%%%%%%%%%	
%%%%%%%%%%%%%%%%%%%%%%%%%%%%%%%%%%%%%%%%%%%%%%%%%%%%%%%%%%%%%%%%%%%%
%%% \BEGIN SECTION
%%%%%%%%%%%%%%%%%%%%%%%%%%%%%%%%%%%%%%%%%%%%%%%%%%%%%%%%%%%%%%%%%%%%
\newpage

\section{Introduction}
\label{sec:one}

The aim of this paper is to present a numerical scheme for simulating McKean-Vlasov Stochastic Differential Equations (MV-SDEs) with drifts of superlinear growth in space and Lipschitz in the measure component, and Lipschitz diffusion coefficients. MV-SDEs differ from standard SDEs by means of the presence of the law of the solution process in the coefficients and their dynamics are of the following type
 
\begin{equation}
\label{Eq:General MVSDE}
\dd X_{t} = \Big( v(t,X_{t}, \mu_{t}^{X}) + b(t,X_{t}, \mu_{t}^{X})\Big)\dd t + \sigma(t,X_{t}, \mu_{t}^{X})\dd W_{t},  \quad X_{0} \in L_{0}^{m}( \bR^{d}).
\end{equation}
where $\mu_{t}^{X}$ denotes the law of the solution process $X$ at time $t$, $W$ is a multidimensional Brownian motion, $v,b, \sigma$ are measurable maps and $X_0$ is a sufficiently integrable initial condition (in $L_{0}^{m}$ for $m\geq 2$). 
These equations were introduced by McKean in the sixties \cite{McKean1966} and have been the target of much research since. There is a rich literature on well-posedness \cite{Sznitman1991,BauerMartin2018SSoM,mishura2020existence,rockner2018well} and we point to the recent summary work \cite{huang2020distribution} highlighting  recent developments in regularity estimates, exponential ergodicity,
long time large deviations (also \cite{dosReisSalkeldTugaut2017,adams2020large}), comparison Theorems and the Vlasov-Fokker-Plank equations associated to MV-SDEs. One particular element of interest is the so-called propagation of chaos (PoC) introduced by Kac \cite{kac1956} and further studied in the MV-SDE literature \cite{Sznitman1991,Meleard1996,Lacker2018,lacker2021hierarchies}. The PoC phenomena states that a MV-SDE is the limit of a certain weakly interacting particle systems (of standard SDEs) as the system's size increases to infinity. Namely, \eqref{Eq:General MVSDE} is the limit, as $N\to \infty$, of the $N$-dimensional system of $\bR^d$-valued interacting particles $X^{i,N}$ , 
\begin{align}
\label{Eq:MV-SDE Propagation}
\dd {X}_{t}^{i,N} 
= \Big( v\Big(t,X_t^{i,N}, \mu^{X,N}_{t}\Big)+ b\Big(t,{X}_{t}^{i,N}, \mu^{X,N}_{t} \Big)\Big) \dd t 
+ \sigma\Big(t,{X}_{t}^{i,N} , \mu^{X,N}_{t} \Big) \dd W_{t}^{i}
, \qquad X^{i,N}_0=X_0^i.
\end{align}
with $\mu^{X,N}_{t}$ being the empirical measure given as $\mu^{X,N}_{t}(\dd x) := \frac{1}{N} \sum_{j=1}^N \delta_{X_{t}^{j,N}}(\dd x)$ and $\delta_{{X}_{t}^{j,N}}$ is the Dirac measure at point ${X}_{t}^{j,N}$, $W^i$ are independent Brownian motions  and with independent and identically distributed initial conditions $X_{0}^{i}$ across $i=1,\cdots,N$. In essence, the law $\mu^X_t$ in $\bR^d$ of \eqref{Eq:General MVSDE} is approximated by the empirical average $\mu^{X,N}_{t}$ generated by the $(\bR^d)^N$-system (a high-dimensional system). 
This methodology, appealing to the interacting particle system, includes a well-known quantified speed of convergence result (summarised in Proposition \ref{eq:momentboundParticiInteractingSystem} below) providing a path for a numerical approximation: from $X$ to $\{X^{i,N}\}_i$ to $\{X^{i,N,\pi}\}_i$ with $X^{\cdot,N,\pi}$ the numerical approximation of each particle in the SDE system.

It is essential to highlight that although $\{X^{i,N}\}_i$ is a high-dimensional SDE and any existing numerical method a priori applies straightforwardly, all rates and results obtained in this straightforwardly way depend on the system's dimension $Nd$. The constants then explode as $N$ increases. Part of the difficulty of the method is to show that such constants, rates and results are independent of $N$ albeit depending on $d$. This step has its non-trivial difficulties which we discuss further below by drawing on prior work \cite{reis2018simulation}.

In this work, we focus on the class of MV-SDE with drifts of superlinear growth in their spatial components \cite{dosReisSalkeldTugaut2017,reis2018simulation,reisinger2020adaptive,bao2020milstein,kumar2020explicit}, encapsulated in the function $v$, where $b,\sigma$ are uniformly Lipschitz (in space and measure), and $\sigma$ is non-constant -- this structural Assumption on $b,v$ is trivial from the theoretical perspective but plays an important role in the numerics. 
This class of MV-SDEs appears in several practical models in science, for example, in neuroscience \cite{Baladron2012,Nnmodelcite} introduce the mean-field FitzHugh-Nagumo model for a neuron networks in the brain; \cite{Bolley2011} discuss individual-based and swarming Cucker-Smale interaction models; and models of battery electrodes  \cite{Dreyer2011phasetransition,Guhlke2018stochmanyparticles}. These equations do not have explicit or closed-form solutions and numerical approximations are needed. Moreover, standard explicit numerical methods suitable for the Lipschitz case fail to converge on the superlinear growth setting. This is exemplified by the `particle corruption' phenomena \cite[{Section 4.1}]{reis2018simulation} for MV-SDE particle system numerics. This phenomenon is akin to the divergence of SDE schemes in the superlinear growth setting as highlighted in the seminal work  \cite{HutzenthalerEtAl2011}.

The numerical approximation of McKean–Vlasov equations in the continuous case was initiated in \cite{bossytalay1997} and has been investigated further in several recent works. We briefly mention a few on numerical schemes for MV-SDEs under the superlinear setting and in the Brownian framework. Tamed Euler schemes appeared first \cite{reis2018simulation}, shortly followed by tamed Milstein schemes \cite{bao2021first,kumar2020explicit} (appeared simultaneously) and Milstein schemes for delay MV-SDEs \cite{bao2020milstein}. In \cite{kumar2020well} a tamed scheme is proposed (and a new wellposedness result) for MV-SDE featuring common noise in the particle system (which \eqref{Eq:MV-SDE Propagation} does not) and where the diffusion is also allowed to grow superlinearly. Adaptive time-stepping methods come as an alternative to taming and in the context of MV-SDE they are proposed in \cite{reisinger2020adaptive} --- these two methods will  henceforth be referred to as the `\textit{Taming}' and the `\textit{Adaptive}' algorithm, respectively. 

Outside the superlinear setting, \cite{hutzenthaler2021multilevel} discusses  computational complexity of MV-SDE algorithms (uniformly Lipschitz drift and constant diffusion coefficient) and we point the reader there for an in-depth overview on the state of the art in that regard. Numerical approximations for MV-SDEs with non-Lipschitz conditions in measure and space exist \cite{ding2020euler} but impose  a linear growth condition on the coefficients.  
The case of simulating MV-SDEs with discontinuous coefficients has been addressed \cite{leobacher2020well}. 
An alternative to the empirical measure approximation of \eqref{Eq:MV-SDE Propagation} is to use a projection-type estimation of the marginal densities \cite{belomestny2018projected} where the error analysis requires differentiability of the coefficients. 
Variance reduction technique have been analysed for the class of MV-SDE, namely, importance sampling \cite{reis2018importance}, antithetic multilevel Monte Carlo sampling \cite{bao2020milstein} and antithetic sampling \cite{bencheikh2019bias}. There also recent progress in the jump-diffusion setting \cite{agarwal2021fourier,biswas2020well}.
\medskip

\textbf{Motivation.} For the superlinear growth case described above, we are motivated by an open question left in \cite{reis2018simulation} regarding implicit-type numerical methods for MV-SDEs. All methods described above are of explicit time-stepping nature which are known to lose some of the geometric properties of the original system. For instance, taming destroys the strict dissipativity of the drift map which then raises questions regarding the stability of the scheme's output across long time horizons (see \cite{zong2014stabilitytaming}). Implicit methods for MV-SDEs are largely unexplored. The notable exceptions are \cite{Malrieu2003,reis2018simulation} which are also starting point for this work. In \cite{Malrieu2003} the authors study MV-SDEs and associated particle systems with drifts of (symmetric) convolution kernel-type (as in \cite{adams2020large}) and constant diffusion coefficient. The theoretical (Bacry-Emery) machinery employed there yields a critical logarithmic Sobolev inequality estimate that is used to establish a concentration inequality for their implicit Euler scheme. Assumption-wise, their setting and the setting of this work do not cover each other but agree over a small class. Further, it is unclear how to extend their methodology to non-constant possibly degenerate diffusion coefficients.

More recently, an implicit method to deal with the superlinear growth was proposed \cite{reis2018simulation} (for general diffusion coefficients and without a concavity assumption) but convergence was shown under stronger restrictions than expected: the measure component of the drift is Lipschitz in Wasserstein-$1$ metric not just in Wasserstein-$2$ (plus uniformly bounded measure dependency). At the core of these difficulties was the use of stopping times arguments which revealed themselves difficult to handle with the measure dependency. The critical point is the proof of Lemma 5.12   (p.41) in \cite{reis2018simulation} and the calculations executed in (p.48-49). Lastly, upon inspection of that proof, we emphasise that general $\theta$-methods for SDEs \cite{higham2002strong} will face the same difficulties.
\medskip

\textbf{Our contributions.} In this work, we revisit the framework of \cite{reis2018simulation} and propose a split-step numerical scheme inspired by the earlier work \cite{higham2002strong}. Our contributions can be summarised as follows, 
\begin{enumerate}[(I)]
    \item \textit{Main results.} We proposed a split-step method (SSM), see \eqref{eq:SSTM:scheme 1}-\eqref{eq:SSTM:scheme 2} below, for this class of MV-SDEs. We prove its convergence and recover the $1/2$-convergence rate in root Mean Square Error (rMSE) under the same general assumptions as Taming \cite{reis2018simulation} or Adaptive \cite{reisinger2020adaptive}. No differentability or non-degeneracy assumptions are imposed and stopping-time arguments are fully avoided. 
    
    We provide the stability analysis of Mean-square contractivity for the SSM (non-constant diffusion coefficient). To the best of our knowledge this has not yet been discussed for MV-SDE schemes in general. The stability of the SSM provides a theoretical foundation for carrying out simulation with larger timestep and we point to positive results by way of numerical simulation with the Cucker-Smale flocking model \cite{csreference} where the SSM outperforms both Taming \cite{reis2018simulation} and Adaptive  \cite{reisinger2020adaptive} algorithms.
    
    In regards to known findings, the SSM here overcomes the limitations of the implicit method proposed in \cite[Section 3.2]{reis2018simulation} and extends its scope of application. In the context of standard SDE simulation    (where the coefficients are independent of the measure), our results lift the differentiability restriction of \cite[Assumption 3.1]{higham2002strong}, allow for time-dependence and can take advantage of (the possible) concavity of the map $v$ (in \eqref{Eq:General MVSDE}) -- this is a mild improvement of known SDE results.
 
    \item \textit{Computational gains.} The scheme is designed so that the superlinear term can be split from the main equation in a way that   optimises/minimises  the computational cost of the inversion method (from the implicit component). This flexibility is understood in the following way: given an MV-SDE it is left to the user to choose which terms form $v$ and which terms form $b$ in \eqref{Eq:General MVSDE} (within restrictions). Moreover, since there is a split $v$ Vs $b$ the user may decide to add \& subtract convenient terms to the drift (see Section \ref{sec:stabilityginzburgLandau} below; also \cite[Eqs.~(37) and (38)]{buckwarbrehier2021FHNmodelandsplittingBSTT2020})  -- this trick allows one to transform a non-dissipative term $x\mapsto v(\cdot,x,\cdot)$ into a dissipative one at the expense of an increase of the Lipschitz constant of $b$.

    The SSM proposed allows to decouple the measure component making it amenable to a parallel implementation (e.g., \cite{2019PCT1}). Namely, one parallelizes the task of solving $N$-times an $\bR^d$-system in opposition to the non-parallelizable task of solving once the $\bR^{dN}$-system. The computational gains of the parallel implementations for the SMM are shown to be on par with Taming \cite{reis2018simulation} and Adaptive  \cite{reisinger2020adaptive} algorithms.

    \item \textit{Comparative study against known literature.} We provide a comparative analysis against the Taming \cite{reis2018simulation} and Adaptive \cite{reisinger2020adaptive} methods, across parallel and non-parallel implementations. The numerical study covers four examples of interest, highlighting a different flavour of these 3 algorithms including stability experiments.
\end{enumerate}
We close this introduction  with two comments. The general setting of \cite{dosReisSalkeldTugaut2017,reis2018simulation,reisinger2020adaptive,bao2020milstein,kumar2020explicit} is based on drifts maps $(t,x,\mu)\mapsto \widehat b(t,x,\mu)$ and it is this map that satisfies a one-sided Lipschitz condition in space and a uniform Lipschitz condition in measure. In \eqref{Eq:General MVSDE}, we specify $\widehat b$ to be $(t,x,\mu)\mapsto \widehat b(t,x,\mu) = v(t,x,\mu)+b(t,x,\mu)$ where the polynomial growth is fully captured by $v$ while $b$ remains uniformly Lipschitz in its variables. Separating $\widehat b$ into $v$ and $b$ is natural non-limiting assumption.  

Settings outside the scope of this work are: non-Lipschitz measure dependency \cite{Malrieu2003,adams2020large},  the  superlinear diffusion case of \cite{kumar2020well} (also \cite[Section 3.1]{reisinger2020adaptive}), common-noise \cite{kumar2020well} or jumps \cite{biswas2020well}. Weak error analysis is left unaddressed.

\medskip

\textbf{Organisation of the paper.} In section 2, notations and necessary  concepts for this work are given. In Section \ref{sec:theSSMresults} we state the SSM, the main Theorem of convergence and the stability results for the scheme. Section 3 provides several numerical examples for comparison with other methods. Examples covering the stability analysis in the superlinear case, notably the Cucker-Smale model, are also given in Section 3. All proofs are postponed to  Section 4. For convenience, \ref{sec:ShortRecap} contains a short description of the Taming and Adaptive algorithms.

\medskip

\textbf{Acknowledgements.} The authors would like to thank the 3 referees for their thorough work and suggestions that led to non trivial improvements.

%%%%%%%%%%%%%%%%%%%%%%%%%%%%%%%%%%%%%%%%%%%%%%%%%%%%%%%%%%%%%
%%%% \END SECTION
%%%%%%%%%%%%%%%%%%%%%%%%%%%%%%%%%%%%%%%%%%%%%%%%%%%%%%%%%%%%%
%%%%%%%%%%%%%%%%%%%%%%%%%%%%%%%%%%%%%%%%%%%%%%%%%%%%%%%%%%%%%
%%%%%%%%%%%%%%%%%%%%%%%%%%%%%%%%%%%%%%%%%%%%%%%%%%%%%%%%%%%%%
%
%
%

%
%
%
%%%%%%%%%%%%%%%%%%%%%%%%%%%%%%%%%%%%%%%%%%%%%%%%%%%%%%%%%%%%%%%%%%%%	
%%%%%%%%%%%%%%%%%%%%%%%%%%%%%%%%%%%%%%%%%%%%%%%%%%%%%%%%%%%%%%%%%%%%
%%% \BEGIN SECTION
%%%%%%%%%%%%%%%%%%%%%%%%%%%%%%%%%%%%%%%%%%%%%%%%%%%%%%%%%%%%%%%%%%%%
% \newpage
\section{The split-step methods for MV-SDEs}
\label{sec:two}

\subsection{Notation and Spaces}
\label{sec:twoNotation}

Let $\bN=\{1,2,\cdots\}$ be the set of natural numbers starting at $1$ and $\bR$ denotes the real numbers where $\bR^+=[0,\infty)$. Also, we denote $\llbracket a,b\rrbracket:= [a,b] \cap \bN =  \{a,\cdots,b\}$,  for any $ a,b\in \bN$ with $a\leq b$. 
% For some $N \in \bN$, define the set $\llbracket 1,N \rrbracket =\{1,2,\cdots,N\}=[1,N]\cap \bN$. 
For $x,y \in \bR^d$ denote the scalar product of vectors by $\langle x, y \rangle$; and the Euclidean distance of $x$ is $|x|=(\sum_{j=1}^d x_j^2)^{1/2}$. 
 The indicator function of a set $\Omega$ is denoted as $\1_\Omega$.

The space of probability measures on $\bR^d$ is denoted by $\cP(\bR^d)$ and its subset of  finite second moment measures is denoted by $\cP_2(\bR^d)$. We recall the definition of the Wasserstein-$2$ distance metrizing $\cP_2(\bR^d)$, 
\begin{align*}
W^{(2)}(\mu,\nu) = \inf_{\pi\in\Pi(\mu,\nu)} \Big(\int_{\bR^d\times \bR^d} |x-y|^2\pi(dx,dy)\Big)^\frac12, \quad \textrm{for any }\mu,\nu\in \cP_2(\bR^d) .
\end{align*}   
where $\Pi(\mu,\nu)$   denotes all probability measure couplings between $\mu$ and $\nu$ on $\bR^d\times \bR^d$, i.e.~$\pi\in\Pi(\mu,\nu)$ if and only if $\pi(\cdot\times \bR^d)=\mu$ and $\pi(\bR^d \times \cdot)=\nu$.

Let $(\Omega, \bF, \cF,\bP)$ be a %n atomless Polish 
filtered probability space with $ \cF_t$ the augmented filtration generated by a standard $l$-dimensional Brownian motion $W=(W^1,\cdots,W^l)$ and with an additionally sufficiently rich sub $\sigma$-algebra $\cF_0$ independent of $W$.  
We use $\bE[\cdot]=\bE^\bP[\cdot]$ to denote the  expectation operator under the measure $\bP$. 
As usual notations, ``i.i.d.'' means ``independent and identically distributed'' and ``w.r.t'' means ``with respect to''.

Fix $T\in (0,\infty)$. We follow the notation from \cite{reis2018simulation}.  Let $p\geq 1$. Define $L_{t}^{p}(\mathbb{R}^{d})$ as the space of $\mathcal{F}_{t}$-measurable $\mathbb{R}^{d}$-valued random variables $X$ that satisfy $\mathbb{E}\Big[\,|X|^{p}\Big]^{1 / p}<\infty$.  We use $\mathbb{S}^{p}$ to denote the space of $\bR^d$-valued $\cF_\cdot$-adapted processes $Z$ satisfying  $\mathbb{E}\Big[\sup _{0 \leqslant t \leqslant T}|Z_t |^{p}\Big]^{1 / p}<\infty$. The cross-variation between two processes $X$ and $Y$ is denoted as   $\langle X,Y \rangle$.

Throughout the text $C\in\mathbb{R^+}$ is a constant that may change from line to line, may depend on the problem's data but is always independent of the constants $h,M,N$ (associated with the numerical scheme and specified below).

\subsection{Framework}

We study MV-SDE \eqref{Eq:General MVSDE} for $m \geq 1$ and we denote the law of $X$ at time $t$ as $\mu_{t}^{X}$. Take $b,v,\sigma$ as measurable functions, $v:[0,T] \times \bR^d \times\cP_2(\bR^d) \to \bR^d$, $b:[0,T] \times \bR^d \times\cP_2(\bR^d) \to \bR^d$ and $\sigma:[0,T] \times \bR^d \times \cP_2(\bR^d) \to \bR^{d\times l}$. Throughout the text, we make the following assumption.

\begin{assumption}
\label{Ass:Monotone Assumption}
Assume that $v,b$ and $\sigma$ are $1/2$-H\"{o}lder continuous in time, uniformly in $x\in \bR^d$ and $\mu\in \cP_2(\bR^d)$. Assume that $b,\sigma$ are uniformly Lipschitz in the sense that there exists $L_b, L_\sigma\geq 0$ such that for all $t \in[0,T]$ all $x, x'\in \bR^d$ and all $\mu, \mu'\in \cP_2(\bR^d)$ we have
\begin{align*}
|b(t, x, \mu)-b(t, x', \mu')|^2
&\leq L_b(|x-x'|^2 + W^{(2)}(\mu, \mu')^2 ),
\\
|\sigma(t, x, \mu)-\sigma(t, x', \mu')|^2&\leq L_\sigma (|x-x'|^2 + W^{(2)}(\mu, \mu')^2 ).
\end{align*}
For $v$, there exist $ L_{v} \in \bR$, $L_{\hat{v}}>0$, $q\in \bN$ and $q>1$ such that for all $t\in[0,T]$, $ x, x'\in \bR^d$, all $\mu, \mu'\in \cP_2(\bR^d)$ we have 
\begin{align*}
\langle x-x', v(t, x,\mu)-v(t, x',\mu) \rangle 
& 
\leq L_{v}|x-x'|^{2},
& \textrm{(One-sided Lipschitz in space)},
\\
|v(t, x,\mu)-v(t, x',\mu)| 
& 
\leq L_{\hat{v}}(1+ |x|^{q} + |x'|^{q}) |x-x'| ,
& \textrm{(Locally Lipschitz in space)},
\\
|v(t, x,\mu)-v(t, x,\mu')|^2 
&
\leq
{L_{\tilde{v}} W^{(2)}(\mu, \mu')^2 },
& \textrm{(Uniformly Lipschitz in measure)}.
\end{align*}
\end{assumption} 
The structural choice of having a drift $\widehat b = v+b$ with only $v$ containing the superlinear growth component, as opposed to a single drift map $\widehat b$ in the style of \cite{reis2018simulation,reisinger2020adaptive}, is negligible and its use is discussed in Remark \ref{remark:removing constraint on Lv}.

Immediate well-known properties can be derived from this assumption.
\begin{remark}[Implied properties]
\label{rem:ImpliedProperties}
Under Assumption \ref{Ass:Monotone Assumption}, define $ \widehat L_v  ={1}/{2}+L_v$ and  
 $C_T=\sup_{t\in[0,T]}|v(t,0,\delta_0)|^2/2$. Let $C>0$, then for all $t \in [0,T]$, $x\in \bR^{d}$ and $\mu\in \cP_2(\bR^d)$ one has 
\begin{align*}
& \langle x,v(t,x,\mu)\rangle 
 = \langle x-0,v(t,x,\mu)-v(t,0,\mu)\rangle+\langle x,v(t,0,\mu)\rangle 
\\
&
\le L_v|x|^2+|x||v(t,0,\mu)|
\le (L_v+\frac{1}{2})|x|^2 +\frac{1}{2}|v(t,0,\mu)-v(t,0,\delta_0)|^2+\frac{1}{2}|v(t,0,\delta_0)|^2
% \\
% &
\le C_T+ \widehat L_v |x|^2+\frac{L_{\tilde{v}} }{2} W^{(2)}(\mu,\delta_0)^2. 
\end{align*}
  where the last step follows using Young's inequality.  Additionally, for $\psi\in\{b,\sigma\}$ one has  
\begin{align*}
\Big\langle x,\psi(t,x,\mu) \Big\rangle 
\le
C\Big(1+|x|^2+W^{(2)}(\mu,\delta_0)^2 \Big)
\quad \textrm{and}\quad 
  |\psi(t,x,\mu)|^2
  \le
   C\Big(1+|x|^2+W^{(2)}(\mu,\delta_0)^2 \Big).
\end{align*}
\end{remark}
The above assumptions cover a larger range of models as highlighted in Section \ref{sec:examples} below. This setting subsumes the standard globally Lipschitz assumptions.  For further examples we point to  \cite{dosReisSalkeldTugaut2017,Bolley2011,gomes2019mean,malrieu2006concentration}.

\begin{Theorem}[Theorem 3.3 in \cite{dosReisSalkeldTugaut2017}]
\label{Thm:MV Monotone Existence}
	Let Assumption \ref{Ass:Monotone Assumption} hold and suppose we have $X_{0} \in L_{0}^{m}(\bR^{d})$ for some fixed $m \ge 2$. 
 	Then, there exists a unique solution $X$ to
    the MV-SDE \eqref{Eq:General MVSDE} and it satisfies $X \in \bS^{m}([0,T])$. 
  
    There exists a constant $C\in\mathbb{R^+}$ such that
 	\begin{align*}
	\mathbb E \Big[ \sup_{t\in[0,T]} |X_{t}|^{m} \Big] 
	\leq C \left(1+ \bE[\, |X_0|^m]\right) e^{C T}.
	\end{align*} 
\end{Theorem}
 
\medskip

\textbf{The interacting particle system approximation.}
In order to approximate of MV-SDE \eqref{Eq:General MVSDE}, we build an interacting particle system  
as follows: 
\begin{enumerate}
    \item For ${i\in \llbracket 1,N\rrbracket}$, take ${X}_{0}^{i,N}=X_{0}^{i}$ as i.i.d.~random initial condition for each particle.
    \item Each particle is driven by its own independent Brownian motion $W^i$ (all are i.i.d.)
    \item The dynamics of the particle system is given by   \eqref{Eq:MV-SDE Propagation}, and we set $\mu^{X,N}_{t}(\dd x) := \frac{1}{N} \sum_{j=1}^N \delta_{X_{t}^{j,N}}(\dd x)$ where $\delta_{x}$ is the Dirac measure at point $x\in \bR^d$.
\end{enumerate}

\medskip

\textbf{Propagation of chaos (PoC).}

Below we show a pathwise PoC result to control the difference between the original MV-SDE and the interacting particle system. For that, we introduce {   the auxiliary equation system of \textit{non interacting particles}}
\begin{align}
	\label{Eq:Non interacting particles}
	\dd X_{t}^{i} = 
	\Big( 
	v(t,X_{t}^{i}, \mu^{X^{i}}_{t})+ b(t, X_{t}^{i}, \mu^{X^{i}}_{t}) \Big)\dd t + \sigma(t,X_{t}^{i}, \mu^{X^{i}}_{t}) \dd W_{t}^{i}, \quad X_{0}^{i}=X_{0}^{i} \, ,\quad t\in [0,T].
\end{align}
This system is just $N$ independent MV-SDEs (each in $\bR^d$). The $X^{i}$s are independent, we have $\mu^{X^{i}}_{t}=\mu^{X}_{t}$ (for all $i \in \llbracket 1,N \rrbracket$) where $\mu^{X}_{t}$ is the law of $X$ solution to  \eqref{Eq:General MVSDE}). Under Assumption \ref{Ass:Monotone Assumption}, we have (see \cite{Sznitman1991,Meleard1996,Lacker2018,Carmona2016})
\begin{align*}
\lim_{N \rightarrow \infty} \sup_{ i\in \llbracket 1,N \rrbracket}
\bE \Big[ 
\sup_{0 \le t \le T} |X_{t}^{i,N} - X_{t}^{i}|^{2}\,
\Big] = 0 \, .
\end{align*}

By showing the following \emph{propagation of chaos} result, we can connect in a quantifiable manner the MV-SDE and the interacting particle system.  
 
The proof can be found in \cite[Theorem 3.1]{reis2018simulation} (and  \cite[Proposition 5.1]{reis2018simulation} for the well-posedness of the particle system \eqref{Eq:MV-SDE Propagation}).

\begin{Proposition}[Propagation of chaos]
	\label{Prop:Propagation of Chaos}
	Let Assumption \ref{Ass:Monotone Assumption} hold and suppose we have for some $m \ge 2$ for all $i \in \llbracket 1,N \rrbracket$ that $X^i_{0} \in L_{0}^{m}(\bR^{d})$. 
	Then there exists a unique solution $\{X^{i,N}\}_i$ to \eqref{Eq:MV-SDE Propagation} in $\bS^{m}([0,T])$ and for any $1\leq p\leq m$ there exists $C_p\in\mathbb{R^+}$ such that

\begin{align}
\label{eq:momentboundParticiInteractingSystem}
    \sup_{ i\in \llbracket 1,N \rrbracket} \bE[ \sup_{t\in[0,T]}  |X^{i,N}_t|^p ] \leq C_p(1+\sup_{ i\in \llbracket 1,N \rrbracket}  \bE[\,|X_0^{i,N}|^p ]).
\end{align}
 	
Moreover,  let $X^{i}\in \bS^m$ satisfy \eqref{Eq:Non interacting particles}
 and assume $m>4$. Then, the convergence rate between MV-SDE \eqref{Eq:Non interacting particles} and the interacting particle system \eqref{Eq:MV-SDE Propagation} is given by
	\begin{align*}
	\sup_{ i\in \llbracket 1,N \rrbracket} 
	\bE[\sup_{0 \le t \le T} |X_{t}^{i} - X_{t}^{i,N}|^{2}] 
	\le C
	\begin{cases}
	N^{-1/2} &, \text{if } d<4,
	\\
	N^{-1/2} \log(N) \quad &, \text{if } d=4,
	\\
	N^{-2/d} &, \text{if } d>4.
	\end{cases} 
	\end{align*}
\end{Proposition}

Under this result, one can approximate MV-SDEs through particle scheme. Thus, by showing the convergence of the numerical methods to the interacting particle scheme, we can obtain the convergence rate between a numerical method and the MV-SDE.
\begin{remark}[Optimising the PoC convergence rate]
\label{rem:BetterConvRate}
We highlight a result from \cite{delarue2018masterCLT} and later reviewed in \cite{reisinger2020adaptive} in the context of numerical methods for MV-SDEs. The PoC rate can be improved for the case $d=4$, namely the $\log(N)$ term can be omitted (under restrictions). This holds under the additional constraint of a constant diffusion coefficient $\sigma$, and a bounded drift with bounded derivatives. It does not cover the superlinear growth case here, nonetheless, it is reasonable to expect that the result can be lifted to match the drift condition in this work and a diffusion coefficient that is bounded (and sufficiently smooth). 
\end{remark}

\subsection{The split-step method (SSM) for MV-SDEs: convergence and stability}
\label{sec:theSSMresults}

The numerical scheme proposed in this work, and dubbed \textit{Split-Step Method} (SSM), improves strongly on the implicit numerical scheme proposed in \cite{reis2018simulation}. It is an enhanced variant of the split-step backward Euler scheme for standard SDEs \cite[Eq. (3.8)-(3.9)]{higham2002strong} and here further optimised for the MV-SDE setting. 

Define the uniform partition of $[0,T]$ as $\pi:=\{t_n:=nh : n\in \llbracket 0,M\rrbracket, h:=T/M \}$ for a prescribed $M\in \bN$. Define recursively the split-step method to approximate \eqref{Eq:MV-SDE Propagation} as follows: for $i\in \llbracket 1,N\rrbracket $ set $\hx_{0}^{i,N}=X^i_0$, then, iteratively over  $n\in \llbracket 0,M-1\rrbracket$ for all $i\in \llbracket 1,N\rrbracket $, 
\begin{align}
\label{eq:SSTM:scheme 1}
Y_{n}^{i,\star,N} &=\hx_{n}^{i,N}+h v\left(t_n,Y_{n}^{i,\star,N},\hm^{X,N}_n\right),  
 \quad 
 \quad 
 \\
\label{eq:SSTM:scheme 2}
\hx_{n+1}^{i,N} &=Y_{n}^{i,\star,N}
            + b(t_n,Y_{n}^{i,\star,N},\hm^{Y,N}_n) h
            +\sigma(t_n,Y_{n}^{i,\star,N},\hm^{Y,N}_n) \Delta W_{n}^i,
\\ \nonumber
\textrm{where}\quad \hm^{X,N}_n(dx):&= \frac1N \sum_{j=1}^N \delta_{\hx_{n}^{j,N}}(dx),
\quad
\hm^{Y,N}_n(dx):= \frac1N \sum_{j=1}^N \delta_{Y_{n}^{j,\star,N}}(dx).
\end{align}

We state immediately the main convergence result between the continuous time extension of the scheme \eqref{eq:SSTM:scheme 1}-\eqref{eq:SSTM:scheme 2} and the solution to the MV-SDE \eqref{Eq:General MVSDE}. The proof is postponed to Section \ref{sec:ProofOfConv}. 
\begin{Theorem}
\label{theo:SSTM convergence rate for MV-sde}
Let the assumptions of Proposition \ref{Prop:Propagation of Chaos} hold.       Let $m>2(q+1)^2 $, where $q$ is the polynomial growth parameter of Assumption \ref{Ass:Monotone Assumption}.  
Take the collection $\{ {\hx}_n^{i,N} \}_{i,n}$ for $n\in\llbracket 0,M\rrbracket$, $i\in\llbracket 1,N\rrbracket$ generated through the scheme \eqref{eq:SSTM:scheme 1}-\eqref{eq:SSTM:scheme 2} under the timestep constraint expressed through the one-sided Lipschitz constant $L_v$ of $v$ (Assumption \ref{Ass:Monotone Assumption}, see Remark \ref{remark:removing constraint on Lv})
\begin{align*}
\begin{cases}
	  h>0 \ \textrm{ and }\ h \le \frac{1}{1+2 L_v}   &, \text{if }L_v> -\frac{1}{2},
	 \\
	 h>0 &, \text{if } L_v\le -\frac{1}{2}.
	 \end{cases} 
\end{align*}
Then, the following assertions hold. There exists a continuous-time extension,  $(\hx_{t}^{i,N})_{t\in[0,T]}$ to the SSM \eqref{eq:SSTM:scheme 1}-\eqref{eq:SSTM:scheme 2} (and given in \eqref{eq: SSTM discrete scheme sde with kappa t} below), satisfying ($C>0$ is a constant independent of $N,M$ but may depend on $T,d$):
\begin{enumerate}
\item Uniformly bounded $p$ moments. Given  $m\ge 2p\geq 1$ there exist constant $C>0$  such that 

\begin{align*}
    \sup_{ i\in \llbracket 1,N \rrbracket}   &\bE\Big[\sup_{0\le t \le T} |\hx_{t}^{i,N}|^{2p}\Big]
    <C\Big(  1+ \sup_{ i\in \llbracket 1,N \rrbracket}  \bE\Big[\, |\hx_{0}^{i,N}|^{2p}\Big] \Big) 
    <\infty.
\end{align*}

\item  Take $(X^{i,N})_i$ as the solution of the interacting particle system \eqref{Eq:MV-SDE Propagation}. Then, scheme \eqref{eq:SSTM:scheme 1}-\eqref{eq:SSTM:scheme 2} converges to $X^{i,N}$ with a strong global convergence rate of ${1}/{2}$ in root mean square error (rMSE) over $[0,T]$, namely, 
\begin{align}
\label{eq:rMSEofSSM}
    \textrm{rMSE}= \sqrt{\sup_{ i\in \llbracket 1,N \rrbracket}   \bE\Bigg[\sup_{0\le t \le T} |X_{t}^{i,N}-\hx_{t}^{i,N} |^2 \Bigg]}\le C h^{\frac{1}{2}}.
\end{align}

\item 
Let $(X^{i})_i$ be the solution of the non-interacting particle system \eqref{Eq:Non interacting particles}. We have 
    \begin{align*}
	\sup_{ i\in \llbracket 1,N \rrbracket}  \bE \Big[ \sup_{0 \le t \le T} |X_{t}^{i} - \hx_{t}^{i,N}|^{2}
	\Big] 
	\le
	 C
	 \begin{cases}
	 h+ N^{-1/2} & \text{if } d<4,
	 \\
	 h+ N^{-1/2} \log(N) \quad & \text{if } d=4,
	 \\
	 h+N^{-2/d} & \text{if } d>4.
	 \end{cases} 
	\end{align*}
\end{enumerate}

\end{Theorem}
We remind the reader about Remark \ref{rem:BetterConvRate} concerning the PoC's convergence rate. 
 
\begin{remark}[The constraint on $L_v$ is soft and removable.]
\label{remark:removing constraint on Lv}

The choice of what $v$ and $b$ are is left to the user. More precisely, to the drift map $\widehat b(t,x,\mu) = v(t,x,\mu) + b(t,x,\mu)$ one can always add and subtract a linear term $\gamma x$ ($\gamma \in \bR$) and set $\widehat b(t,x,\mu) = \Big(v(t,x,\mu)-\gamma x) + \Big(b(t,x,\mu) + \gamma x\Big)$. This means that the one-sided Lipschitz constant $L_v$ becomes $L_v-\gamma$ and hence negative if $\gamma$ is sufficiently large.

We show below that this operation is not free of cost. Concretely, there is an implication in terms of the scheme's stability since for $L_v$ to become negative the Lipschitz constant $L_b$ increases proportionally. In Section \ref{sec:stabilityginzburgLandau} we discuss this in view of a numerical example and via a mean-square stability result we provide in Theorem \ref{theo:SSTM:stabilty} for the SSM  \eqref{eq:SSTM:scheme 1}-\eqref{eq:SSTM:scheme 2}.

Lastly, to the best of our knowledge, the SSM scheme \eqref{eq:SSTM:scheme 1}-\eqref{eq:SSTM:scheme 2} is not of the usual form SSM schemes are presented in the literature.  Usually there is no structural separation of $v+b$ and one sets $b=0$. Consequently there is no drift component in \eqref{eq:SSTM:scheme 2} (only a diffusion part), thus a discussion on the benefits/drawbacks of adding/subtracting of a $\gamma x$-term seems generally absent. 
\end{remark}

Remark \ref{remark:1-hLv is a constant like 1/2} provides more details on the choice of $h$.  
The constraint of $L_v<-1/2$ is not sharp. In fact, it can be replaced by  $L_v<-\varepsilon$ for some $\varepsilon \in (0,1)$ at the expense of another constant growing proportionally to $1/\varepsilon$. We choose for simplicity $\varepsilon=1/2$, see Remark \ref{remark:1-hLv is a constant like 1/2} and the definition of $\widehat L_v$ in Remark \ref{rem:ImpliedProperties} for further details. This issue is relevant in case one sets $b=0$ as is usual in the SSM literature (and the trick of Remark \ref{remark:removing constraint on Lv} cannot be applied).

Theorem \ref{theo:SSTM convergence rate for MV-sde} shows the strong convergence rate of the SSM is $1/2$ (rMSE) which is the same as Taming \cite{reis2018simulation} and Adaptive \cite{reisinger2020adaptive}. Also, the complexity of the particle system is of order $N^2$ in the worst situation, however, in several examples of section \ref{sec:examples} the complexity for the calculation of the interaction term is of order $N$.

After the convergence study of Theorem \ref{theo:SSTM convergence rate for MV-sde} we introduce the notion of mean-square contractivity and study the stability of the SSM  \eqref{eq:SSTM:scheme 1}-\eqref{eq:SSTM:scheme 2}.

\begin{definition}[Mean-square contractivity]
\label{def:definition of mean-square stable numerics X Y }
Suppose that $X_0\in L_0^m(\mathbb{R}^d)\nonumber$ and $Z_0\in L_0^m(\mathbb{R}^d)\nonumber$ for some sufficiently large $m\in \bN$. 
Take two numerical solutions of the same numerical scheme $\hx_n^{i,N}$ and  $\hz_n^{i,N}$ of \eqref{Eq:General MVSDE} with $\hx_0^i$ and $\hz_0^i$ being i.i.d.~copies of $X_0$ and $Z_0$ respectively. The scheme is called mean-square contractive if we have 
\begin{align*}
\lim_{n\rightarrow \infty}\sup_{ i\in \llbracket 1,N \rrbracket}  \bE\Big[\, |\hx_{n}^{i,N}-\hz_{n}^{i,N} |^2 \Big]=0.
\end{align*}
\end{definition}
The next result shows when the SSM is mean-square contractive.
\begin{Theorem}
\label{theo:SSTM:stabilty}
Let the assumptions of Theorem \ref{theo:SSTM convergence rate for MV-sde} hold. 

Assume a further form of the Lipschitz condition of $b,\sigma$. Namely let $L_b,L_{\tilde{b} },L_\sigma,L_{\tilde{\sigma} }\ge0$ be such that
% \xbox{ Need to change}
	\begin{align*}
	 |b(t, x, \mu)-b(t, x', \mu')|^2&\leq L_b |x-x'|^2 + L_{\tilde{b} } W^{(2)}(\mu, \mu')^2,
	 \\
	 |\sigma(t, x, \mu)-\sigma(t, x', \mu')|^2&\leq L_\sigma |x-x'|^2 + L_{\tilde{\sigma} } W^{(2)}(\mu, \mu')^2.
	\end{align*}
for all $t\in[0,T]$, $x,x'\in\bR^d$ and $\mu,\mu'\in \cP_2(\bR^d)$. Suppose that $X_0\in L_0^m(\mathbb{R}^d)\nonumber$ and $Z_0\in L_0^m(\mathbb{R}^d)\nonumber$ for a sufficiently large $m\in \bN$, and let $\hx_0^i$ and  $\hz_0^i$ be i.i.d.~copies of $X_0$ and $Z_0$ respectively. 

Set $h>0$. Define two families $\{(X_n^{i,N},Y^{i,\star,N}_n)\}_{i,n}$ and $\{(Z_n^{i,N},G^{i,\star,N}_n)\}_{i,n}$ as the output of the SSM \eqref{eq:SSTM:scheme 1}-\eqref{eq:SSTM:scheme 2}
w.r.t correspond empirical measure pairs $\{\hm^{X,N}_n,\hm^{Y,N}_n\}_{n}$ and $\{\hm^{Z,N}_n,\hm^{G,N}_n\}_{n}$
with input initial conditions $\{X_0^{i,N}\}_i$ and $\{Z_0^{i,N}\}_i$ respectively.

Then, for any $n\in \bN$, 
\begin{align}
\label{eq:stable: beta h equation}
\sup_{ i\in \llbracket 1,N \rrbracket}  \bE\Big[\, |\hx_{n}^{i,N}-\hz_{n}^{i,N} |^2 \Big]
\le    
 (1+\beta h )^n \sup_{ i\in \llbracket 1,N \rrbracket}\bE\Big[ |\hx^{i,N}_{0}-\hz^{i,N}_{0} |^2\Big],
\end{align}
where 
% \xbox{Need to change}\cmark
\begin{align}
\label{SSTM:beta formula}
\beta=
\frac{(2L_v+A+1+L_{\tilde{v}} )+h(L_{\tilde{v}}A +B)+h^2 B L_{\tilde{v}}  }{1-  h (2L_v+1)}
\quad\textrm{and}\quad 
A=( 2\sqrt{L_b}+2\sqrt{L_{\tilde{b}}}+L_\sigma+L_{\tilde {\sigma}} )
\quad
B=L_b+L_{\tilde{b}}.
%\Bigg]
\end{align}
Under the choice of $h$ stated in Theorem \ref{theo:SSTM convergence rate for MV-sde}, the quantity $1+\beta h$ is always positive. 
If $L_v< -(1+L_{\tilde{v}}+A)/2\le -\frac{1}{2}$ and for a sufficient small $h$ then $\beta<0$ and thus the SSM is Mean-square contractive.
\end{Theorem}

The proof of this result is postponed to Section \ref{sec:ProofOfStability}. We illustrate immediately the scope of our findings with several numerical results. The reasoning behind the specification of the Lipschitz constants of $b$ and $\sigma$ will become apparent in the examples Section \ref{sec:contractivityCase}. The main motivation is to account for the contribution of the spatial and measure components separately as a way to explore the flexibility allowed by the scheme in choosing the drift coefficients $v$ and $b$.

The equivalent definition of \textit{Mean-square contractivity} (Definition \ref{def:definition of mean-square stable numerics X Y }) for the initial MV-SDE \eqref{Eq:General MVSDE} is the so-called \textit{Exponential mean-square stability inequality} defined next. We will make use of this definition in Section \ref{sec:contractivityCase} below.
\begin{definition}[Exponential mean-square contractive solutions]
\label{def:definition of generate exponential mean-square contractive solutions }
Let $X,Y$ be two solutions to \eqref{Eq:General MVSDE} with initial conditions $X_{0},Y_{0} \in L_{0}^{m}( \bR^{d})$ respectively. If $X,Y$ satisfy 
\begin{align*}
    \bE\Big[\,|X_t-Y_t|^2\Big] \le e^{\eta t} \bE\Big[ |X_0-Y_0|^2\Big]
\end{align*}
for some real number $\eta<0$, then the MV-SDE \eqref{Eq:General MVSDE} is said to generate \textit{exponential mean-square contractive solutions}.
\end{definition}

%%%%%%%%%%%%%%%%%%%%%%%%%%%%%%%%%%%%%%%%%%%%%%%%
%%%%%%%%%%%%%%%%%%%%%%%%%%%%%%%%%%%%%%%%%%%%%%%%%%%
%%%%%%%%%%%%%%%%%%%%%%%%%%%%%%%%%%%%%%%%%%%%%%%%%
% \newpage

%%%%%%%%%%%%%%%%%%%%%%%%%%%%%%%%%%%%%%%%%%%%%%%%
%%%%%%%%%%%%%%%%%%%%%%%%%%%%%%%%%%%%%%%%%%%%%%%%%%%
%%%%%%%%%%%%%%%%%%%%%%%%%%%%%%%%%%%%%%%%%%%%%%%%%

%%%%%%%%%%%%%%%%%%%%%%%%%%%%%%%%%%%%%%%%%%%%%%%%
%%%%%%%%%%%%%%%%%%%%%%%%%%%%%%%%%%%%%%%%%%%%%%%%%%%
%%%%%%%%%%%%%%%%%%%%%%%%%%%%%%%%%%%%%%%%%%%%%%%%%
% \newpage
\section{Examples of interest}
\label{sec:examples}

We now illustrate our numerical scheme \eqref{eq:SSTM:scheme 1},\eqref{eq:SSTM:scheme 2} through several examples of interest. Moreover, alongside the SSM simulations we also provide a comparative analysis with two other numerical schemes: Taming  \cite{reis2018simulation} and Adaptive timestepping \cite{reisinger2020adaptive}. For convenience, \ref{sec:ShortRecap} contains a brief description of the algorithms including convergence results and conditions.
 
Since the exact solution to each example is unknown we use a proxy for the true solution in order to compute approximation errors. Concretely, we compare SSM/Taming/Adaptive results with SSM/Taming/Adaptive results at a lower value of timestep $h$ as a proxy (each method is compared with its own approximation of the true solution but at a much smaller timestep). 
Within each example, all the methods will use same Brownian motion paths. 
\color{black}

We consider the weak error ($k=1$) and the strong error ($k=2$) between the true solution $X_T$ and the approximation $\hat X_T$ as follow
\begin{align*}
\epsilon_k
=
\begin{cases}
\Big| \bE\Big[\, X_T-\hat{X}_T\Big]\Big|_2
\approx \Big|\frac1N \sum_{j=1}^N X_T^j - \hat{X}_T^j|\Big|_2 
 , \qquad k=1,
 \\
\Big( \bE\Big[\, |X_T-\hat{X}_T|^k\Big]\Big)^{1/k}
\approx \Big(\frac1N \sum_{j=1}^N |X_T^j - \hat{X}_T^j|^k\Big)^{1/k}  , \qquad k=2.
\end{cases}
\end{align*}
Our main theorem covers only the strong convergence result, nonetheless, we also present the weak convergence rate estimation. 

We study several examples. The stochastic Ginzburg Landau example is a well-studied one and it provides a comparison example to other methods.
The second example is a multi-dimensional FitzHugh-Nagumo model of McKean-Vlasov type from neuroscience which shows that the split-step method can properly deal with the superlinear term in a complex system.  
% For these two examples we discuss the runtime of the serial and parallel implementation. 
For the first example, we discuss the parallel implementation, for the second example, we discuss the accuracy  w.r.t runtime.

The third example lies outside the scope of our assumptions by featuring a non-Lipschitz measure dependency, but all the methods still work. Proving the convergence of the method under this setting is left for future research.

In the last part, we discuss the stability of the SSM as understood by Theorem \ref{theo:SSTM:stabilty}. We first look at a linear case to compare the conditions for mean-square contractivity between MV-SDEs and the numerical scheme. Then, the non-linear Ginzburg Landau type equation is used to illustrate the mean square contractivity for the split-step method. At last, the Cucker-Smale flocking model (a degenerate MV-SDE) shows the split-step method has better properties compared to the other methods under larger choices of timestep $h$.

\begin{remark}[Parallel implementation]
\label{rem:parallelImplementation}
To implement the SSM \eqref{eq:SSTM:scheme 1}-\eqref{eq:SSTM:scheme 2} in parallel, at each timestep, we first solve step \eqref{eq:SSTM:scheme 1} of the SSM distributed by the cores, then calculate the empirical measure of the particle system in the first core, and finally the second step \eqref{eq:SSTM:scheme 2} of the SSM is executed also in parallel. 
To implement Taming \eqref{Eq:Tamed MVSDE} and Adaptive \eqref{Eq:Adaptive MVSDE} in parallel, at each timestep $h$ one needs to first communicate the empirical measure of the particle system between different cores and then each core can calculate its particles dynamics independently.

A priori one should expect the taming to be the fastest of the three algorithms. We highlight that due to the nature of the empirical measure, after each timestep the processors need to communicate to update the calculation of the empirical measure, this leads to a well-known loss of parallelization power \cite{bernal2017}. 
\end{remark}

\subsection{Example: the stochastic Ginzburg Landau equation}
\label{exam:StochGinzburg}
The first example we consider is the mean-field perturbation of the classic one-dimensional stochastic Ginzburg Landau equation, namely,  for all $t\in[0,T]$, 
\begin{align}
\label{eq:MV-GL-SDE}
dX_t&=\Big(\frac{(\sigma')^2}{2}X_t-X_t^3+c\bE[X_t]\Big)dt+\sigma' X_tdW_t,
\quad X_0=x_0 \in\mathbb{R}.
\\
\label{eq:MV-GL-SDE-veresion1}
\textrm{where in SSM version 1:}\quad v(t,x,\mu)&=-x^3,
\quad b(t,x,\mu)=\frac{(\sigma')^2}{2} x+c\int_\bR x \mu(dx) , \quad \sigma(t,x,\mu)= \sigma' x,
\\
\label{eq:MV-GL-SDE-veresion2}
\textrm{in SSM version 2:}\quad v(t,x,\mu)&=-x^3+c\int_\bR x \mu(dx),
\quad b(t,x,\mu)=\frac{(\sigma')^2}{2} x , \quad \sigma(t,x,\mu)= \sigma' x.
\end{align}
where $\sigma',c$ are constants. 

\begin{figure}[h!bt]
    \centering
        \begin{subfigure}{.33\textwidth}
			\centering
 			\includegraphics[scale=0.33]{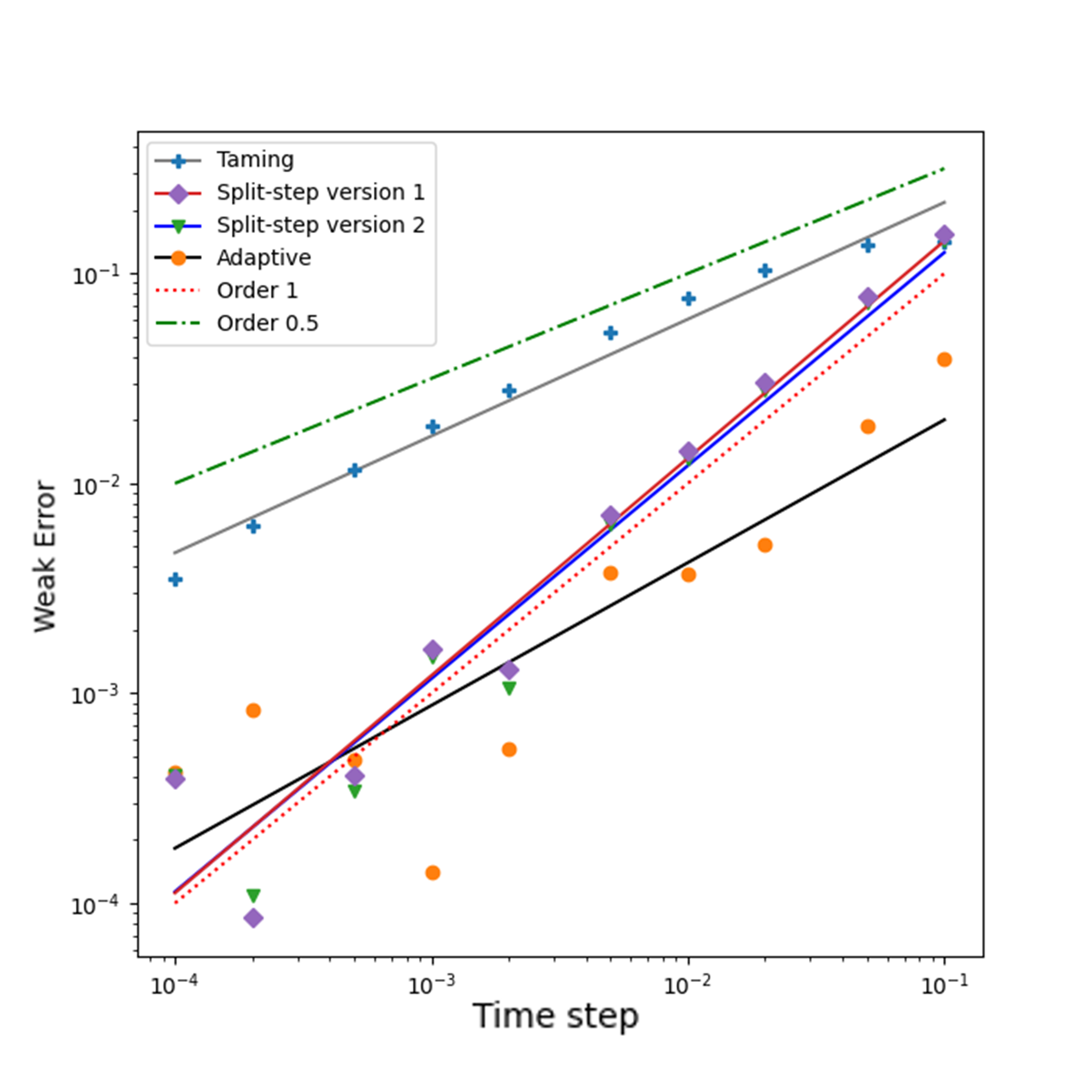}
			\caption{Weak Error w.r.t Time step}
			\label{fig:11}
		\end{subfigure}%
		\begin{subfigure}{.32\textwidth}
			\centering
 			\includegraphics[scale=0.33]{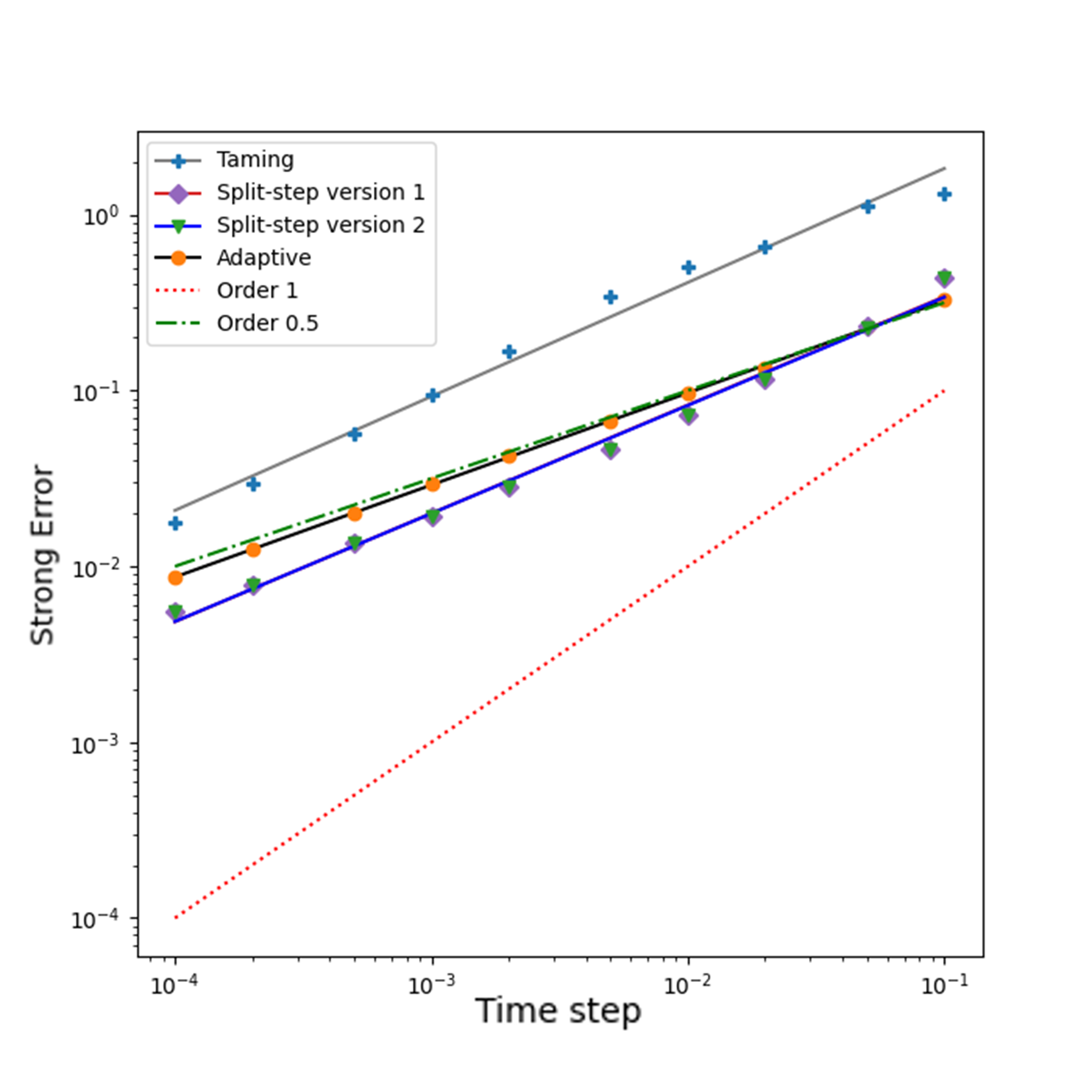} 
			\caption{Strong Error (rMSE) w.r.t Time step}
			\label{fig:12}
		\end{subfigure}
		\begin{subfigure}{.32\textwidth}
			\centering
 			\includegraphics[scale=0.33]{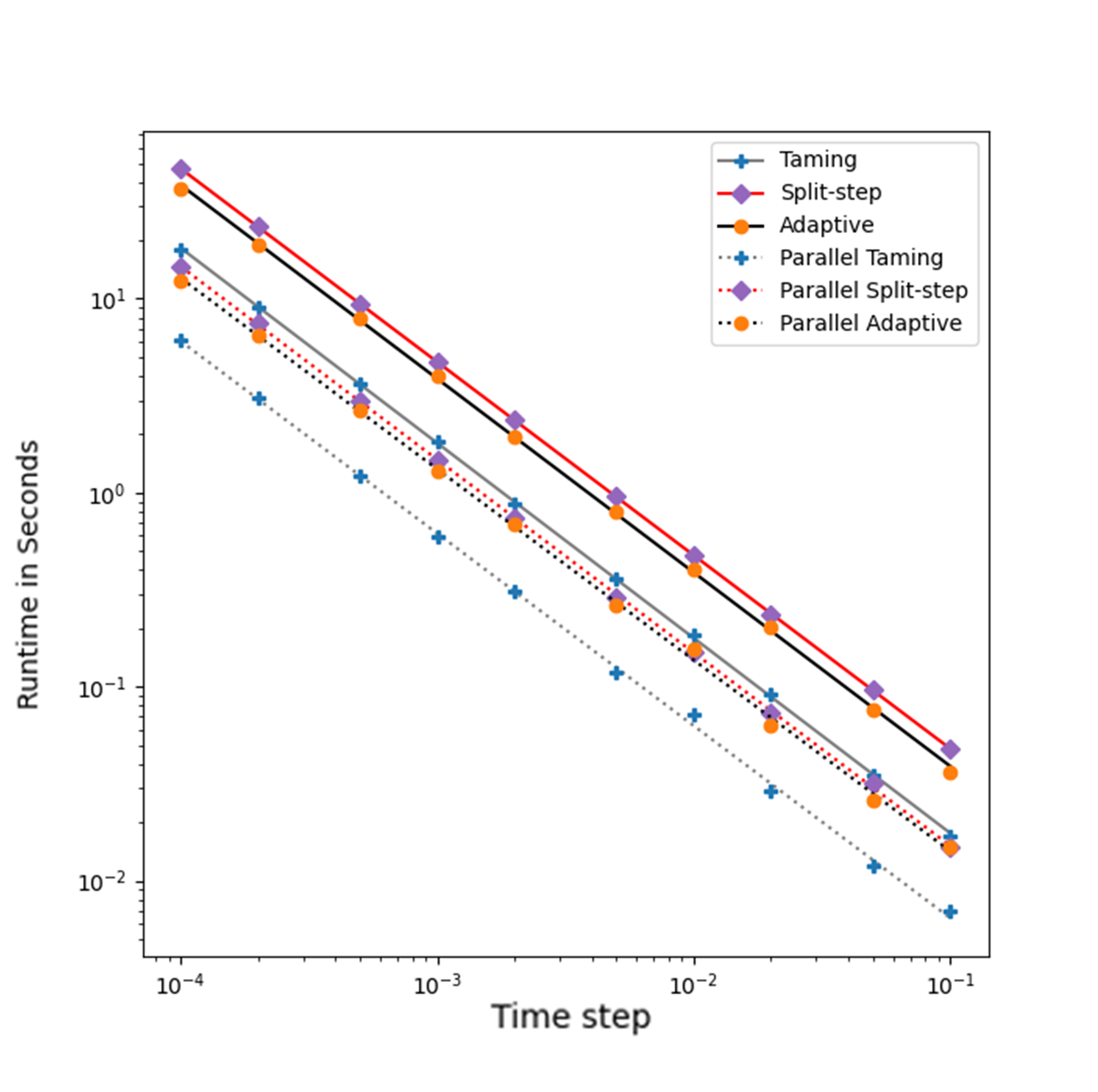} 
			\caption{Algorithm runtime w.r.t Time step}
			\label{fig:13}
		\end{subfigure}
    \caption{ Simulations of the stochastic Ginzburg Landau equation with $N=1000$ particles. (a) Weak error of different methods. (b) Strong error of different methods. (c) Runtime of different methods in serial and in parallel. }
    \label{fig:gl example}
\end{figure}

This equation is a toy one-dimensional MV-SDE that we use for a methodological comparison analysis. It features in \cite[Section 4.1]{reis2018simulation} and \cite[Section 4.2]{reisinger2020adaptive}, thus a comparison is insightful under the same choice of coefficients. 
Namely, we take $\sigma'=1.5$, $x_0=1$, $c=0.5$ , $T=1$ and $N=1000$ particles. The results are shown in  Figure \ref{fig:gl example}. The timestep are $h\in \{ 10^{-4},2\times10^{-4}, 5\times10^{-4}, 10^{-3},\dots, 10^{-1} \}$. The true solution is calculated with $h=10^{-5}$ (for each scheme). We have tested both versions of SSM, \eqref{eq:MV-GL-SDE-veresion1} and \eqref{eq:MV-GL-SDE-veresion2}, and both show similar results (we present only one). 

Taming is implemented with $\alpha=0.5$ while Adaptive under the choice  $\mathbf{h}^{\delta}(x)=h\min(1,|x|^{-2})$. Fig \ref{fig:gl example}(a) shows the Weak error rate of the Taming to roughly be $1/2$ where for other methods it is 1.0. 
Fig \ref{fig:gl example}(b) shows the rMSE rate of all the methods to be around $1/2$ with Taming's error being about one-order of magnitude higher than the other two errors (also observed in \cite[Section 4]{reisinger2020adaptive}). 
The two versions of the SSM have similar behaviour for this model. 
% We show the running time for version 1 in the next figure. 
Fig \ref{fig:gl example}(c) depicts running times of 3 methods (version 1 of the SSM), the top 3 lines are the standard implementations (non-parallel) and the bottom 3 lines are the parallel implementation with 4 cores. 

In both the parallel and non-parallel implementation, Taming is the fastest while the SSM takes a slightly longer time than the other methods but with a performance comparable to Adaptive. In a parallel implementation with $4$ cores we reach a reduction to nearly $27\%$ in relation to the non-parallel implementation's computational time. 
{ 
In this example, to reach the same strong error level Taming takes nearly $7$-times more time than SSM; Adaptive is similar to SSM.}

\subsection{Example: the FitzHugh-Nagumo model}
\label{sec:ExampleFHN}

This is a three-dimensional ($d=3$) MV-SDE \eqref{Eq:General MVSDE} defined with  $v:[0,T]\times\mathbb{R}^3\rightarrow \mathbb{R}^3$, $b:[0,T]\times\mathbb{R}^3\times \cP_2(\mathbb{R}^3)\rightarrow \mathbb{R}^3$, $\sigma:[0,T]\times\mathbb{R}^3\times \cP_2(\mathbb{R}^3)\rightarrow \mathbb{R}^{3\times 3}$ for $x=(x_1,x_2,x_3)\in\mathbb{R}^3$, $z\in \mathbb{R}$ and $\mu^{x_3}$ is the marginal measure in $x_3$, 
\begin{align*}
 v(t,x)=
    \begin{pmatrix}
-(x_1)^3/3
\\
0
\\
0
\end{pmatrix}
,\quad
b(t,x,\mu)=
    \begin{pmatrix}
x_1-x_2+I-\int_{\mathbb{R}^3} J(x_1-V_{rev})z d\mu^{x_3}(z)
\\
c(x_1+a-b x_2)
\\
a_r \frac{T_{max}(1-x_3)}{1+\exp(-\lambda(x_1-V_T))}-a_dx_3
\end{pmatrix},
\end{align*}
\begin{align*}
    \sigma(t,x,\mu)=
    \begin{pmatrix}
\sigma_{ext}& 0 & -\int_{\mathbb{R}^3} \sigma _J(x_1-V_{rev})z d\mu^{x_3}(z)
\\
0 & 0& 0
\\
0 & \sigma_{32}(x) &0
\end{pmatrix},\quad
x_0\sim \mathcal{N}
    \begin{pmatrix}
\begin{pmatrix}
V_0
\\
w_0
\\
y_0
\end{pmatrix}
,
\begin{pmatrix}
&\sigma_{V_0} & 0& 0
\\
& 0 &\sigma_{w_0}& 0
\\
& 0& 0 &\sigma_{y_0}
\end{pmatrix}
\end{pmatrix},
\end{align*}
where
\begin{align*}
    \sigma_{32}(x):=\mathds{1}_{\{x_3\in(0,1)\}}\sqrt{a_r\frac{T_{max}(1-x_3)}{1+\exp(-\lambda(x_1-V_T) )}+a_d x_3} \times \Gamma \exp\Bigg[ -\frac{\Lambda}{1-(2x_3-1)^2} \Bigg].
\end{align*}
and $I,J,V_{rev},V_{T},T_{\max},\Gamma,\Lambda,\sigma_{ext},\sigma_{J},\sigma_{32},a,b,c,a_r,a_d,\lambda$ are constants.

All the parameters are the same as in \cite[Section 4.3]{reis2018simulation} ($V_0=0,w_0=$1/2$,\dots$) -- see also \cite[Section 4.4]{reisinger2020adaptive}. The split-structure of the SSM enable us to flexibly deal with the superlinear term and the Lipschitz terms separately -- we make judicious choices regarding $v$ and $b$ and thus we can optimise the solver of the implicit part. Otherwise, we would have been forced to use a general purpose solver which is costlier.
\begin{figure}[h!bt]
    \centering
        \begin{subfigure}{.33\textwidth}
			\centering
 			\includegraphics[scale=0.33]{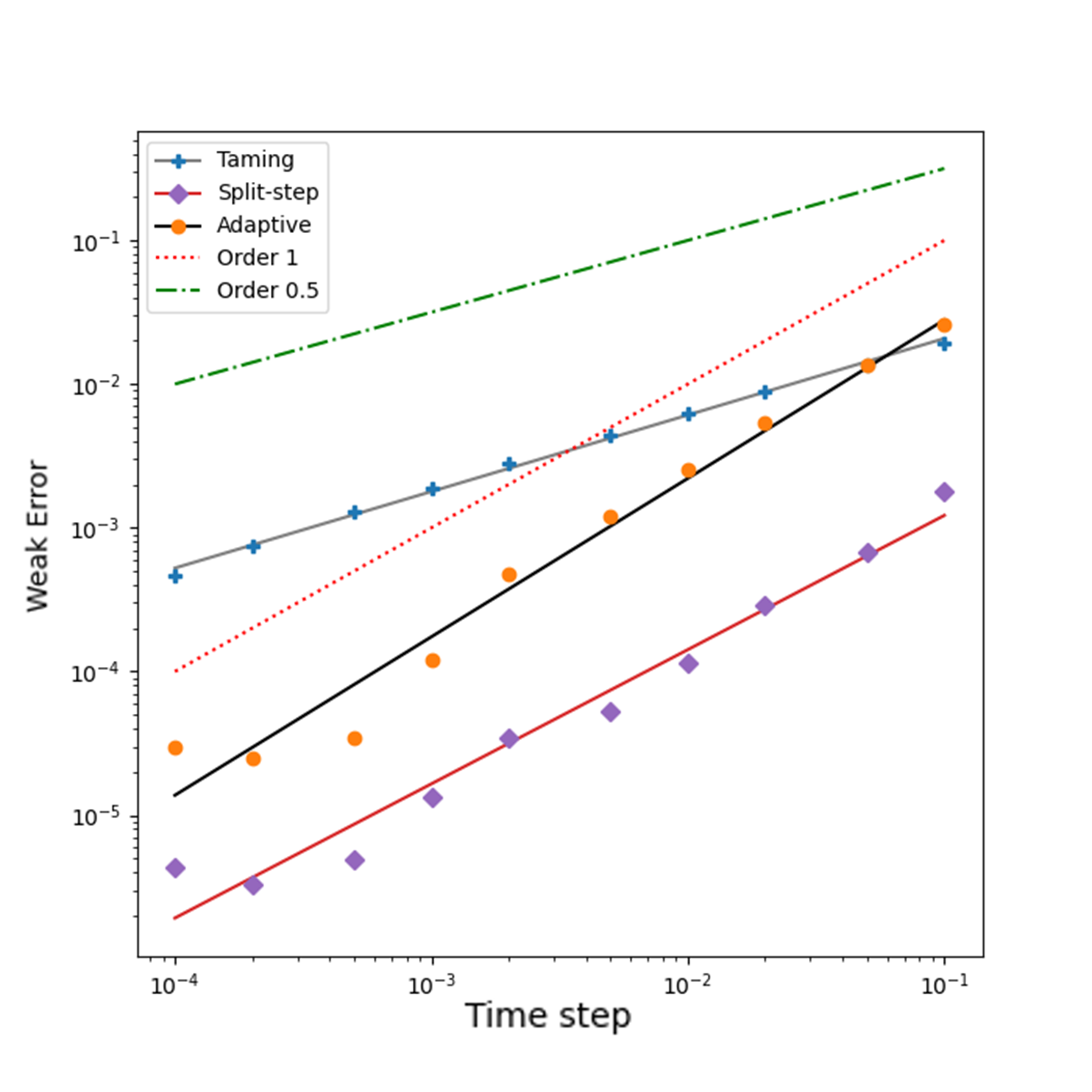}
			\caption{Weak Error w.r.t Time step}
			\label{fig:21}
		\end{subfigure}%
		\begin{subfigure}{.33\textwidth}
			\centering
 			\includegraphics[scale=0.33]{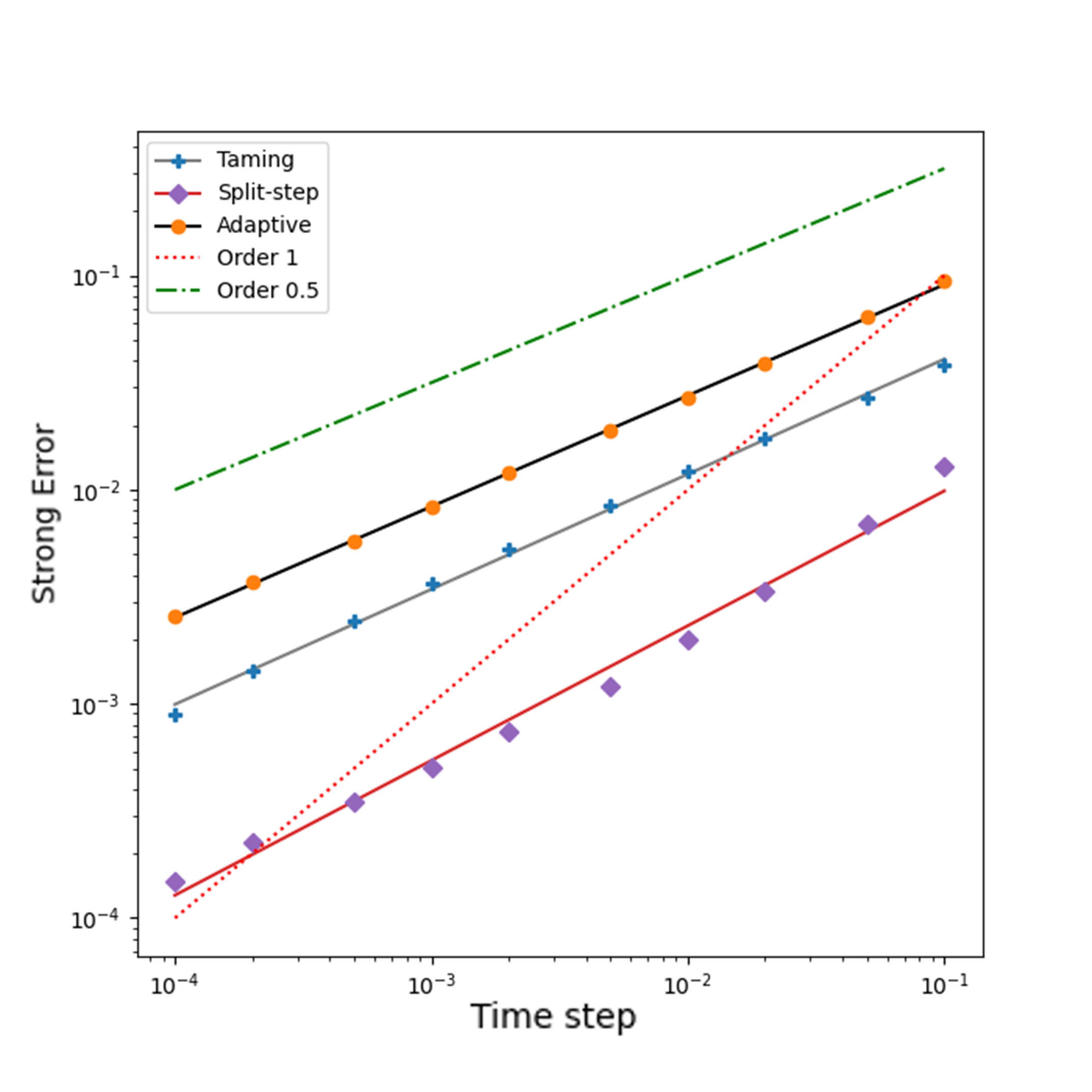} 
			\caption{Strong Error (rMSE) w.r.t Time step}
			\label{fig:22}
		\end{subfigure}
        \begin{subfigure}{.33\textwidth}
			\centering
 			\includegraphics[scale=0.33]{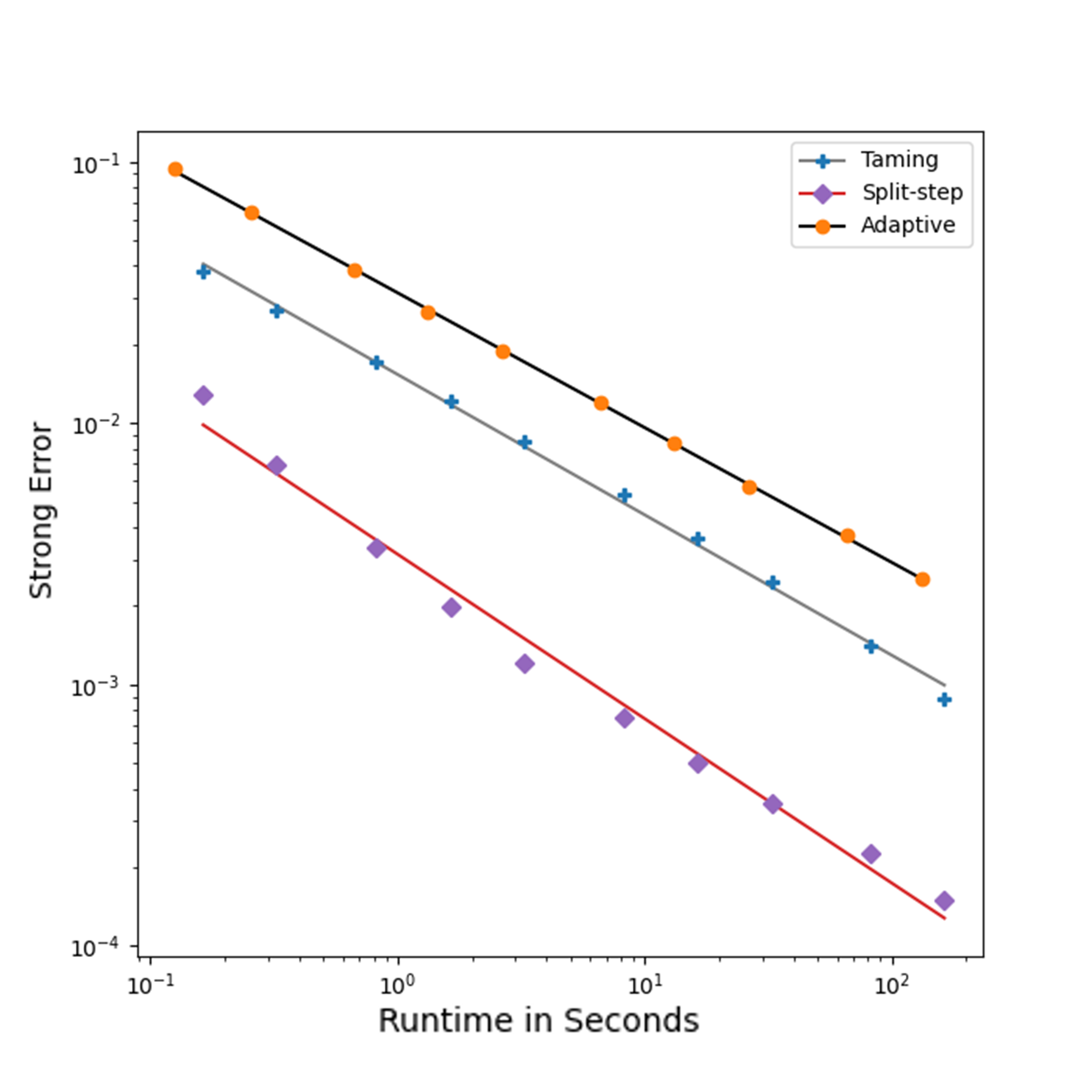}
			\caption{Strong Error (rMSE) w.r.t Algorithm runtime }
			\label{fig:23}
		\end{subfigure}%
        \caption{Simulations of the FitzHugh-Nagumo model with $N=1000$ particles.  (a) Weak error of different methods. (b) Strong error of different methods. (c) Strong error compare to different Algorithm runtime.}
    \label{fig:Simulation_Nn}
\end{figure}

For the simulation, we take $N=1000$, $T=2$, the time step is taken from $h\in \{ 10^{-4},2\times10^{-4}, 5\times10^{-4}, 10^{-3},\dots, 10^{-1} \}$ and the true solution is calculated with $h=10^{-5}$. Taming is implemented with $\alpha=1/2$ and Adaptive with $\mathbf{h}^{\delta}(x)=h\min(1,|x|^{2}/|b(t,x,\mu)|^2)$.  
Fig \ref{fig:Simulation_Nn}(a) shows the weak error rate of Taming to roughly be $1/2$ with other methods being 1.0 (implementing Taming with $\alpha=1$ yields a weak convergence rate of order $1.0$, we do not present the result). 
Fig \ref{fig:Simulation_Nn}(b) shows the strong error rate of all the methods to roughly be $1/2$, the error of SSM is an order of magnitude lower than the others. 
 {   Fig \ref{fig:Simulation_Nn}(c) shows that, to reach the same strong error level Taming takes nearly $80$-times more time than SSM; Adaptive takes nearly $10$-times more time than SSM.
 }

\subsection{Example: Polynomial drift (non-Lipschitz measure dependency but still of one-sided Lipschitz type)}

We present an example from \cite[{Section 3.3}]{bencheikh2019bias} that falls outside the theoretic framework of this work. Take the one-dimensional MV-SDE for $t\in[0,T]$ and $\gamma\in\bR$
\begin{align}
\label{eq:poly1}
    d X_{t}&=\Big(\gamma X_{t}+\mathbb{E}\left[X_{t}\right]-X_{t} \mathbb{E}\left[\,|X_{t}|^{2}\right]\Big) d t+X_{t} d W_{t}, \quad \text{with } X(0)=x_{0}\in \bR,
    \\
    \nonumber
    \textrm{and set:}\quad v(t,x)&=\gamma x-x\int_\bR |x|^2 \mu(dx),
\quad b(t,x,\mu)=\int_\bR x \mu(dx) , \quad \sigma(t,x,\mu)= x.
\end{align}
The dynamics of the interacting particle system \eqref{Eq:MV-SDE Propagation}, for ${ i\in \llbracket 1,N \rrbracket}$, $X^{i,N}\in \mathbb{R} $, is 
\[
d X_{t}^{i, N}=\Big[\gamma X_{t}^{i, N}+\frac{1}{N} \sum_{j=1}^{N} X_{t}^{j, N}-X_{t}^{i, N} \frac{1}{N} \sum_{j=1}^{N}|X_{t}^{j, N}|^{2}\Big] d t+X_{t}^{i, N} d W_{t}^{i}.
\]
where $(W^{i})_{i}$ are independent Brownian motions.
Take $\gamma=-1,~T=1$. We have a superlinear measure component in \eqref{eq:poly1} and none of the schemes is applicable (insofar as existing theoretical results allow). If the measure component was fixed (with finite 2nd moment), then, the drift would satisfy a one-sided Lipschitz condition -- the convergence of this type of schemes is left for future work.
\begin{figure}[h!bt]
    \centering
        \begin{subfigure}{.49\textwidth}
			\centering
 			\includegraphics[scale=0.48]{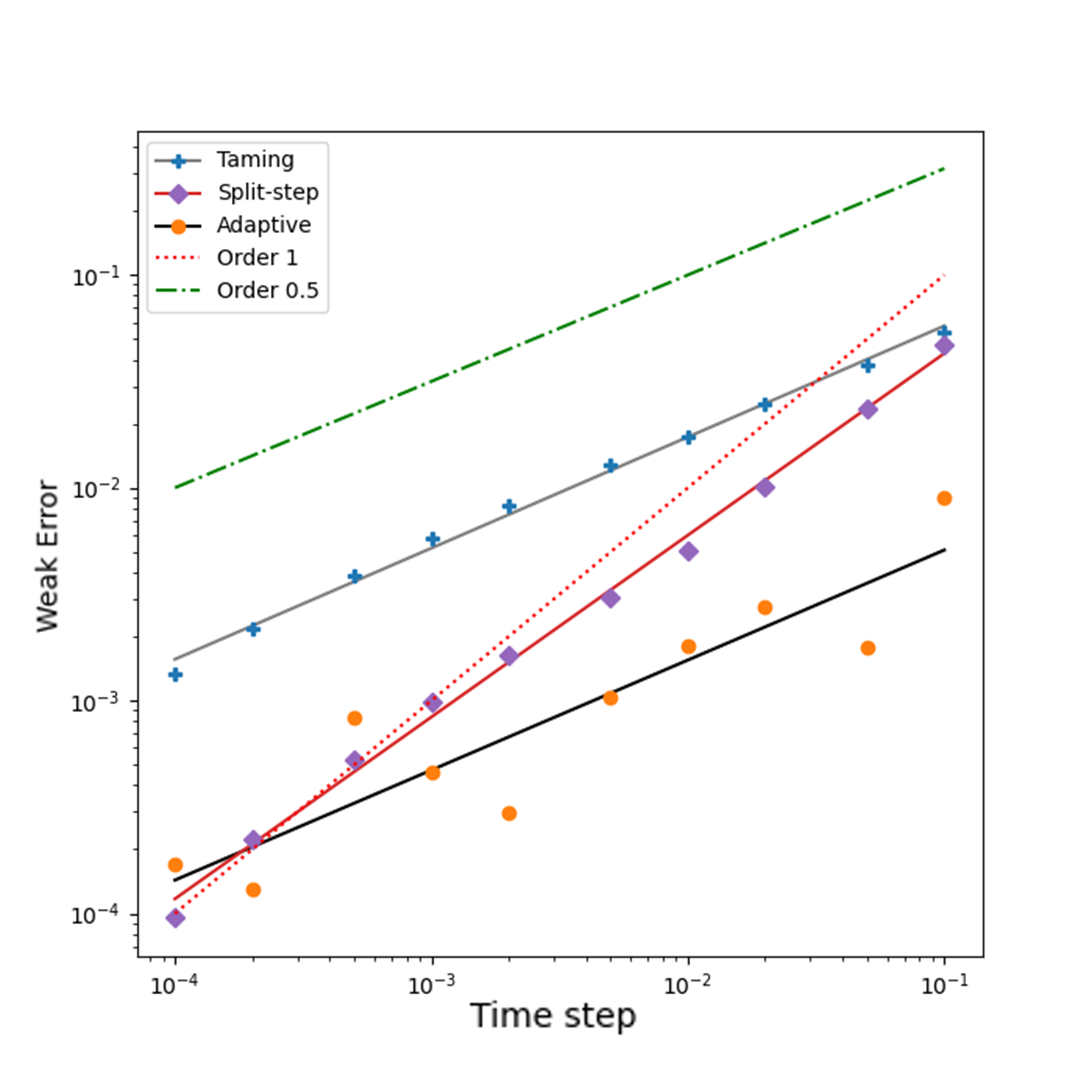}
			\caption{Weak Error w.r.t Time step}
			\label{fig:poly21}
		\end{subfigure}%
		\begin{subfigure}{.49\textwidth}
			\centering
 			\includegraphics[scale=0.48]{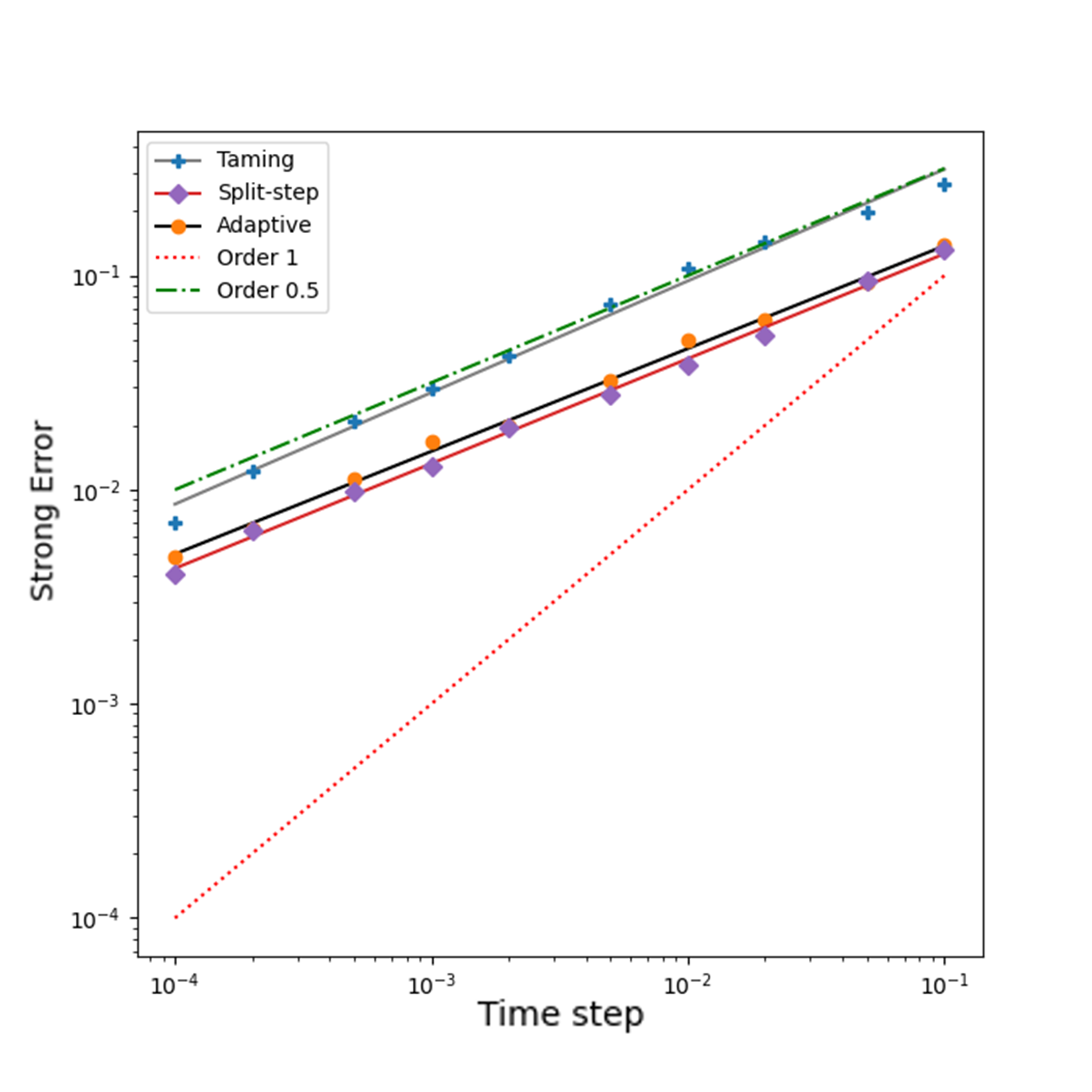}
			\caption{Strong Error (rMSE) w.r.t Time step}
			\label{fig:poly22}
		\end{subfigure}
    \caption{Simulations of the polynomial drift model with $N=1000$ particles. (a) Weak error of different methods. (b) Strong error of different methods.}
    \label{fig:Simulation_poly}
\end{figure}

The results are shown in  Figure   \ref{fig:Simulation_poly}. We take $N=1000,~T=2$ and the timestep is taken from $h\in \{ 10^{-4},2\times10^{-4}, 5\times10^{-4}, 10^{-3},\dots, 10^{-1} \}$. The true solution is calculated under $h=10^{-6}$. Taming is applied with $\alpha=0.5$ while Adaptive  under the choice  $\mathbf{h}^{\delta}(x)=h\min(1,|x|^{-2})$.  
Figure \ref{fig:Simulation_poly}(a) shows the weak error rate of Taming with $\alpha=0.5$ to roughly be $1/2$ with other methods being 1.0. 
Figure \ref{fig:Simulation_poly}(b) shows strong error rate of all the methods to roughly be $1/2$ (as expected).

\subsection{Stability of the SSM: linear, non-linear and the Cucker-Smale model case}
\label{sec:contractivityCase}

Recall from Theorem \ref{theo:SSTM:stabilty} the expression \eqref{SSTM:beta formula} for $\beta$. For the remainder of this section, let $t \in[0,T]$, we define $X_0,Z_0\in L_0^m(\mathbb{R}^d)\nonumber$ and $\hx_0^i,\hz_0^i,~i\in \llbracket 1,N\rrbracket$ as i.i.d.~samples from $X_0$ and $Z_0$ respectively, $X_n^{i,N}$ and $Z_n^{i,N}$ are defined as in Theorem \ref{theo:SSTM:stabilty} as outputs of our SSM scheme \eqref{eq:SSTM:scheme 1} and \eqref{eq:SSTM:scheme 2} (with the corresponding initial conditions). If $\beta<0$ we then have  $\bE[|X_n^{i,N}|^2]=0$ as $n\rightarrow \infty$.

\subsubsection[Linear case: an Ornstein-Uhlenbeck McKean-Vlasov SDE]{Linear case: an   Ornstein-Uhlenbeck McKean-Vlasov SDE} 
\label{section:OU example}
For the MV-SDE (see e.g. \cite[{Section 2.1}]{bencheikh2019bias}), for all $t \in[0,T]$ and $x_0\in\bR$ 
\begin{align}\label{eq:stablity:OU}
d X_t&=\Big(\rho  X_t+\lambda \bE[ X_t] \Big) d t+\eta ~dW_t, \quad X_0=x_{0},
\\
\textrm{set:}\quad 
v(t,x,\mu)&=\rho x,
\quad b(t,x,\mu)=\lambda \int_\bR x \mu(dx) , \quad \sigma(t,x,\mu)= \eta. 
\end{align}
where $\rho,\lambda,\eta $ are constants. The first and second moments of $X$ are respectively given by  $\mathbb{E}\left[X_{t}\right]=x_{0} \exp ((\rho+\lambda) t)$ and $\mathbb{E}\left[X_{t}^{2}\right]=
x_{0}^{2} \exp (2(\rho+\lambda) t)+\frac{\eta^{2}}{2 \rho}(\exp (2 \rho t)-1)$.  

Let $X, Z$ be two solution of \eqref{eq:stablity:OU} with $X_0$ and $Z_0$ as initial condition respectively, then by direct calculation 
\begin{align*}
% \bE\Big[X_t-Z_t \Big]&=\bE\Big[ X_0-Z_0 \Big]e^{(\rho+\lambda) t }
% \\
\bE\Big[|X_t-Z_t|^2 \Big]&=\frac{1}{2\lambda}e^{2(\rho+\lambda) t} +  \bE\Big[ |X_0-Z_0|^2 \Big]e^{2\rho t }. 
\end{align*}
Let $\rho\le 0$ and $\rho+\lambda<0$  then from Definition \ref{def:definition of generate exponential mean-square contractive solutions }, \eqref{eq:stablity:OU} generates exponential mean-square contractive solutions.
The parameters of this example are $L_v=\rho,~L_{\tilde{b}}=\lambda^2$, $L_{\tilde{v}}=L_b=L_{\sigma}=L_{\tilde{\sigma}}=0$. Plugging  these into \eqref{SSTM:beta formula} and in order to make $\beta<0$, we need to choose $h$ satisfying
\begin{align*}
    h<  -\frac{2\rho+2\lambda+1 }{\lambda^2} .
\end{align*}
From Definition \ref{def:definition of mean-square stable numerics X Y }, the split-step method \eqref{eq:SSTM:scheme 1}-\eqref{eq:SSTM:scheme 2} is mean-square contractive.
So, for the split-step method,  $h$ exists when $\rho+\lambda<-1/2$. Both this condition and  the condition for the SDE to generate exponential mean-square contractive solutions need the constraint $\rho+\lambda<0$. Thus, the condition for a mean-square stable numerical solution is slightly stronger than the condition for the SDE to generate exponential mean-square contractive solutions.

\subsubsection{Nonlinear case I: a stochastic Ginzburg Landau type equation}
\label{sec:stabilityginzburgLandau}

We illustrate the stability of the SSM scheme via the stochastic Ginzburg Landau type equation (in the style of that in Section \ref{exam:StochGinzburg}), we consider the following one-dimensional MV-SDE for all $t \in[0,T]$ 
\begin{align}
\label{eq:satblity:MV-GL-SDE}
dX_t&=\Big(-\frac{5}{2}X_t-\frac{1}{4}X_t^3+\bE[X_t]\Big)dt+ X_tdW_t,\qquad X_0=1,
\\ \nonumber 
\textrm{set:}\quad v(t,x,\mu)&=-\frac{5}{2}x-\frac{1}{4}x^3-\gamma x,
\quad b(t,x,\mu)=\int_\bR x \mu(dx) +\gamma x \ \textrm{ for }\ \gamma\in\bR, \quad \sigma(t,x,\mu)= x.
\end{align}
The parameters of this example are $L_v=-{5}/{2}-\gamma$, $L_{\tilde{b}}=L_{\sigma}=1$, and $L_b=\gamma^2,~L_{\tilde{v}}=L_{\tilde{\sigma}}=0$  -- with these parameters it is known \cite{zhang2021existence} that the system is conservative and the solution satisfies $X_t\to 0$ a.s.~as $t \to \infty$.   
Plugging these into the mean-square stability $\beta$ constant \eqref{SSTM:beta formula} and,  when $\gamma=0$, one must have small $h$ in order to make $\beta<0$. We now employ the split-step method, under different choices of $h$ and initial values. Set the number of particles to be $N=1000$.

\begin{figure}[h!bt]
    \centering
        \begin{subfigure}{.45\textwidth}
			\centering
 			\includegraphics[scale=0.33]{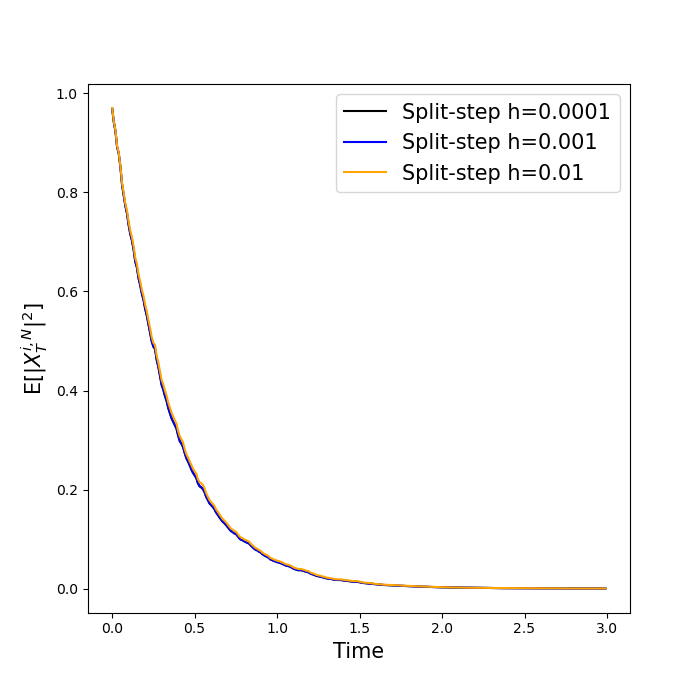}
			\caption{$(h,t) \mapsto \bE[ |X_t^{i,N}|^2]$}
			\label{fig:stab1}
		\end{subfigure}%
		\begin{subfigure}{.45\textwidth}
			\centering
 			\includegraphics[scale=0.33]{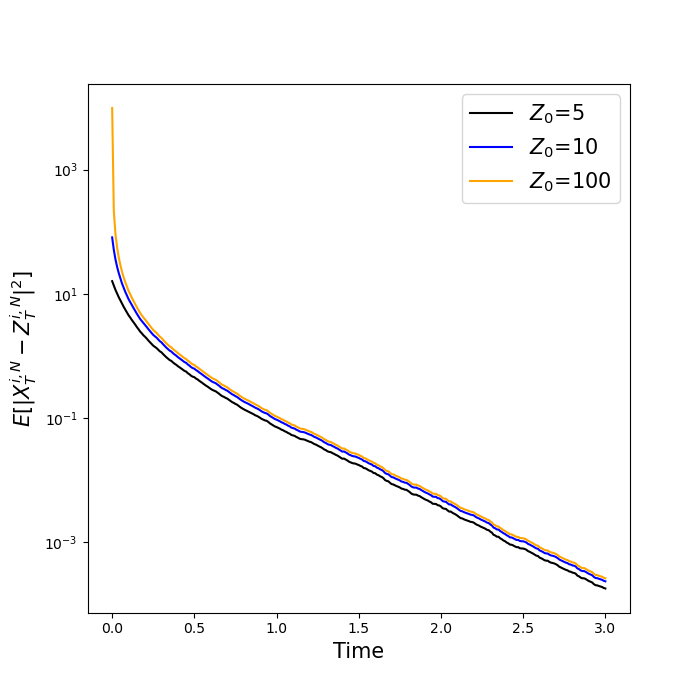}
			\caption{$(Z_0,t)\mapsto \bE[ |X_t^{i,N}-Z_t^{i,N}|^2]$}
			\label{fig:stab2}
		\end{subfigure}
		\\
		\begin{subfigure}{.45\textwidth}
			\centering
 			\includegraphics[scale=0.33]{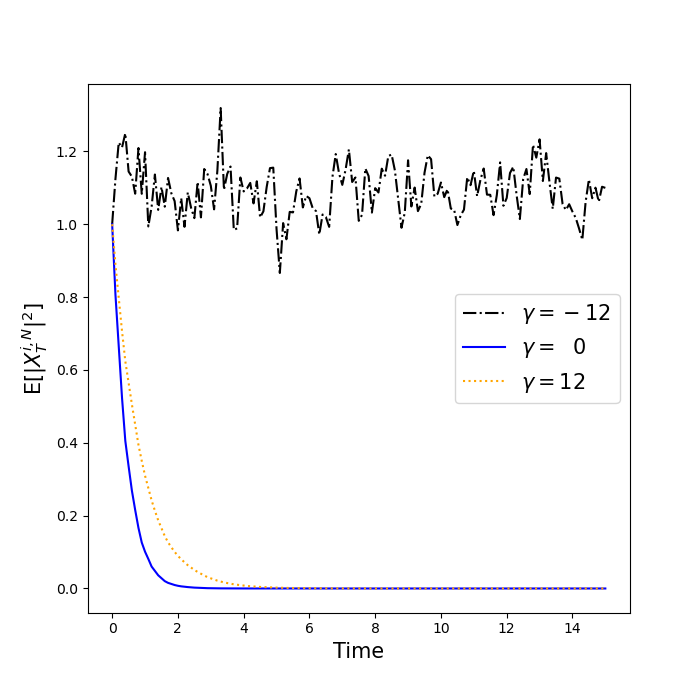}
			\caption{$(\gamma,t) \mapsto \bE[ |X_t^{i,N}|^2]$}
			\label{fig:stab3}
		\end{subfigure}
		\begin{subfigure}{.45\textwidth}
			\centering
			\includegraphics[scale=0.26]{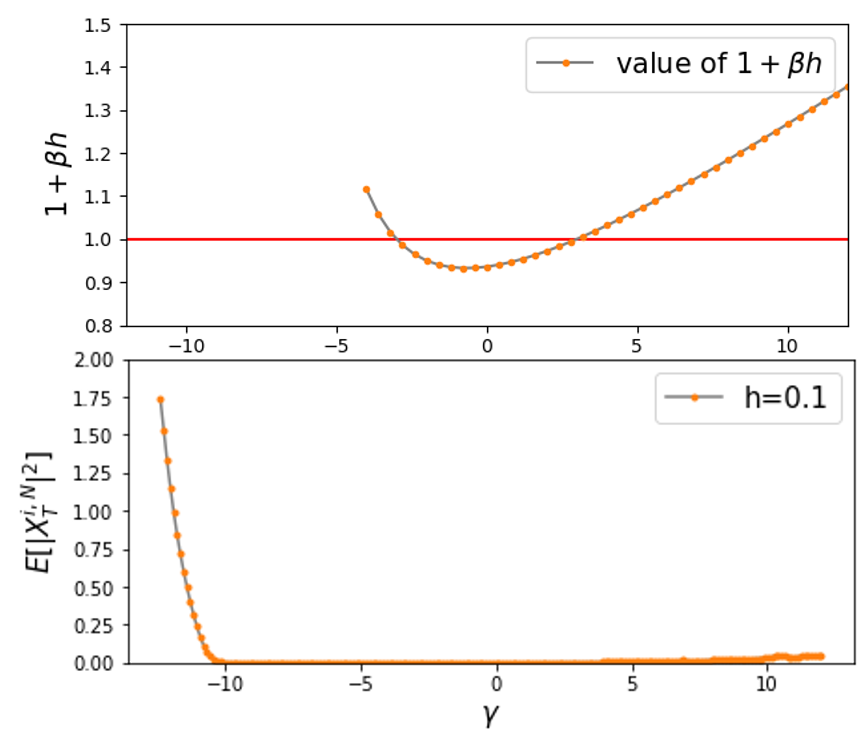}
			\caption{top: $\gamma \mapsto 1+\beta(\gamma) h$, bottom: $\gamma \mapsto \bE[\,|X_T^{i,N}|^2]$  }
			\label{fig:stab4}
		\end{subfigure}
    \caption{Simulations of the stochastic Ginzburg Landau type equation \eqref{eq:satblity:MV-GL-SDE} with $N=1000$ particles and $T=3.0$. 
    (a) shows $[0,3]\ni t\mapsto \bE[ |X_t^{i,N}|^2]$  under three different timesteps $h\in\{10^{-2},10^{-3},10^{-4}\}$ highlighting  mean-square stability. 
    (b) shows the mean square differences between $X,Z$ when $h=0.01$ for fixed $X_0=1$,  $Z_0\in\{5,10,100\}$ and $t\in[0,T]$ highlighting mean-square contractivity.
    \\
    Highlighting mean-square stability/instability of approximation as map of $\gamma$ under fixed $h=0.1$.  
    (c) this shows $[0,T]\ni t\mapsto \bE[ |X_t^{i,N}|^2]$ (for $T=15$) under three different $\gamma \in\{-12,0,12\}$. (d) (top) shows $\gamma \mapsto 1+\beta(\gamma) h$ where $\beta$ is given in Theorem \ref{theo:SSTM:stabilty} and (bottom) $[-12,12]\ni \gamma \mapsto \bE[ |X_T^{i,N}|^2],~T=3$. As $\gamma \geq 5$ the method starts showing an error increase (bottom) which can be matched to $\beta(\gamma)>0$ (top) and hence loss of stability.
    }
    \label{fig:Simulation_Stability}
\end{figure}

Set $\gamma=0$, $T=3$ and $X_0=1$. Figure \ref{fig:Simulation_Stability}(a) shows that $\bE[ |X_T^{i,N}|^2]$ decreases to zero under different values of $h$ (but small). Figure \ref{fig:Simulation_Stability}(b) shows mean square contraction property between $X$ and $Z$ highlights an exponential decay (where $Z$ solves \eqref{eq:satblity:MV-GL-SDE} for $Z_0=5,10,100$.)

Fix $h=0.1$. Figure \ref{fig:Simulation_Stability}(c) shows that when $\gamma=-12$ the scheme performs poorly, and this follows from the conditions of Theorem \ref{theo:SSTM:stabilty} not being satisfied. For  $\gamma\in\{0,12\}$, the scheme shows contraction as $t\to \infty$, but it is much slower for  $\gamma=12$.  
Figure \ref{fig:Simulation_Stability}(d, lower graph)  shows what happens when one shifts ``slope from $v$ to $b$'' via the linear term $\gamma x$ (see Remark \ref{remark:removing constraint on Lv}). 
We have now $L_v=-\gamma-5/2$, thus, when $\gamma\in (-5,5)$ the figure shows \textit{contraction}, with $X_T^{i,N}\approx 0$ as expected. There is a  significant change for $\gamma<-10$ where the approximation is not converging to the correct value. For $\gamma\geq 5$ and higher (recall that $h=0.1$ is fixed) it seems the contraction is happening (although at a  slower pace) but in Figure \ref{fig:Simulation_Stability}(d, upper graph) one sees that $\gamma \to 1+\beta(\gamma)h$ is now above $1.0$ which does not guarantee contraction (in the sense of Theorem \ref{theo:SSTM:stabilty}). 

Figure \ref{fig:Simulation_Stability}(c) and (d) highlight the trade off and care needed between: \textit{(i)} making $L_v$ negative via $\gamma$ and thus removing the constraint on $h$ imposed in Theorem \ref{theo:SSTM convergence rate for MV-sde}, and, \textit{(ii)} ensuring the stability of the scheme as imposed by Theorem \ref{theo:SSTM:stabilty}.

\subsubsection{Nonlinear case II: the two-dimensional Cucker-Smale flocking model}

This example (see \cite[{Section 2}]{csreference}) highlights the stability of the split-step method. It is stable under larger timestep $h$ by using the implicit step for the superlinear part. The explicit methods (Taming and Adaptive) fail to have acceptable results at this level of $h$.

Applied our settings, this is a two-dimensional MV-SDE define under $v,b,\sigma:[0,T]\times\mathbb{R}^2\times \cP_2(\mathbb{R}^2)\rightarrow \mathbb{R}^2$ for $x=(V,X)\in\mathbb{R}^2~,z\in\mathbb{R}$, $\mu^V$ is the measure of $V$ as:
\begin{align*}
 v(t,x,\mu)=
    \begin{pmatrix}
-(V)^3
\\
0
\end{pmatrix}
,\ \ 
b(t,x,\mu)=
\begin{pmatrix}
1+\lambda \int_{\mathbb{R}} (V-z) d\mu^V(z)
\\
V
\end{pmatrix},
\ \textrm{ and }\ 
    \sigma(t,x,\mu)=
    \begin{pmatrix}
\sigma' \int_{\mathbb{R}} (V-z) d\mu^V(z)
\\
0 
\end{pmatrix}.
\end{align*}
where $\lambda,\sigma'$ are constants.
The dynamics of the particle system follows easily 
\begin{align*}
d V_{t}^{i, N}&=\Big[1-(V_{t}^{i, N})^3+ \frac{\lambda}{N} \sum_{j=1}^{N} (V_{t}^{j, N}-V_{t}^{i, N}) \Big] d t+\frac{\sigma'}{N} \sum_{j=1}^{N} (V_{t}^{j, N}-V_{t}^{i, N})d W_{t}^{i},
% \\
\quad d X_{t}^{i, N}= V_{t}^{i, N} dt.
\end{align*}
where $i\in \llbracket 1,N \rrbracket,V^{i, N},X^{i,N}\in \mathbb{R} $, $(W^{i})_{i}$ are independent Brownian motions.

\begin{figure}[h!bt]
    \centering
		\begin{subfigure}{.32\textwidth}
			\centering
 			\includegraphics[scale=0.4]{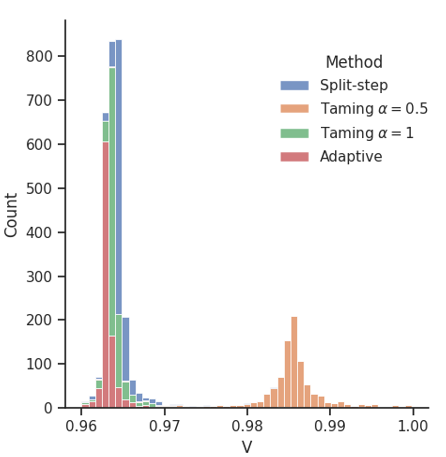}  
			\caption{Histogram of $V$ at $T=1$}
			\label{fig:cs2}
		\end{subfigure}
		\begin{subfigure}{.32\textwidth}
			\centering
 			\includegraphics[scale=0.4]{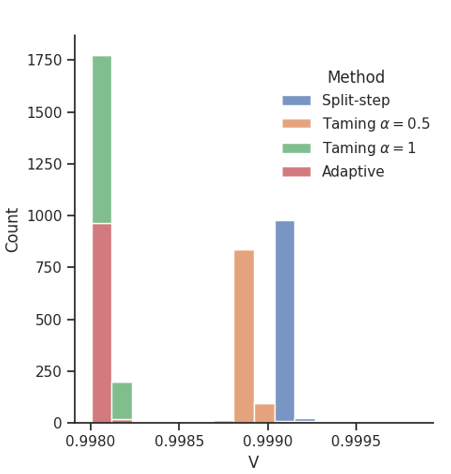} 
			\caption{Histogram of $V$ at $T=2$}
			\label{fig:cs1}
		\end{subfigure}
		\begin{subfigure}{.32\textwidth}
			\centering
 			\includegraphics[scale=0.3]{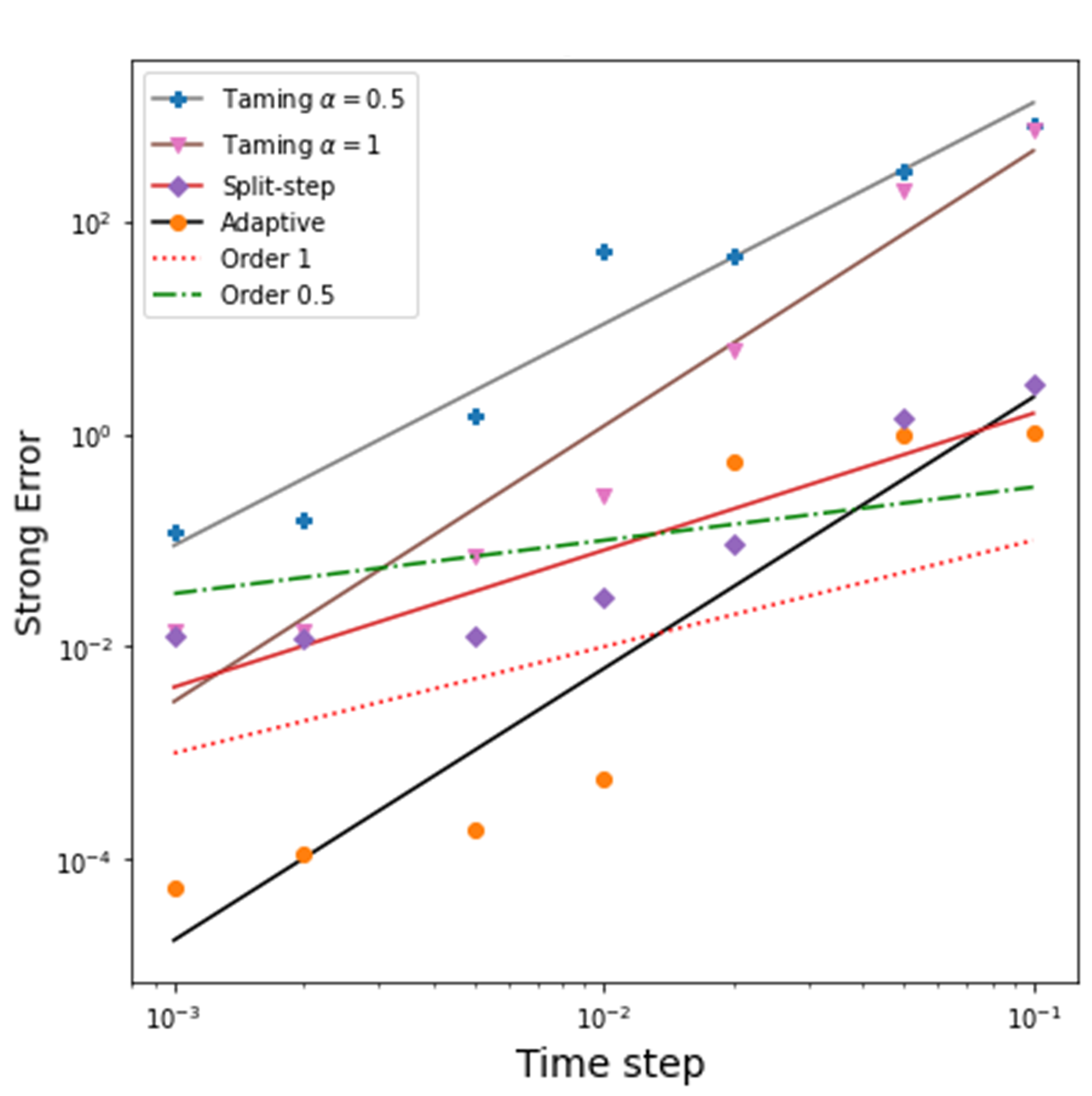}  
			\caption{Strong Error of $V$ w.r.t Time step }
			\label{fig:cs4}
		\end{subfigure}
    \caption{Simulations of the Cucker-Smale type flocking model.  (a,b,c) Histogram of $V$ at different time $T=1,2$ of different methods with $h=10^{-3}$.  (d) Strong error of different methods at $T=2$.
    }
    \label{fig:Simulation_cs}
\end{figure}

Take $\lambda=2,~\sigma'=4,~T=2$.  
With this choice of parameters the solution process $V_t$ converges to $1$ as $t\to \infty$ \cite{csreference}. 
  $V_{0}^{i, N}\sim \mathcal{N}(1,2)$ are i.i.d. samples from standard normal distribution, the timesteps $h$ are $h\in \{ 10^{-3},2\times10^{-3}, \dots, 10^{-1} \}$, particles $N=1000$. The true solution is calculated under $h=10^{-4}$. Taming is run with $\alpha=0.5$ and $1$ and Adaptive with $\mathbf{h}^{\delta}(x)=h\min(1,|x|^{-2})$. 
Figure \ref{fig:Simulation_cs} (a,b) show the distribution of $V$ at time $T=1,2$. All four methods have same initial distribution (and same filtration) nonetheless there is a slight skew between the final results. Taming with $\alpha=0.5$ has a different distribution than hte other three methods at $T=1$, later, the SSM clusters at a different point than Adaptive and Taming method with $\alpha=1$, but the deviation is very small ($<10^{-2}$). Consider the strong error graph (c), the two Taming methods fail to have acceptable result with larger timestep; while SSM and Adaptive are at a similar position. 
For Adaptive, the error rate is nicely behaved but there is a jump at $h=0.02$. The split-step method error rate decrease is stable as $h$ decreases.

\subsection{Discussion}
\label{sec:numericsdiscussion}

We discuss some comparative advantages between the methods starting with generalist comments. All schemes have the same convergence error rate $\textrm{rMSE} \approx Ch^\frac12$. Taming is by far the easiest to implement, with Adaptive the most complex requiring tuning the $\mathbf{h}^{\delta}$ map for each case (see \ref{sec:ShortRecap}). The SSM requires an implicit solver and ad-hoc choices of $v$ and $b$ for efficiency. Taming is the fastest algorithm with SSM and Adaptive running times being comparable with each other. In the way we presented the SSM: all methods are amenable to an efficient parallel implementation (under the caveat of processor communication \cite{bernal2017}); moreover, in view of Remark \ref{remark:removing constraint on Lv}, the SSM does not have any (real) restriction on the time-stepping although one solves an implicit method.

From the numerical examples, we see that 
\begin{enumerate}
    \item the strong error of the SSM is consistently one order of magnitude smaller than that of Taming. Under same choice of timestep, Taming is the fastest with SSM comparable to Adaptive.  {  However, to reach the same strong error level, SSM takes less computational time than Adaptive and significantly less than Taming. }
    
    \item Compared to Adaptive, the SSM has in general no worse convergence than it and no clear domination of one over the other emerged. 
    Implementation wise (at the level of computing the rMSE), to keep the same filtration for different timestep choices, the Brownian motion paths for Adaptive with function $\mathbf{h}^{\delta}(x)$ is much harder to generate (requiring sub-simulation from Brownian bridges) than the SSM with a fixed  timestep. As a rule of thumb, Adaptive does on average a double amount of timesteps than SSM or Taming \cite{reisinger2020adaptive}.

    \item From the numerical examples and at the level of the strong error, the SSM performs better than Taming and Adaptive at larger time steps $h$ (via comparative lower errors). 
\end{enumerate}
We have not investigated the effect of dimensionality, and we suspect that the running time gap of Taming between SSM or Adaptive will widen. The SSM we present has the extra advantage of flexibility in the way of how $v$ and $b$ are chosen. This means that a layer of optimisation can be added to the implicit solver. Lastly, and partially addressed here with the stability analysis, do the schemes preserve the finer properties of the underlying dynamics? Are they geometrically ergodic? Do they preserve oscillatory dynamics, such as amplitudes, frequencies and phases of oscillations? Even for large time steps? 
It is known that explicit Euler type schemes face difficulties in regards to this, with implicit or splitting methods being more stable \cite{buckwarbrehier2021FHNmodelandsplittingBSTT2020}.

\section{Proof of the convergence result for the split-step method (SSM)}
\label{sec:proofOfSSM}

Throughout this section Assumption \ref{Ass:Monotone Assumption} is assumed to hold for all results.

\subsection[Proof of the main convergence result]{Proof of the main convergence result, Theorem \ref{theo:SSTM convergence rate for MV-sde}}
\label{sec:ProofOfConv}

For all auxiliary results next, we assume the conditions of Theorem \ref{theo:SSTM convergence rate for MV-sde} are in force and we thus do not state them.

\subsubsection*{Preliminary results}
As a first step, we state a result that allows us to re-write \eqref{eq:SSTM:scheme 1} and \eqref{eq:SSTM:scheme 2} as a map of $\hat X^{i,N}$ without the presence of the $Y^{i,\star,N}$. We present first a new general version of \cite[{Lemma 3.4}]{higham2002strong} where the differentiability Assumption is lifted and the maps are allowed to depend on time (and the measure component).

\begin{lemma}
\label{lemma:SSTM:new functions def and properties1}
Let $v$ be as in Assumption \ref{Ass:Monotone Assumption}. Choose $h>0$ satisfying $1-h(2L_v+1)>0$. Then for $t\in[0,T],~c,d\in \mathbb{R}^d$, $ \mu\in \cP_2(\bR^d)$ the implicit equation,  
 
with $d, \mu,t$ fixed and $c$ unknown,  
\begin{align}\label{sstmmod func def}
    c=d+h  v(t,c,\mu).
\end{align}
has a unique solution in $c$. Define the functions $v_h$ and $F_h$ as 
\begin{align}
\label{eq:aux:sstmmod func def}
    v_h(t,d,\mu)=v\Big(t,F_{h}(t,d,\mu),\mu \Big)
    \quad \textrm{with}\quad 
    [0,T]\times \bR^d \times \cP_2(\bR^d) \ni (t,d,\mu)\mapsto  F_{h}(t,d,\mu)=c \in \bR^d
    .
\end{align}
We then have for all $t\in[0,T]$, $x,y \in\mathbb{R}^d $, $ \mu,\mu^x,\mu^y\in \cP_2(\bR^d)$ the following four inequalities,  
\begin{align}
|v_h(t,x,\mu)|\le
&
\frac{|v(t,x,\mu)| }{1-h  L_v},\label{sstmmod func prop1}
\\
|F_{h}(t,x,\mu)-F_{h}(t,y,\mu) |^2\le
&
\frac{|x-y|^2 }{1-2h  L_v},\label{sstmmod func prop2}
\\
|F_{h}(t,x,\mu^x)-F_{h}(t,y,\mu^y) |^2\le
&
\frac{1 }{1-h  (2L_v+1)}\Big(|x-y|^2+L_{\tilde{v}} h \big(W^{(2)}(\mu^x,\mu^y)\big)^2  \Big),\label{sstmmod func prop222}
\\
\langle x-y,v_h(t,x,\mu)-v_h(t,y,\mu)  \rangle\le
&
\frac{L_v }{1-2h  L_v}|x-y|^2.\label{sstmmod func prop3}
\end{align}

For $x_i, y_i\in \bR^d,~i\in \llbracket 1,N\rrbracket$ and $ \mu^x, \mu^y\in \cP_2(\bR^d)$ being the empirical measures associated with the collections $\{x_i\}_i, \{y_i\}_i$, define the maps 
\begin{align*}
&b_h(t,x_i,\mu^{F_{h,x,\mu^x}})=b(t,F_{h}(t,x_i,\mu^x),\mu^{F_{h,x,\mu^x}})
\quad \textrm{and}\quad   \sigma_h(t,x_i,\mu^{F_{h,x,\mu^x}})=\sigma(t,F_{h}(t,x_i,\mu^x),\mu^{F_{h,x,\mu^x}}),
\\
&\textrm{where }\ \mu^{F_{h,x,\mu^x}}(dx)=\frac{1}{N}\ \sum_{j=1}^N \delta_{F_{h}(t,x_j,\mu^y)}(dx) 
\quad \textrm{and}
\quad 
    \quad    \mu^{F_{h,y,\mu^y}}(dx)=\frac{1}{N}\ \sum_{j=1}^N \delta_{F_{h}(t,y_{j},\mu^y)}(dx). 
\end{align*}
then, $ b_h$ and $ \sigma_h$ are satisfy
\begin{align*}
    \Big|b_h(t,x_i,\mu^{F_{h,x,\mu^x}})-b_h(t,y_i,\mu^{F_{h,y,\mu^y}})\Big|^2
    \le\frac{L_b}{1-h  (2L_v+1)} \Bigg( |x_i-y_i|^2
    +\frac{2L_{\tilde{v}} h+1}{N} \sum_{j=1}^N  | x_j-y_{j}|^2\Bigg),
    \\
    \Big|\sigma_h(t,x_i,\mu^{F_{h,x,\mu^x}})-\sigma_h(t,y_i,\mu^{F_{h,y,\mu^y}})\Big|^2
    \le
    \frac{L_\sigma}{1-h  (2L_v+1)} \Bigg( |x_i-y_i|^2
    +\frac{2L_{\tilde{v}} h+1}{N} \sum_{j=1}^N  | x_j-y_{j}|^2\Bigg). 
\end{align*}
Lastly, $v_h\rightarrow v$, $b_h\rightarrow b$, and $ \sigma_h\rightarrow \sigma$ uniformly over the compacts of $[0,T]\times \bR^d \times \cP_2(\bR^d)$ as $h \rightarrow 0^+$. 
\end{lemma}
We observe that to establish \eqref{sstmmod func prop1} one only needs $1-hL_v>0$, in other words: the condition $1-2hL_v>0$ is not sharp for that result. Nonetheless, inequalities \eqref{sstmmod func prop2} and \eqref{sstmmod func prop3} are critical for our work and hence we    write   one single constraint.

\begin{proof}

Existence and uniqueness for \eqref{sstmmod func def} can be proved via a strict monotonicity contraction argument. Namely, fix some $t\in[0,T]$, from \eqref{sstmmod func def} one defines the operator $A:\bR^{d}\to\bR^d$ as $A(u)=u-hv(t,u,\mu)$ for $u\in \bR^d$. Following \cite[Definition 25.2 (p.500)]{Zeidler1990B}, the operator $A$ is continuous and strongly monotone (uniformly in $t$) under the restrictions $h>0$ and $1-h L_v>0$. This follows by directly injecting the one-sided condition of $v$ (from Assumption \ref{Ass:Monotone Assumption}) in the definition of \textit{strongly monotone operator}. Finally, from \cite[Theorem 26.A (p.557)]{Zeidler1990B} we conclude that the operator $A$ is invertible and the inverse map is Lipschitz continuous. Thus,  \eqref{sstmmod func def} has a unique measurable inverse given by $F_h$ from \eqref{eq:aux:sstmmod func def}. See also \cite[(p.2596)]{LionnetReisSzpruch2015}. 
%%%%%%%%%%%%%%%%%%%%

We now determine the Lipschitz constant of $v_h$ and $F_h$ of \eqref{eq:aux:sstmmod func def}. For \eqref{sstmmod func prop1}, suppose $c,d\in \mathbb{R}^d$, $\mu\in \cP_2(\bR^d)$ satisfy $c=d+h  v(t,c,\mu)$ then, from Assumption \ref{Ass:Monotone Assumption}:
\begin{align*}
c&=d+  h\Big(v(t,c,\mu)-v(t,d,\mu)\Big)+h  v(t,d,\mu),
\\
\Rightarrow ~ |c-d|^2
&
=\Big\langle c-d,v(t,c,\mu)-v(t,d,\mu) \Big\rangle  h+\langle c-d,v(t,d,\mu) \rangle  h
% \\
% &
\le
h    L_v|c-d|^2 +\langle c-d,v(t,d,\mu) \rangle  h
\\
& \Leftrightarrow 
(1-hL_v)~|c-d|^2  \le |c-d|~|v(t,d,\mu)|~h.
\end{align*}
Since $c=d+h  v_h(t,d,\mu)$, we have by re-arranging the terms and plugging the inequality above 
\begin{align*}
     |v_h(t,d,\mu)|=\frac{1}{  h}|c-d|\le \frac{|v(t,d,\mu)| }{1-h  L_v}.
\end{align*}

For \eqref{sstmmod func prop2}, suppose $c_1,c_2,d_1,d_2\in \mathbb{R}^d,~\mu\in \cP_2(\bR^d)$ satisfy $c_1=d_1+ h  v(t,c_1,\mu)$ and $c_2=d_2+h  v(t,c_2,\mu)$, then
\begin{align*}
\Big|F_{h}(t,d_1,\mu)-F_{h}(t,d_2,\mu) \Big|^2
=|c_1-c_2|^2
&= \langle c_1-c_2 ,d_1-d_2\rangle+ \Big\langle c_1-c_2, v(t,c_1,\mu)- v(t,c_2,\mu) \Big\rangle h 
\\
&
\le \frac{1}{2} |c_1-c_2|^2 +\frac12 |d_1-d_2|^2 +h    L_v|c_1-c_2|^2 .
\end{align*}

For \eqref{sstmmod func prop222}, suppose $x,y,\hat{x},\hat{y}\in \mathbb{R}^d,~\mu^x,\mu^y\in \cP(\bR^d)$ satisfy $\hat{x}=x+h  v(t,\hat{x},\mu^x),~\hat{y}=y+h  v(t,\hat{y},\mu^y)$. We then have
\begin{align*}
\Big|F_{h}(t,x,\mu^x)-&F_{h}(t,y,\mu^y) \Big|^2
=|\hat{x}-\hat{y}|^2
= \langle \hat{x}-\hat{y} ,x-y\rangle+ \Big\langle \hat{x}-\hat{y}, v(t,\hat{x},\mu^x)- v(t,\hat{y},\mu^y) \Big\rangle h 
\\
&
\le \frac{1}{2} |\hat{x}-\hat{y}|^2 +\frac12 |x-y|^2 +\Big\langle \hat{x}-\hat{y}, v(t,\hat{x},\mu^x)- v(t,\hat{y},\mu^x) \Big\rangle h +\Big\langle \hat{x}-\hat{y}, v(t,\hat{y},\mu^x)- v(t,\hat{y},\mu^y) \Big\rangle h 
\\
&
\le  \frac{1}{2} |\hat{x}-\hat{y}|^2 +\frac12 |x-y|^2 + (L_v+\frac{1}{2})| \hat{x}-\hat{y} |^2 h+\frac{1}{2}  L_{\tilde{v}}  \big(W^{(2)}(\mu^x,\mu^y)\big)^2.
\end{align*}

To prove \eqref{sstmmod func prop3} we use the same notation/identities used to prove Inequality \eqref{sstmmod func prop2} above. We have that 
\begin{align*}
\Big\langle d_1-d_2, (d_1-d_2)+h  \Big( v_h(t,d_1,\mu)-v_h(t,d_2,\mu) \Big) \Big\rangle
&=\langle  d_1-d_2 , c_1-c_1\rangle
\le \frac{1}{2}|d_1-d_2 | ^2+\frac{1}{2}\frac{|d_1-d_2|^2 }{1-2 h   L_v} ,
\end{align*}

and thus
\begin{align*}
\Big\langle  d_1-d_2, v_h(t,d_1,\mu)-v_h(t,d_2,\mu)  \Big\rangle
&\leq \frac{L_v }{1-2 h   L_v}     |d_1-d_2|^2.
\end{align*}

We now address the Lipschitz property of $b_h,\sigma_h$.  
Since they are  of the same nature, we provide  only the proof for $b_h$ as that for $\sigma_h$ is identical.
Let $i\in \llbracket 1,N\rrbracket$, using the definition of $b_h$, then Assumption \ref{Ass:Monotone Assumption} followed by  \eqref{sstmmod func prop222} we have  
 
\begin{align*}
    \Big|b_h(t,x_i,\mu^{F_{h,x}})-b_h(t,y_i,\mu^{F_{h,y}})\Big|^2
    &\le L_b\Big(|F_{h}(t,x_i,\mu^x)-F_{h}(t,y_i,\mu^y)|^2
    +\big(W^{(2)}(\mu^{F_{h,x}},\mu^{F_{h,y}})\big)^2 \Big)
    \\
    & \le L_b\Big(|F_{h}(t,x_i,\mu^x)-F_{h}(t,y_i,\mu^y)|^2
    +\frac{1}{N} \sum_{j=1}^N | F_{h}(t,x_j,\mu^x)-F_{h}(t,y_{j},\mu^y) |^2\Big)
    \\
     & \le \frac{L_b}{1-h  (2L_v+1)} \Bigg( |x_i-y_i|^2
    +\frac{2L_{\tilde{v}} h+1}{N} \sum_{j=1}^N  | x_j-y_{j}|^2\Bigg).
\end{align*}
The convergence result in the final statement follows straightforwardly from \cite[Lemma 3.4]{higham2002strong}. This convergence result is applied with fixed $N$ and the parameter of the convergence is $h$ (not $N$). One only needs to apply their arguments over $[0,T]\times \bR^{dN}$ where the measures $\mu\in \cP_2(\bR^d)$ are taken to be compactly supported on the compact where the family of points $\{x_i\}_i$ is contained.

\end{proof}

After having introduced Lemma \ref{lemma:SSTM:new functions def and properties1}, we can finally address the continuous-time extension of the SSM  \eqref{eq:SSTM:scheme 1}-\eqref{eq:SSTM:scheme 2} as referenced in Theorem \ref{theo:SSTM convergence rate for MV-sde}. The SSM can be written as a continuous time SDE via linear interpolation of the iterates, namely, for $t\in[t_n,\tnp]$,  $i\in\llbracket 1,N \rrbracket,~\hat{X}_{0}^i
    \in L_0^m(\mathbb{R}^d)$: 
\begin{align}
\label{eq: SSTM discrete scheme sde with kappa t}
    d\hat{X}_{t}^{i,N}
    &
    =
    \Big( v_h\big(\kt,\hx_{\kt}^{i,N},\hat{\mu}^{N}_{\kt}\big)
    +b_h\big(\kt,\hx_{\kt}^{i,N},\widetilde {\mu}_{\kt}^{N}\big)  \Big) dt 
    + \sigma_h\big(\kt,\hx_{\kt}^{i,N},\widetilde {\mu}_{\kt}^{N}\big) dW_t^i,
    \\
    \nonumber
    \textrm{where }\quad
    \hat{\mu}^{N}_{\kt}(dx)&:= \frac{1}{N} \sum_{j=1}^N \delta_{\hx_{\kt}^{j,N} }(dx)
    \quad\textrm{and}\quad
    \widetilde {\mu}^{N}_{\kt}(dx):= \frac{1}{N} \sum_{j=1}^N \delta_{F_{h}(\kt,\hx_{\kt}^{j,N}, \hat{\mu}^{N}_{\kt})}(dx).
\end{align}
where $\kappa(t)=\sup\Big\{t_n: t_n\le t,\ n\in \llbracket 0,M \rrbracket \Big\}$ and $\widetilde {\mu}^{N}_{t_n}=\widetilde {\mu}^{N}_{n}$.

\subsubsection*{Moment bounds}

We now employ the results of Lemma \ref{lemma:SSTM:new functions def and properties1} to establishing a domination of $|Y_{n}^{i,\star,N}|$ by $|\hx_{n}^{i,N}|$.
\begin{lemma}\label{lemma：SSTM:compare betwen Y star and X}
Choose $h$ as in Theorem \ref{theo:SSTM convergence rate for MV-sde} and recall
$C_T,\widehat L_v$ as defined in Remark \ref{rem:ImpliedProperties}. 

Then, $|Y_{n}^{i,\star,N}|$ of \eqref{eq:SSTM:scheme 1} satisfies for any $i\in \llbracket 1,N \rrbracket ,\  n\in \llbracket 0,M-1 \rrbracket$,
\begin{align}
\label{eq:Y sqaure le X square}
    |Y_{n}^{i,\star,N}|^2 
    &
    \le |\hx_{n}^{i,N}|^2(1+\frac{2\widehat L_v  }{1-2\widehat L_v  h}h)
    + \frac{ L_{\tilde{v}}  h }{1-2\widehat L_v  h} \Bigg( \frac{1}{N}  \sum_{j=1}^N |\hx_{n}^{j,N}|^2 \Bigg)
    +\frac{2C_T  }{1-2\widehat L_v  h} h.
\end{align}
\end{lemma}

\begin{remark}[More on the stepsize $h$]
\label{remark:1-hLv is a constant like 1/2}
The bound for $h$ is necessary, for example, if $v(t,x,\mu)=10x$, then the implicit solution gives $y=x/(1-10h)$ and one sees that the $h<0.1$ condition is critical. 

Recall the stepsize constraint on $h$ in Theorem \ref{theo:SSTM convergence rate for MV-sde}. By inspection of the proof of this Lemma and the definition of $\widehat L_v$ in Remark \ref{rem:ImpliedProperties}, the reader will find that the constraint $L_v<-1/2$ needed to ensure $1-h\widehat L_v>0 $ is not sharp and can be replaced by some number $\varepsilon\in (0,1)$, i.e., $L_v<-\varepsilon$.  The lack of sharpness arises from the choice of $\widehat L_v$ in Remark \ref{rem:ImpliedProperties}. There we used the Cauchy-Schwarz inequality where we could have used a Young type inequality from which the parameter $\varepsilon$ would have arisen. We choose $1/2$ for ease of presentation.

However, through the split-step structure, one can use the ``add and subtract a linear component'' in the drift before the split-step is executed to make $L_v<0$ and thus remove the constraint on $h$ -- see Section \ref{section:OU example}.
\end{remark}

\begin{proof}

From \eqref{eq:SSTM:scheme 1} and Remark \ref{rem:ImpliedProperties}, for any $i,n$ and any $t_n\in\pi$  we have using Cauchy-Schwarz and the properties of $v$ that
\begin{align*}
|Y_{n}^{i,\star,N}|^2 &=\langle Y_{n}^{i,\star,N},\hx_{n}^{i,N} \rangle + \langle  Y_{n}^{i,\star,N}, v(t_n,Y_{n}^{i,\star,N},\hat{\mu}^{N}_n	) \rangle   h
\\
&\le \frac{1}{2}|Y_{n}^{i,\star,N}|^2 + \frac{1}{2}|\hx_{n}^{i,N}|^2 +    h\Bigg(C_T+\widehat L_v |Y_{n}^{i,\star,N}|^2+\frac{L_{\tilde{v}}}{2}\Big( \frac{1}{N}  \sum_{j=1}^N |\hx_{n}^{j,N}|^2 \Big)  \Bigg).
\end{align*}
\end{proof}

\begin{remark}
\label{rem:constants independent of h funcproperty y star x}
Under Assumption \ref{Ass:Monotone Assumption} and under the choice of $h$ in Theorem \ref{theo:SSTM convergence rate for MV-sde}, $1/(1-2\widehat L_v   h)$ {  is bounded above by some constant independent of $h$}. We can thus claim that there exist constants $C>0,~\widetilde C\in \bR$ depending on $L_v,C_T, \widehat L_v$ but independent of $h$ such that from Lemma \ref{lemma:SSTM:new functions def and properties1} and Lemma \ref{lemma：SSTM:compare betwen Y star and X} we have (for any $t,x,y,\mu$)
\begin{align*}
&|v_h(t,x,\mu)|\le
C{|v(t,x,\mu)| }, 
\quad
&|F_{h}(t,x,\mu)-F_{h}(t,y,\mu) |^2\le
C |x-y|^2,
\\
&\Big\langle x-y,v_h(t,x,\mu)-v_h(t,y,\mu)  \Big\rangle\le
\widetilde C|x-y|^2,
\quad
&|Y_{k}^{i,\star,N}|^2 
\le |\hx_{k}^{i,N}|^2(1+C h)+h \frac{C}{N} \sum_{j=1}^N |\hx_{k}^{j,N}|^2 +C h.
\end{align*}
\end{remark}
After this remark emphasising the independence of the constants in $h$, we are in a position to prove the moment bounds for the output of the SSM. 

\begin{Proposition}[Moment bounds of SSM] 
\label{prop:SSTM: discrete moment bounded}
Choose $h$ as in Theorem \ref{theo:SSTM convergence rate for MV-sde}. Then,  there exists a constant $C\in\mathbb{R^{+}}$ such that for any $ i\in \llbracket 1,N \rrbracket ,\  n\in \llbracket 0,M \rrbracket$,  for $\frac{m}{2} \ge p\ge 1$, the output $\hx_{n}^{i,N}$ of the scheme \eqref{eq:SSTM:scheme 1}-\eqref{eq:SSTM:scheme 2} satisfies, 

\begin{align*}
    \sup_{ i\in \llbracket 1,N \rrbracket} 
    \bE \Big[\sup_{n\in \llbracket 0,M \rrbracket } |\hx_{n}^{i,N}|^{2p}  \Big]
    <C \Big(  1 
              +  \sup_{ i\in \llbracket 1,N \rrbracket} 
                \bE\Big[\, |\hx_{0}^{i,N}|^{2p}\Big]
    \Big) <\infty.
\end{align*}
\end{Proposition}
\begin{proof}
Let $i\in \llbracket 1,N \rrbracket$, ${n\in \llbracket 0,M-1 \rrbracket }$ and recall \eqref{eq:SSTM:scheme 2}.  
Using Assumption \ref{Ass:Monotone Assumption}, Remark \ref{rem:constants independent of h funcproperty y star x}, we have by taking squares and applying Young's inequality: 
\begin{align*}
     | \hx_{n+1}^{i,N}|^2
    & =
    |Y_{n}^{i,\star,N}|^2
     +\Big| b(t_n,Y_{n}^{i,\star,N},\hm_{n}^{Y,N})h+ \sigma(t_n,Y_{n}^{i,\star,N},\hm_{n}^{Y,N})\Delta W_n^i\Big|^2
    \\
    &\qquad +2 \Big\langle Y_{n}^{i,\star,N}, b(t_n,Y_{n}^{i,\star,N},\hm_{n}^{Y,N})h+ \sigma(t_n,Y_{n}^{i,\star,N},\hm_{n}^{Y,N})\Delta W_n^i \Big\rangle
    \\
    &
    \le  
    |\hx_{n}^{i,N}|^2(1+C h)+ \frac{Ch}{N} \sum_{j=1}^N |\hx_{n}^{j,N}|^2 +C h
    + 2 | b(t_n,Y_{n}^{i,\star,N},\hm_{n}^{Y,N})|^2h^2+2| \sigma(t_n,Y_{n}^{i,\star,N},\hm_{n}^{Y,N})|^2(\Delta W_n^i)^2
    \\
    &\qquad  + 2\Big\langle Y_{n}^{i,\star,N}, \sigma(t_n,Y_{n}^{i,\star,N},\hm_{n}^{Y,N})\Delta W_n^i \Big\rangle
    + 2 C\Big(1+|Y_{n}^{i,\star,N}|^2+\frac{1}{N} \sum_{j=1}^N |Y_{n}^{j,\star,N}|^2\Big)h
    \\
    & 
    \le 
    |\hx_{n}^{i,N}|^2 +C\Big(1+|\hx_{n}^{i,N}|^2+\frac{1}{N} \sum_{j=1}^N |\hx_{n}^{j,N}|^2 \Big) \Big( |\Delta W_k^i|^2+h\Big)  
    +2\Big\langle Y_{n}^{i,\star,N}, \sigma(t_n,Y_{n}^{i,\star,N},\hm_{n}^{Y,N})\Delta W_n^i \Big\rangle.
\end{align*}
where we used $W^{(2)}(\hm_{n}^{Y,N}, \delta_0)^2\le \frac{1}{N} \sum_{j=1}^N |Y_{n}^{j,\star,N}|^2$.  By backward induction from $n+1$ to zero, we have (after some simplification)
\begin{align*}
    | \hx_{n+1}^{i,N}|^2 
  \leq 
  |\hx_{0}^{i,N}|^2 
   & +    
    2\sum_{k=0}^{n} \Big\{  \Big\langle Y_{k}^{i,\star,N}, \sigma(t_k,Y_{k}^{i,\star,N},\hm_{k}^{Y,N})\Delta W_k^i \Big\rangle  \Big\}
    +C \sum_{k=0}^{n} \Big\{ |\Delta W_k^i|^2 \Big\}
    +C \sum_{k=0}^{n} \Big\{ |\hx_{k}^{i,N}|^2 |\Delta W_k^i|^2 \Big\}
\\
    &
    +C \sum_{k=0}^{n} \Big\{\frac{1}{N} \sum_{j=1}^N |\hx_{k}^{j,N}|^2 |\Delta W_k^i|^2 \Big\}
    +C \sum_{k=0}^{n} \Big\{ h \Big\}
    +C \sum_{k=0}^{n} \Big\{ |\hx_{k}^{i,N}|^2  h \Big\}
    +C \sum_{k=0}^{n} \Big\{ \frac{1}{N} \sum_{j=1}^N |\hx_{k}^{j,N}|^2 h \Big\}.
\end{align*}
Taking power $p\geq 2$ on both sides, expectations and re-organising the terms, we have (with $C_p$ independent of $h,N$) 
\begin{align*}
    \bE \Big[   | \hx_{n+1}^{i,N}|^{2p} \Big]
    \le & 
    C_p \Bigg\{   \bE \Big[   |\hx_{0}^{i,N}|^{2p} \Big]
    + 
     \bE\Big[  \Big(    \sum_{k=0}^{n} \Big\{  \langle Y_{k}^{i,\star,N}, \sigma(t_k,Y_{k}^{i,\star,N},\hm_{k}^{Y,N})\Delta W_k^i \rangle  \Big\}\Big)^p\Big]
     + 
    \bE\Big[ \Big(   \sum_{k=0}^{n}\Big\{   |\Delta W_k^i|^2\Big\}\Big)^p\Big]
    \\
    & 
    +
    \bE\Big[\Big(      \sum_{k=0}^{n}\Big\{ |\hx_{k}^{i,N}|^2 |\Delta W_k^i|^2\Big\}\Big)^p\Big]
    +
    \bE\Big[ \Big(     \sum_{k=0}^{n}\Big\{ \frac{1}{N}\sum_{j=1}^N |\hx_{k}^{j,N}|^2 |\Delta W_k^i|^2\Big\}\Big)^p\Big]
    \\
    &+
    \bE\Big[\Big(   p   \sum_{k=0}^{n}\Big\{ h\Big\}\Big)^p\Big]
    + 
    \bE\Big[ \Big(     \sum_{k=0}^{n}\Big\{ |\hx_{k}^{i,N}|^2 h\Big\}\Big)^p \Big] 
    +
    \bE\Big[ \Big(     \sum_{k=0}^{n}\Big\{ \frac{1}{N}\sum_{j=1}^N |\hx_{k}^{j,N}|^2 h\Big\}\Big)^p \Big] \Bigg\}.
 \end{align*}

There are 8 terms to be estimated, but in essence only 3 arguments are needed. We present them only for the most complex terms since for the remaining ones it is just a simplification of the arguments presented. We present them in the form of supremum over $n\in \llbracket 0,M-1 \rrbracket$.
We start with the last term of the 2nd line: apply Jensen's inequality twice after scaling the outer summation and then tower property to take advantage of the conditional independence between $\hx_{k}^{j,N}$ and $\Delta W_k^i$, namely, (recall $h=T/M$) 
\begin{align*}
       \bE \Big[ \sup_{n\in \llbracket 0,M-1 \rrbracket}\Big(     \sum_{k=0}^{n} \Big\{\frac{1}{N}\sum_{j=1}^N |\hx_{k}^{j,N}|^2 |\Delta W_k^i|^2\Big\} \Big)^p  \Big] 
       & \le \bE \Big[  \frac{1}{M} \sum_{k=0}^{M-1}\Big\{ \frac{1}{N}\sum_{j=1}^N |\hx_{k}^{j,N}|^{2p} |\Delta W_k^i|^{2p}\Big\}  M^p \Big]\\
      &  =h \bE \Big[  \sum_{k=0}^{M-1} \frac{1}{N}\sum_{j=1}^N |\hx_{k}^{j,N}|^{2p}  \Big] C T^p  .
\end{align*}

We now address the second term in the 1st line. Using Burkholder–Davis–Gundy (BDG) inequality and Jensen's inequality as above, the Lipschitz property of $\sigma$, Jensen's inequality and the domination of $Y_{k}^{i,\star,N}$ by $\hx_{k}^{i,N}$ in Lemma \ref{lemma：SSTM:compare betwen Y star and X} gives

\begin{align*}
    \bE \Bigg[ \sup_{n\in \llbracket 0,M-1 \rrbracket} \Bigg(   \sum_{k=0}^{n} \Big\{  \langle Y_{k}^{i,\star,N}
    ,\sigma(t_k&,Y_{k}^{i,\star,N},\hm_{k}^{Y,N})\Delta W_k^i \rangle  \Big\}\Bigg)^p\Bigg]
    \le
    C    \bE \Bigg[ \Bigg(\sum_{k=0}^{M-1} \Big\{  |\langle Y_{k}^{i,\star,N}
    ,\sigma(t_k,Y_{k}^{i,\star,N},\hm_{k}^{Y,N})\rangle|^2 h  \Big\}\Bigg)^\frac{p}{2} \Bigg]
    \\
    \le& 
        C    \bE \Bigg[ 
        \Bigg(
        \frac1M\sum_{k=0}^{M-1} \Big\{  1+|Y_{k}^{i,\star,N}|^4+\frac{1}{N} \sum_{j=1}^N |Y_{k}^{j,\star,N}|^4   \Big\}h 
        \Bigg)^\frac{p}{2}
        \Bigg] M^\frac{p}{2}
    \\
         \le& 
        C    \bE \Bigg[ 
        1
        +
        \Big[\sum_{k=0}^{M-1} |\hx_{k}^{i,N}|^{2p}\Big] h
        +
         \Big[\sum_{k=0}^{M-1} \frac{1}{N} \sum_{j=1}^N |\hx_{k}^{j,N}|^{2p} \Big] h
           \Bigg].
\end{align*}
Finally, we address the last term in the 3rd line of the inequality. The result follows by applying Jensen's inequality
\begin{align*}
    \bE \Big[  \sup_{n\in \llbracket 0,M-1 \rrbracket} \Big(\sum_{k=0}^{n} \frac{1}{N}\sum_{j=1}^N |\hx_{k}^{j,N}|^2 h\Big)^p  \Big] 
    &
    = 
    \bE \Big[  \Big(\frac{1}{M} \sum_{k=0}^{M-1} \frac{1}{N}\sum_{j=1}^N |\hx_{k}^{j,N}|^2 h\Big)^p M^p \Big]
    \leq 
    h \bE \Big[  \sum_{k=0}^{M-1} \frac{1}{N}\sum_{j=1}^N |\hx_{k}^{j,N}|^{2p}  \Big]T^p.
\end{align*}
Collecting the several inequalities and injecting them back in the initial inequality, we conclude that 
\begin{align*}
    \bE \Big[\sup_{n\in \llbracket 0,M-1 \rrbracket}| \hx_{n+1}^{i,N}|^{2p} \Big]
    \le 
    & C \Bigg\{  
    1
    + \bE\Big[\, |\hx_{0}^{i,N}|^{2p} \Big]
    + \bE \Big[    \sum_{k=0}^{M-1}  |\hx_{k}^{i,N}|^{2p}h  \Big] 
    +  \bE \Big[  \sum_{k=0}^{M-1} \frac{1}{N}\sum_{j=1}^N |\hx_{k}^{j,N}|^{2p} h \Big] \Bigg\}.
\end{align*}
Taking supremum and using that the particles are conditional i.d.d.~(for fixed $k$)
\begin{align*}
    \sup_{ i\in \llbracket 1,N \rrbracket} \bE \Big[\sup_{n\in \llbracket 0,M \rrbracket}|\hx_{n}^{i,N}|^{2p}  \Big] 
    &\le C \Bigg\{  1
                   + \sup_{ i\in \llbracket 1,N \rrbracket} \bE\Big[\, |\hx_{0}^{i,N}|^{2p}\Big]
                   +\sum_{k=0}^{M-1}   \sup_{ i\in \llbracket 1,N \rrbracket} \bE \Big[\sup_{0\le n \le k}|\hx_{n}^{i,N}|^{2p}  \Big] h    \Bigg\}.
\end{align*}
The proof finishes after applying the discrete Gr\"onwall's inequality to the inequality, and using that the $\hx_{0}^{i,N}$ are i.i.d.
\end{proof}

We now provide (moment) estimates for the continuous-time extension of the SSM.
\begin{Proposition}
\label{prop:SSTM：moment bound for the big theorm time extensions }
Choose $h$ as in Theorem \ref{theo:SSTM convergence rate for MV-sde}. 
Take $(\hx_{t}^{i,N})_{t\in[0,T]}$ as the map satisfying \eqref{eq: SSTM discrete scheme sde with kappa t}, i.e., the continuous time extension of the SSM.
 Then,      for any for $\frac{m}{2} \ge p\ge 1$, there exist $C\in\mathbb{R^{+}}$:
\begin{align} \label{eq:SSTM:momentbound for split-step time extension process}
    \sup_{ i\in \llbracket 1,N \rrbracket}   \bE\Big[\sup_{0\le t \le T} |\hx_{t}^{i,N}|^{2p}\Big]<C\Big(  1+   \sup_{ i\in \llbracket 1,N \rrbracket}  \bE\Big[\, |\hx_{0}^{i,N}|^{2p}\Big] \Big) <\infty.
\end{align}
\end{Proposition}

\begin{proof}
Let $i\in \llbracket 1,N \rrbracket$ and ${n\in \llbracket 0,M-1 \rrbracket }$. 
From \eqref{eq: SSTM discrete scheme sde with kappa t}, set $t_n+s=t$, for all $t \in [0,T], n\in \llbracket 0,M-1 \rrbracket$, then 
\begin{align}\label{eq:1:prop:SSTN:moment bound for the big theorm time extensions }
    \hx_{t}^{i,N}
    = \hx_{n}^{i,N}+v_h(t_n,\hx_{n}^{i,N},\hat{\mu}^{N}_n	)s
       +b_h(t_n,\hx_{n}^{i,N},\wm_n^{N}  )s
       +\sigma_h(t_n,\hx_{n}^{i,N},\wm_n^{N})(W_{t_n+s}^i-W_{t_n}^i).
\end{align}
Plugging \eqref{eq:SSTM:scheme 1} in \eqref{eq:1:prop:SSTN:moment bound for the big theorm time extensions } gives 
\begin{align*}
    \hx_{t}^{i,N}= \hx_{n}^{i,N}(1-\frac{s}{  h})+\frac{s}{  h}Y_{n}^{i,\star,N}+b_h(t_n,\hx_{n}^{i,N},\wm_{n}^{N})s+\sigma_h(t_n,\hx_{n}^{i,N},\wm_{n}^{N})(W_{t_n+s}^i-W_{t_n}^i).
\end{align*}
Since $s\le h $, by the Lipschitz conditions on $b_h$ and $\sigma_h$, and Lemma \ref{lemma：SSTM:compare betwen Y star and X}, we have
\begin{align*}
     |\hx_{t}^{i,N}|^2 
     &
     \le C\Bigg\{1+| \hx_{n}^{i,N}|^2+ \frac{1}{N} \sum_{j=1}^N |\hx_{n}^{j,N}|^2+|b_h(t_n,\hx_{n}^{i,N},\wm_{n}^{N})|^2 + \Big|\sigma_h(t_n,\hx_{n}^{i,N},\wm_{n}^{N})(W_{t_n+s}^i-W_{t_n}^i)\Big|^2\Bigg\}
     \\
     &\le C\Bigg\{1+| \hx_{n}^{i,N}|^2+ \frac{1}{N} \sum_{j=1}^N |\hx_{n}^{j,N}|^2 + \Big|\sigma_h(t_n,\hx_{n}^{i,N},\wm_{n}^{N})(W_{t_n+s}^i-W_{t_n}^i)\Big|^2\Bigg\}.
\end{align*}
Taking supremum over time and expectations on both sides yields 
\begin{align*}
    \bE\Big[\sup_{0\le t \le T}|\hx_{t}^{i,N}|^{2p}\Big]
    & =\bE\Big[ \sup_{ n\in \llbracket 0,M-1 \rrbracket}   \sup_{0\le s \le h}| \hx_{t_n+s}^{i,N}|^{2p}\Big]
    \le  C\bE\Bigg[1+ \sup_{ n\in \llbracket 0,M-1 \rrbracket}\Big\{\Big
    (| \hx_{n}^{i,N}|^{2p}+
      \frac{1}{N} \sum_{j=1}^N |\hx_{n}^{j,N}|^2 \Big)
      + I_h^{i,n}
    \Big\}  \Bigg].
\end{align*}
where $I_h^{i,n}$ is given by $I_h^{i,n}:=\sup_{0\le s \le h}  \Big|\sigma_h(t_n,\hx_{n}^{i,N},\wm_{n}^{N})(W_{t_n+s}^i-W_{t_n}^i)\Big|^{2p}$. Using the BDG inequality, Jensen's inequality, the Lipschitz condition on $\sigma_h$ and Proposition \ref{prop:SSTM: discrete moment bounded}, gives
\begin{align*}
     \bE [I_h^{i,n} ]
    &\le C \bE \Big[   \Big(|\sigma_h(t_n,\hx_{n}^{i,N},\wm_{n}^{N})|^{2} h\Big)^{p}  \Big]
    \le 
    C h^p \bE \Big[ 1 +|\hx_{n}^{i,N}|^{2p}+\frac{1}{N} \sum_{j=1}^N |\hx_{n}^{j,N}|^{2p}  \Big]
     \le C h^p.
\end{align*}
Take supremum over $i\in \llbracket 1,N \rrbracket$, then by Proposition \ref{prop:SSTM: discrete moment bounded} it follows that
\begin{align*}
    \sup_{ i\in \llbracket 1,N \rrbracket} \bE\Big
[\sup_{0\le t \le T}|\hx_{t}^{i,N}|^{2p}\Big
]
      &\le \sup_{ i\in \llbracket 1,N \rrbracket}  \bE\Big
[C(1+\sup_{0\le n \le M} | \hx_{n}^{i,N}|^{2p}+  Ch^p )\Big
]
      \le C\Big(  1+ 
       \sup_{ i\in \llbracket 1,N \rrbracket}  
      \bE\Big[\, |\hx_{0}^{i,N}|^{2p}\Big] \Big) <\infty.
\end{align*}
\end{proof}

The last result in this block concerns incremental (in time) moment bounds of $\hx^{i,N}$.
\begin{Proposition}
\label{prop:SSTM: local increment error}
There exists  $C\in \mathbb{R^{+}}$ such that for any $p\geq 2$,     with $m\ge (q+1)p$, 
\begin{align}
    \sup_{ i\in \llbracket 1,N \rrbracket}  
    \bE\Big[  \sup_{0\le t \le T}|\hx_{t}^{i,N}-\hx_{\kt}^{i,N} |^p \Big]
    \le C h^{\frac p2}.
\end{align}
\end{Proposition}
\begin{proof}
From \eqref{eq: SSTM discrete scheme sde with kappa t}, for all $t \in [0,T]$ such that $t\in[t_n,t_{n+1}$] set $s\in[0,h]$ such that $t_n+s=t$. Then 
\begin{align*}
    \hx_{t}^{i,N}= \hx_{n}^{i,N}+v_h(t_n,\hx_{n}^{i,N},\hm_{n}^{N})s+b_h(t_n,\hx_{n}^{i,N},\wm_{n}^{N})s+\sigma_h(t_n,\hx_{n}^{i,N},\wm_{n}^{N})(W_{t}^i-W_{t_n}^i).
\end{align*}
Thus, we have a constant $C_p$ only depending on $p$ such that 
\begin{align*}
    |\hx_{t}^{i,N}- \hx_{n}^{i,N}|^p
    \le
    &  
    C_p \Big[|v_h(t_n,\hx_{n}^{i,N},\hm_{n}^{N})|^p h^p
    +|b_h(t_n,\hx_{n}^{i,N},\wm_{n}^{N})|^p h^p
    +|\sigma_h(t_n,\hx_{n}^{i,N},\wm_{n}^{N})|^p|(W_{t}^i-W_{t_n}^i)|^p  \Big].
\end{align*}
Take $\sup$ over time and expectations on both sides. Using Assumption \ref{Ass:Monotone Assumption}, one deals with the last term with the BDG inequality (using conditional expectations via $\bE[\cdot|\cF_\tn]$) and Jensen's inequality. From Lemma \ref{lemma:SSTM:new functions def and properties1} and Propositions \ref{prop:SSTM: discrete moment bounded} and \ref{prop:SSTM：moment bound for the big theorm time extensions }, there exists a positive $C$ independent of $h,N,M$ such that
\begin{align*}
\sup_{ i\in \llbracket 1,N \rrbracket}  
\bE\Big[  \sup_{0\le t \le T}|\hx_{t}^{i,N}-\hx_{\kt}^{i,N} |^p \Big]
\le C (h^p+h^{\frac{p}{2}}) \Big(1+ \sup_{ i\in \llbracket 1,N \rrbracket}   \bE\Big[\sup_{0\le n \le M} |\hx_{n}^{i,N}|^{(q+1)p}\Big]\Big) \le C h^{\frac{p}{2}}.
\end{align*}
where $q$ follows from the polynomial growth property of $v$ in Assumption \ref{Ass:Monotone Assumption}.
\end{proof}

\subsubsection*{Local errors}

After having discussed moment bounds, we now discuss the local error.
\begin{Proposition}
\label{prop:SSTM:diff functions between orgins and h.s}
Let the assumptions of Theorem \ref{theo:SSTM convergence rate for MV-sde} hold. Take the functions $v,b,\sigma$ and the corresponding functions $v_h,b_h,\sigma_h$ as defined in Lemma \ref{lemma:SSTM:new functions def and properties1}. 

Then, there exist positive constants $C_1,C_2,C_3$ and $q'=2(q+1)^2$, such that for all $t \in[0,T]$, $i\in \llbracket 1,N \rrbracket$, $ z_i \in \bR^d$, and the collection $\{z_i\}_i$, we have 
\begin{align}
    |v_h(t,z_i,\mu^z)-v(t,z_i,\mu^z)|^2
    &
    \le C_1\Big(1+|z_i|^{q'}+\frac{1}{N} \sum_{j=1}^N |z_j|^{q'}\Big) h^2,
    \\
        |b_h(t,z_i,\mu^{F_{h,z,\mu^z}})- b(t,z_i,\mu^z)|^2
    &\le C_2\Big(1+|z_i|^{q'}+\frac{1}{N} \sum_{j=1}^N |z_j|^{q'}\Big) h^2,
    \\
        |\sigma_h(t,z_i,\mu^{F_{h,z,\mu^z}})-\sigma(t,z_i,\mu^z)|^2
    &\le C_3\Big(1+|z_i|^{q'}+\frac{1}{N} \sum_{j=1}^N |z_j|^{q'} \Big) h^2.
\end{align}
where $\mu^z$ and $\mu^{F_{h,z,\mu^z}}$ are the two empirical measures associated with $\{z_i\}_i$ and $ \{F_h(t,z_i,\mu^z)\}_i$ respectively,  i.e.,
\begin{align*}
   \mu^z(dz)=\frac{1}{N}\ \sum_{j=1}^N \delta_{z_j}(dz),
    \quad\textrm{and}\quad
    \mu^{F_{h,z,\mu^z}}(dz)=\frac{1}{N}\ \sum_{j=1}^N \delta_{\{F_{h}(t,z_j,\mu^z)\}}(dz).
\end{align*}
\end{Proposition}
\begin{proof} 
Recall the estimates given in Lemma \ref{lemma:SSTM:new functions def and properties1}. Using the identity  $F_{h}(t,z_i,\mu^z)=z_i+h v_h(t,z_i,\mu^z)$, \eqref{sstmmod func prop1},  Assumption \ref{Ass:Monotone Assumption}, Young's inequality and Jensen's inequality, we have 
\begin{align*}
\Big| v_h(t,z_i,\mu^z) -v(t,z_i,\mu^z) \Big|^2
&= 
\Big|v \Big(t,z_i+   h v_h(t,z_i,\mu^z),\mu^z \Big)-v(t,z_i,\mu^z)\Big|^2
\\
&
\le 
C\Big(1+|z_i+    h v_h(t,z_i,\mu^z)|^q+|z_i|^q\Big)^2h^2 |v_h(t,z_i,\mu^z)|^2
\\
&
\leq 
C\Big(1+|z_i|^{2q}+ |z_i|^{2q(q+1)}h^{2q} +h^{2q}\frac{1}{N} \sum_{j=1}^N |z_j|^{2q} \Big) \frac {h^2}{(1-hL_v)^2} \Big(1+|z_i|^{2q+2}+\frac{1}{N} \sum_{j=1}^N |z_j|^{2}\Big)
\\
&
\leq 
C\Big(1+|z_i|^{2(q+1)(q+1)} 
+\frac{1}{N} \sum_{j=1}^N |z_j|^{4q+2}\Big)h^2.
\end{align*}
As in Lemma \ref{lemma:SSTM:new functions def and properties1} we show only the result for $b_h$ as the computation is the same for $\sigma_h$ (and overall very close to that for $v_h$). Using the definition of $b_h$, the Lipschitz property of $b$ and the definition of $\mu^{F_{h,z,\mu^z}},\mu^z$, using similar calculations as above, by Young's inequality and Jensen's inequality, we have 
\begin{align*}
\Big| b_h(t,z_i,\mu^{F_{h,z,\mu^z}})-b(t,z_i,\mu^z) \Big|^2
& =
\Big|b\Bigg(t,F_{h}(t,z_i,\mu^z),\mu^{F_{h,z,\mu^z}}\Bigg)-b(t,z_i,\mu^z)\Big|^2 
\\ &
\le L_b\Big(h^2|v_h(t,z_i,\mu^z)|^2+ \big(W^{(2)}(\mu^{F_{h,z,\mu^z}} ,\mu^z ) \big)^2\Big)
\\ &
\le L_b \Big(h^2|v_h(t,z_i,\mu^z)|^2+\frac{1}{N} \sum_{j=1}^N | F_{h}(t,z_j,\mu^z)-z_j|^2    \Big)
\\
&\le L_b\Big( h^2|v_h(t,z_i,\mu^z)|^2+\frac{1}{N} \sum_{j=1}^N h^2|v_h(t,z_j,\mu^z)|^2  \Big).
\end{align*}
Applying Inequality \eqref{sstmmod func prop1} and the growth in $v$ from Assumption \ref{Ass:Monotone Assumption} (as in the previous proof), we have the claim.

\end{proof}

\begin{Proposition}\label{prop:SSTM: converge rate of big thoerm}
Let the assumptions of Theorem \ref{theo:SSTM convergence rate for MV-sde} holds   with $m\ge 2(q+1)^2$. 
Let $i\in\llbracket 1,N\rrbracket$ and take $X^{i,N}$ to be the solution to the interacting particle system \eqref{Eq:MV-SDE Propagation} and let $\hx^{i,N}$ be the continuous-time extension of the SSM given by \eqref{eq: SSTM discrete scheme sde with kappa t}. We then have
\begin{align*}
     \sup_{ i\in \llbracket 1,N \rrbracket}   \bE\Big[\sup_{0\le t \le T} |X_{t}^{i,N}-\hx_{t}^{i,N} |^2 \Big]\le C h.
\end{align*}
\end{Proposition}

\begin{proof}
Take $i\in\llbracket 1,N\rrbracket$, $t\in[0,T]$. 
From \eqref{eq:momentboundParticiInteractingSystem} and \eqref{eq:SSTM:momentbound for split-step time extension process}, both $X^{i,N}$ and $\hx^{i,N}$ have bounded $2p$-moments ($p\ge 2$). 
Define the auxiliary quantity $\Delta X^i:=X^{i,N}-\hx^{i,N}$. It\^o's formula applied to $|X_{t}^{i,N}-\hx_{t}^{i,N}|^2=|\Delta X^i_t|^2$ yields
\begin{align}
    |\Delta X^i_t|^2 \label{eq:Conv of big theom:term 1 }
    =&2\int_0^t \Big\langle  v(s,X_{s}^{i,N},\mu_{s}^{N})-v_h(\kappa(s),\hx_{\kappa(s)}^{i,N},\hat{\mu}^{N}_{\kappa(s)}),\Delta X^i_s \Big\rangle ds  
    \\\label{eq:Conv of big theom:term 2 }
        &+2\int_0^t \Big\langle  b(s,X_{s}^{i,N},\mu_{s}^{N})-b_h(\kappa(s),\hx_{\kappa(s)}^{i,N},\wm_{\kappa(s)}^{N}),\Delta X^i_s\Big\rangle ds
    \\    \label{eq:Conv of big theom:term 3 }
    &+\int_0^t \Big| \sigma(s,X_{s}^{i,N},\mu_{s}^{N})-\sigma_h(\kappa(s),\hx_{\kappa(s)}^{i,N},\wm_{\kappa(s)}^{N})\Big|^2 ds
    \\\label{eq:Conv of big theom:term 4 }
    &+2\int_0^t \Big\langle \Delta X^i_s,\Big(\sigma(s,X_{s}^{i,N},\mu_{s}^{N}) - \sigma_h(\kappa(s),\hx_{\kappa(s)}^{i,N},\wm_{\kappa(s)}^{N})\Big) dW^i_s\Big\rangle.
\end{align}
We analyse the components term by term. Namely, for \eqref{eq:Conv of big theom:term 1 } 
\begin{align*}
     &\Big\langle  v(s,X_{s}^{i,N},\mu_{s}^{N})-v_h(\kappa(s),\hx_{\kappa(s)}^{i,N},\hat{\mu}^{N}_{\kappa(s)}),\Delta X^i_s \Big\rangle
     =
    \Big\langle  v(s,X_{s}^{i,N},\mu_{s}^{N})-v(s,\hx_{s}^{i,N},\hat{\mu}^{N}_{s}),\Delta X^i_s \Big\rangle
    \\
    & \qquad
    +    
    \Big\langle  v(s,\hx_{s}^{i,N},\hat{\mu}^{N}_{s})-v(\kappa(s),\hx_{\kappa(s)}^{i,N},\hat{\mu}^{N}_{\kappa(s)}),\Delta X^i_s \Big\rangle
    +    
    \Big\langle  v(\kappa(s),\hx_{\kappa(s)}^{i,N},\hat{\mu}^{N}_{\kappa(s)})-v_h(\kappa(s),\hx_{\kappa(s)}^{i,N},\hat{\mu}^{N}_{\kappa(s)}),\Delta X^i_s \Big\rangle.
\end{align*}
From the Assumption \ref{Ass:Monotone Assumption}, take supremum over $[0,T]$ and expectations, by the Young's inequality, we have
\begin{align}
\nonumber
    &
    \bE \Bigg[  \sup_{0\le t \le T} \int_0^t   
    \Big\langle  v(s,X_{s}^{i,N},\mu_{s}^{N})-v(s,\hx_{s}^{i,N},\hat{\mu}^{N}_{s}),\Delta X^i_s \Big\rangle  ds\Bigg]
    \\\nonumber
    &\le \bE \Bigg[  \sup_{0\le t \le T} \int_0^t    
    \Big\langle  v(s,X_{s}^{i,N},\mu_{s}^{N})-v(s,\hx_{s}^{i,N},\mu_{s}^{N})+v(s,\hx_{s}^{i,N},\mu_{s}^{N})-v(s,\hx_{s}^{i,N},\hat{\mu}^{N}_{s}),\Delta X^i_s \Big\rangle  ds\Bigg]
    \\\nonumber
    &\le \bE \Bigg[  \sup_{0\le t \le T} \int_0^t  \Big[   
    L_v |\Delta X^i_s|^2+ \frac12 \Big|v(s,\hx_{s}^{i,N},\mu_{s}^{N})-v(s,\hx_{s}^{i,N},\hat{\mu}^{N}_{s})\Big|^2+\frac12 |\Delta X^i_s|^2 \Big] ds\Bigg]
    \\
    \label{eq:se:v-v argument 1}
    &\le \bE \Bigg[  \sup_{0\le t \le T} \int_0^t  \Big[   
    (L_v+\frac12) |\Delta X^i_s|^2+\frac{ L_{\tilde{v}}}{2} W^{(2)}(\mu_{s}^{N},\hat{\mu}^{N}_{s} ) \Big] ds\Bigg]
    \le
    C  \bE \Bigg[  \int_0^T  
    \Big(
    |\Delta X^i_s|^2+\frac{1}{N} \sum_{j=1}^N |\Delta X^j_s|^2
    \Big)  ds \Bigg]
    .
\end{align}

By the $1/2$-H\"older regularity in time and the assumption on $v$, the particles being i.i.d.~and the Cauchy-Schwarz inequality we have 
\begin{align}     
    \nonumber
     &\bE \Bigg[ \sup_{0\le t \le T} \int_0^t   \Big\langle   v(s,\hx_{s}^{i,N},\hat{\mu}^{N}_{s})-v(\kappa(s),\hx_{\kappa(s)}^{i,N},\hat{\mu}^{N}_{\kappa(s)}),\Delta X^i_s \Big\rangle ds  \Bigg]
    \\\nonumber
    &
    \le
    \frac{1}{2} \bE \Bigg[ \int_0^T   \Big|v(s,\hx_{s}^{i,N},\hat{\mu}^{N}_{s})-v(\kappa(s),\hx_{\kappa(s)}^{i,N},\hat{\mu}^{N}_{\kappa(s)})\Big|^2 ds\Bigg]
    + \frac{1}{2} \bE \Bigg[\int_0^T    |\Delta X^i_s|^2 ds\Bigg]
    \\ \nonumber
    &\le
    C h+
    C\bE\Bigg[ 
    \int_0^T  \Bigg(
    \Big(1+|\hx_{s}^{i,N}|^{2q}+ |\hx_{\kappa(s)}^{i,N}|^{2q}\Big) |\hx_{s}^{i,N}-\hx_{\kappa(s)}^{i,N}|^2 +\frac{1}{N}\sum_{j=1}^N  |\hx_{s}^{j,N}-\hx_{\kappa(s)}^{j,N}|^2 \Bigg)
    ds\Bigg] 
    + 
    \frac{1}{2} \bE \Bigg[\int_0^T    |\Delta X^i_s|^2 ds\Bigg]
     \\ \nonumber
    &\le
    C h+ C
    \int_0^T  \Bigg[ 
    \sqrt{ 
    \bE\Big[ \Big(1+|\hx_{s}^{i,N}|^{2q}+ |\hx_{\kappa(s)}^{i,N}|^{2q}\Big)^2\Big]
    \bE\Big[|\hx_{s}^{i,N}-\hx_{\kappa(s)}^{i,N}|^4 \Big]
     }+\frac{1}{N}\sum_{j=1}^N  \bE\Big[|\hx_{s}^{j,N}-\hx_{\kappa(s)}^{j,N}|^2\Big]  ~\Bigg]~ ds
    \\
    \label{eq:se:v-v argument 2} %\nonumber
    &\qquad \quad
    +\frac{1}{2} \bE \Bigg[\int_0^T    |\Delta X^i_s|^2 ds\Bigg]
    \quad \le\quad 
    C h+   \frac{1}{2} \bE \Bigg[\int_0^T    |\Delta X^i_s|^2 ds\Bigg].
\end{align}
where in the last inequality we used H\"older's inequality on the product term in combination with Proposition \ref{prop:SSTM：moment bound for the big theorm time extensions } and \ref{prop:SSTM: local increment error}   with $m\ge 2(q+1)^2$ to guarantee the error satisfies {  $ \bE\big[\,|\hx_{s}^{i,N}-\hx_{\kappa(s)}^{i,N}|^4 \big]\le Ch^2$.} We now make use of Proposition \ref{prop:SSTM:diff functions between orgins and h.s} and arguments similar to those above to deal with the last term of the initial inequality
\begin{align}
\nonumber
    \bE \Bigg[  \sup_{0\le t \le T} &\int_0^t    	
    \Big\langle  v(\kappa(s),\hx_{\kappa(s)}^{i,N},\hat{\mu}^{N}_{\kappa(s)})-v_h(\kappa(s),\hx_{\kappa(s)}^{i,N},\hat{\mu}^{N}_{\kappa(s)}),\Delta X^i_s \Big\rangle ds \Bigg]
     \\
     \nonumber
    & 
    \le
    \frac{1}{2}  \bE \Bigg[  \int_0^T     	\Big[\Big| v(\kappa(s),\hx_{\kappa(s)}^{i,N},\hat{\mu}^{N}_{\kappa(s)})-v_h(\kappa(s),\hx_{\kappa(s)}^{i,N},\hat{\mu}^{N}_{\kappa(s)})\Big|^2+ |\Delta X^i_s|^2\Big] ds\Bigg]
    \\ \label{eq:se:v-vh argument}
    % \nonumber
    &\le  \bE \Bigg[   \int_0^T     
    C\Big[ 1+|\hx_{\kappa(s)}^{i,N}|^{q'} +  \frac{1}{N} \sum_{j=1}^N |\hx_{\kappa(s)}^{j,N}|^{q'}\Big] h^2 ds \Bigg]	
    + \frac{1}{2} \bE \Bigg[\int_0^T    |\Delta X^i_s|^2 ds\Bigg]
    \le
    C h^2+   \frac{1}{2} \bE \Bigg[\int_0^T    |\Delta X^i_s|^2 ds\Bigg]
    .
\end{align}
where $q'$ defined in Proposition \ref{prop:SSTM:diff functions between orgins and h.s} such that $m\ge 2(q+1)^2=q'$ as to guarantee $\bE\big[|\hx_{\kappa(s)}^{i,N}|^{q'}\big]\le C$. We now proceed to estimate the $b$ components. Using Young's inequality, for \eqref{eq:Conv of big theom:term 2 }
\begin{align}
\nonumber
\Big\langle  b(s,X_{s}^{i,N}&,\mu_{s}^{N})-b_h(\kappa(s),\hx_{\kappa(s)}^{i,N},\wm_{\kappa(s)}^{N}),\Delta X^i_s\Big\rangle
 \\
 \nonumber
 &
\le
\frac{1}{2}  \Big|b(s,X_{s}^{i,N},\mu_{s}^{N})-b(\kappa(s),\hx_{\kappa(s)}^{i,N},\wm_{\kappa(s)}^{N})\Big|^2+\frac{1}{2} |\Delta X^i_s|^2 
\\
\label{eq:se:b term}
&\le  
 \Big|b(s,X_{s}^{i,N},\mu_{s}^{N})-b(\kappa(s),\hx_{\kappa(s)}^{i,N},\hm_{\kappa(s)}^{N})\Big|^2
 +\Big|b(\kappa(s),\hx_{\kappa(s)}^{i,N},\hm_{\kappa(s)}^{N})-b_h(\kappa(s),\hx_{\kappa(s)}^{i,N},\wm_{\kappa(s)}^{N})\Big|^2
 +|\Delta X^i_s|^2 .
\end{align}

For the first term above, the Lipschitz condition on $b$, \eqref{sstmmod func prop2} and Proposition \ref{prop:SSTM: local increment error}, yield
\begin{align*}
    \Big|b(s,X_{s}^{i,N},\mu_{s}^{N})- &b(\kappa(s),\hx_{\kappa(s)}^{i,N},\hm_{\kappa(s)}^{N})\Big|^2
        \le C\Bigg[~ h +|X_{s}^{i,N}- \hx_{\kappa(s)}^{i,N}|^2+ \frac{1}{N} \sum_{j=1}^N |X_{s}^{j,N}- \hx_{\kappa(s)}^{j,N}|^2  ~   \Bigg]
       \\
        &\le C\Big[~ h 
        +|\Delta X^i_s|^2
        +|\hx_{s}^{i,N}- \hx_{\kappa(s)}^{i,N}|^2
        + \frac{1}{N} \sum_{j=1}^N |\Delta X^j_s|^2
        +\frac{1}{N} \sum_{j=1}^N |\hx_{s}^{j,N}- \hx_{\kappa(s)}^{j,N}|^2 ~   \Big].
\end{align*}

and similarly we obtain: 
\begin{align}
\nonumber
     \Big|\sigma(s,X_{s}^{i,N},\mu_{s}^{N})-\sigma_h(\kappa(s),\hx_{\kappa(s)}^{i,N},\wm_{\kappa(s)}^{N})\Big|^2
    \le &
    C\Big[~ h 
        +|\Delta X^i_s|^2
        +|\hx_{s}^{i,N}- \hx_{\kappa(s)}^{i,N}|^2
        + \frac{1}{N} \sum_{j=1}^N |\Delta X^j_s|^2
        +\frac{1}{N} \sum_{j=1}^N |\hx_{s}^{j,N}- \hx_{\kappa(s)}^{j,N}|^2 ~   \Big]
        \\
        \label{eq:se:sigma term}
        &
        + 2
        \Big|\sigma(\kappa(s),\hx_{\kappa(s)}^{i,N},\hm_{\kappa(s)}^{N})
        -
        \sigma_h(\kappa(s),\hx_{\kappa(s)}^{i,N},\wm_{\kappa(s)}^{N})
        \Big|^2.
\end{align}

Consider the last term \eqref{eq:Conv of big theom:term 4 }, take expectations, using the BDG inequality, Cauchy-Schwarz inequality and Proposition \ref{prop:SSTM:diff functions between orgins and h.s}, 
\begin{align}
\nonumber
\bE\Bigg[ & \sup_{0\le t \le T} \int_0^t \Big\langle \Delta X^i_s,\Big(\sigma(s,X_{s}^{i,N},\mu_{s}^{N}) - \sigma_h(\kappa(s),\hx_{\kappa(s)}^{i,N},\wm_{\kappa(s)}^{N})\Big) dW^i_s\Big\rangle  ~\Bigg]
\\\nonumber
\le& \bE\Bigg[~ C\Big(  \int_0^T  |\Delta X^i_s|^2~\Big|\sigma(s,X_{s}^{i,N},\mu_{s}^{N}) - \sigma_h(\kappa(s),\hx_{\kappa(s)}^{i,N}, \wm_{\kappa(s)}^{N})\Big|^2 ds \Big)^{\frac{1}{2}}~\Bigg]
\\\nonumber
\le& \bE\Bigg[~ \Big(  \frac{1}{4} \sup_{0\le t \le T}|\Delta X^i_t|^2 ~C\int_0^T  ~\Big|\sigma(s,X_{s}^{i,N},\mu_{s}^{N}) - \sigma_h(\kappa(s),\hx_{\kappa(s)}^{i,N}, \wm_{\kappa(s)}^{N})\Big|^2 ds \Big)^{\frac{1}{2}}~\Bigg]
\\\nonumber
\leq& 
\frac{1}{4} \bE\Big[ \sup_{0\le t \le T} |\Delta X^i_t|^2 ~\Big] 
+ 
C\bE\Bigg[~   \int_0^T \Big|\sigma(s,X_{s}^{i,N},\mu_{s}^{N}) - \sigma_h(\kappa(s),\hx_{\kappa(s)}^{i,N}, \wm_{\kappa(s)}^{N})\Big|^2 ds ~\Bigg]
\\\nonumber
\leq& 
\frac{1}{4} \bE\Big[ \sup_{0\le t \le T} |\Delta X^i_t|^2 ~\Big] 
+ 
C\bE\Bigg[~   \int_0^T \Big(
        h 
        +|\Delta X^i_s|^2
        +|\hx_{s}^{i,N}- \hx_{\kappa(s)}^{i,N}|^2
        + \frac{1}{N} \sum_{j=1}^N |\Delta X^j_s|^2
        +\frac{1}{N} \sum_{j=1}^N |\hx_{s}^{j,N}- \hx_{\kappa(s)}^{j,N}|^2
        ds \Big)
~\Bigg]
\\\nonumber
&\qquad +
\frac{1}{2}\bE\Bigg[~   \int_0^T \Big(
        \Big|\sigma(\kappa(s),\hx_{\kappa(s)}^{i,N},\hm_{\kappa(s)}^{N})
        -
        \sigma_h(\kappa(s),\hx_{\kappa(s)}^{i,N},\wm_{\kappa(s)}^{N})\Big|^2
        ds \Big)
~\Bigg]
\\
\label{eq:se:dw argument}
% \nonumber
\leq& 
\frac{1}{4} \bE\Big[ \sup_{0\le t \le T} |\Delta X^i_t|^2 ~\Big] 
+ 
C\bE\Bigg[~   \int_0^T \Big(
        h 
        +|\Delta X^i_s|^2
        + \frac{1}{N} \sum_{j=1}^N |\Delta X^j_s|^2
        ds \Big)
~\Bigg].
\end{align}
where the last inequality follows by Proposition \ref{prop:SSTM：moment bound for the big theorm time extensions }, \ref{prop:SSTM: local increment error} and \ref{prop:SSTM:diff functions between orgins and h.s}. Similarly, for the terms \eqref{eq:se:b term} and \eqref{eq:se:sigma term}, by Proposition \ref{prop:SSTM:diff functions between orgins and h.s} and similar arguments as in \eqref{eq:se:v-vh argument}, we conclude that 
\begin{align}
\nonumber
\bE\Bigg[  \sup_{0\le t \le T} & \Big(
2\int_0^t \Big\langle  b(s,X_{s}^{i,N},\mu_{s}^{N})-b_h(\kappa(s),\hx_{\kappa(s)}^{i,N},\wm_{\kappa(s)}^{N}),\Delta X^i_s\Big\rangle ds
+\int_0^t \Big| \sigma(s,X_{s}^{i,N},\mu_{s}^{N})-\sigma_h(\kappa(s),\hx_{\kappa(s)}^{i,N},\wm_{\kappa(s)}^{N})\Big|^2 ds \Big)
~\Bigg]
\\ 
\label{eq:se:b sigma arguments}
% \nonumber
& \leq 
C\bE\Bigg[~   \int_0^T \Big(
        h 
        +|\Delta X^i_s|^2
        + \frac{1}{N} \sum_{j=1}^N |\Delta X^j_s|^2
        ds \Big)
~\Bigg].
\end{align}
Gathering all inequalities \eqref{eq:se:v-v argument 1}, \eqref{eq:se:v-v argument 2}, \eqref{eq:se:v-vh argument}, \eqref{eq:se:dw argument} and \eqref{eq:se:b sigma arguments} together, taking supremum on $i$, since the particles are i.i.d., we conclude (where $h$ is the leading term)
\begin{align*}
\sup_{ i\in \llbracket 1,N \rrbracket}\bE\Big[ \sup_{0\le t \le T}|X_{t}^{i,N}-\hx_{t}^{i,N} |^2 ~\Big] 
&\le    
\sup_{ i\in \llbracket 1,N \rrbracket}\bE \Big[~C \int_0^T \Big( h+|X_{s}^{i,N}-\hx_{s}^{i,N}|^2 +\frac{1}{N} \sum_{j=1}^N |X_{s}^{j,N}-\hx_{s}^{j,N}|^2 \Big) ds    ~\Big]
\\
&\le    
Ch +C\int_0^T \sup_{ i\in \llbracket 1,N \rrbracket}\bE \Big[ \sup_{0\le u \le s}|X_{u}^{i,N}-\hx_{u}^{i,N}|^2 ~\Big] ds.
\end{align*}
Gr\"onwall's lemma delivers the final result.
\end{proof}

Now, the proof of the main theorem is concluded as follow.
\begin{proof}[Proof of Theorem \ref{theo:SSTM convergence rate for MV-sde}]
In relation to Points 1, 2 and 3 in the theorem's statement: 
Point 1 follows from Proposition \ref{prop:SSTM：moment bound for the big theorm time extensions }; Point 2 follows from Proposition \ref{prop:SSTM: converge rate of big thoerm}; the last point follows by a straightforward combination of Proposition \ref{Prop:Propagation of Chaos} and Proposition \ref{prop:SSTM: converge rate of big thoerm}.
\end{proof}

\subsection[Proof of the stability Theorem]{Proof of the stability Theorem, Theorem \ref{theo:SSTM:stabilty}}
\label{sec:ProofOfStability}

\begin{proof}[Proof of Theorem \ref{theo:SSTM:stabilty}]
Let $i\in \llbracket 1,N\rrbracket$ and $n\in \mathbb{N}$. From \eqref{sstmmod func prop222} in Lemma \ref{lemma:SSTM:new functions def and properties1}, since the particles are identically distributed, we have
\begin{align*}
     \bE \Big[        |Y_{n}^{i,\star,N}-G_{n}^{i,\star,N}|^2\Big]
    &\le 
    \bE \Bigg[  
    \frac{1 }{1-h  (2L_v+1)}\Big( 
    {   
    |\hx_{n}^{i,N}-\hz_{n}^{i,N}|^2 }+L_{\tilde{v}} h W^{(2)}(\hm^{X,N}_n,\hm^{Z,N}_n)  \Big)\Bigg]
    \\
    &
    \le 
    \frac{1+L_{\tilde{v}} h }{1-h  (2L_v+1)}
     \bE \Big[        |\hx_{n}^{i,N}-\hz_{n}^{i,N}|^2\Big].
\end{align*}
By definition of the SSM, one also has $(Y_{n}^{i,\star,N}-\hx_{n}^{i,N})=hv\left(t_n,Y_{n}^{i,\star,N},\hm^{X,N}_n\right)$ and $(G_{n}^{i,\star,N}-\hz_{n}^{i,N})=hv\left(t_n,G_{n}^{i,\star,N},\hm^{Z,N}_n\right)$. 
Thus, from Assumption \ref{Ass:Monotone Assumption}, for any $n$, the Cauchy-Schwarz inequality yields 
\begin{align*}
&\bE\Big[|\hx_{n+1}^{i,N}-\hz_{n+1}^{i,N} |^2   \Big] 
\\
& %\Big( \Big)
\le \bE \Bigg[ ~|Y_{n}^{i,\star,N}-G_{n}^{i,\star,N}|^2
+2\Big\langle  Y_{n}^{i,\star,N}-G_{n}^{i,\star,N},b(t_n,Y_{n}^{i,\star,N},\hm^{Y,N}_n)-b(t_n,G_{n}^{i,\star,N},\hm^{G,N}_n)     \Big\rangle h
\\
&
\qquad +\Big|b(t_n,Y_{n}^{i,\star,N},\hm^{Y,N}_n)-b(t_n,G_{n}^{i,\star,N},\hm^{G,N}_n)\Big|^2h^2 
+\Big|\sigma(t_n,Y_{n}^{i,\star,N},\hm^{Y,N}_n)-\sigma(t_n,G_{n}^{i,\star,N},\hm^{G,N}_n)\Big|^2(\Delta W_{n}^i)^2   \Bigg]
\\
&
\le  (1+L_\sigma h+L_b h^2 )\bE \Bigg[
|Y_{n}^{i,\star,N}-G_{n}^{i,\star,N}|^2
 \Bigg]
 +
 (L_{\tilde {\sigma}} h+L_{\tilde{b}} h^2 )\bE \Bigg[
\frac1N \sum_{j=1}^N |Y_{n}^{j,\star,N}-G_{n}^{j,\star,N}|^2   \Big]
 \Bigg]
\\
&
\qquad
+ 
2h \sqrt{ \bE \Big[|Y_{n}^{i,\star,N}-G_{n}^{i,\star,N}|^2\Big] }
  \sqrt{
\bE \Bigg[\Big|  b(t_n,Y_{n}^{i,\star,N},\hm^{Y,N}_n)-b(t_n,G_{n}^{i,\star,N},\hm^{Y,N}_n)  \Big|^2\Bigg] ~ } 
\\
&
\qquad
+ 
2h \sqrt{ \bE \Big[|Y_{n}^{i,\star,N}-G_{n}^{i,\star,N}|^2\Big] }
\sqrt{  
\bE \Bigg[\Big|  b(t_n,G_{n}^{i,\star,N},\hm^{Y,N}_n) -  b(t_n,G_{n}^{i,\star,N},\hm^{G,N}_n)
\Big|^2\Bigg] ~} 
\\
&
\le  \Big(1+(2\sqrt{L_b} +2\sqrt{L_{\tilde{b}}}+L_\sigma +L_{\tilde {\sigma}})h+(L_b +L_{\tilde{b}} )h^2 \Big)\bE \Big[
|Y_{n}^{i,\star,N}-G_{n}^{i,\star,N}|^2
 \Big].
\end{align*}
  
where we used the tower property of the expectation with $\cF_{t_n}$-conditional expectations to deal with the Brownian increment term (it holds that 
$\bE[ \,|\Delta W_{n}^i|^2| \cF_{t_n}]=h$ after using that all  $Y_{n}^{j,\star,N},~G_{n}^{j,\star,N}$ are $\cF_{t_n}$-adapted), 
 the Cauchy-Schwarz inequality and that the particles are i.i.d.

Taking supremum over $i$ and using \eqref{sstmmod func prop222} yields  
\begin{align*}
\sup_{ i\in \llbracket 1,N \rrbracket}   \bE\Big[ |\hx_{n+1}^{i,N}-\hz_{n+1}^{i,N} |^2 \Big]
\le   
( 1+\beta h)
\sup_{ i\in \llbracket 1,N \rrbracket}
\bE\Big[ |\hx_{n}^{i,N}-\hz_{n}^{i,N}|^2   \Big].
\end{align*}
where $\beta$ and $\alpha$ are exactly given by \eqref{SSTM:beta formula}.  
A straightforward induction argument leads to 
\begin{align*}
\sup_{ i\in \llbracket 1,N \rrbracket}  \bE\Big[ |\hx_{n}^{i,N}-\hz_{n}^{i,N} |^2 \Big]
\le &     
\big( 1+\beta h\big)^n
\sup_{ i\in \llbracket 1,N \rrbracket} \bE \Big[ |\hx^{i,N}_{0}-\hz^{i,N}_{0} |^2\Big].
\end{align*}

Recall \eqref{SSTM:beta formula} for the expression for $\beta$. 
If $L_v\ge -(1+L_{\tilde{v}}+A)/2$ then $1+\beta h > 1$ and hence $ \lim_{n\rightarrow \infty}   \big( 1+\beta h\big)^n \neq 0$, this implies that the SSM is not Mean-square contractive.  
On the other hand, since $(1+\beta h)$ is always positive then when $\beta<0$ and $(1+\beta h)\in (0,1)$ $\Leftrightarrow  L_v< -(1+L_{\tilde{v}}+A)/2< -\frac{1}{2}$ with sufficient small $h$ and consequently the SSM is guaranteed to be Mean-square contractive
\begin{align*}
\lim_{n\rightarrow \infty}\sup_{ i\in \llbracket 1,N \rrbracket}  \bE\Big[ |\hx_{n}^{i,N}-\hz_{n}^{i,N} |^2 \Big]
\le &  
\lim_{n\rightarrow \infty}   ( 1+\beta h)^n
\sup_{ i\in \llbracket 1,N \rrbracket} \bE \Big[ |\hx^{i,N}_{0}-\hz^{i,N}_{0} |^2\Big]   =0.
\end{align*}

\end{proof}

%%%%%%%%%%%%%%%%%%%%%%%%%%%%%%%%%%%%%
%%%%%%%%%%%%%%%%%%%%%%%%%%%%%%%%%%%%%%%%%%%%%%%%%%%
%%%%%%%%%%%%%%%%%%%%%%%%%%%%%%%%%%%%%%%%%%%%%%%%%%%%%%%%%%%%%%
%%%%%%%%%%%%%%%%%%%%%%%%%%%%%%%%%%%%%

%%%%%%%%%%%%%%%%%%%%%%%%%%%%%%%%%%%%%
%%%%%%%%%%%%%%%%%%%%%%%%%%%%%%%%%%%%%%%%%%%%%%%%%%%
%%%%%%%%%%%%%%%%%%%%%%%%%%%%%%%%%%%%%%%%%%%%%%%%%%%%%%%%%%%%%%

%%%%%%%%%%%%%%%%%%%%%%%%%%%%%%%%%%%%%%%%%%%%%%%%%%%%%%%%%%%%%
%%%%%%%%%%%%%%%%%%%%%%%%%%%%%%%%%%%%%%%%%%%%%%%%%%%%%%%%%%%%%
%%%%%%%%%%%%%%%%%%%%%%%%%%%%%%%%%%%%%%%%%%%%%%%%%%%%%%%%%%%%%
%%%% \BEGIN APPENDIX
%%%%%%%%%%%%%%%%%%%%%%%%%%%%%%%%%%%%%%%%%%%%%%%%%%%%%%%%%%%%%
% \newpage
\appendix

\section{Short description of the Taming and Adaptive time-stepping method}
\label{sec:ShortRecap}

We provide a brief review of Taming \cite{reis2018simulation} and Adaptive time-stepping  \cite{reisinger2020adaptive} for superlinear growth MV-SDEs in the context of the notation set in Section \ref{sec:one} and \ref{sec:two}. Each method approximates \eqref{Eq:General MVSDE} through the interacting particle system \eqref{Eq:MV-SDE Propagation} as described next. Table \ref{table with Schemes} summarises strong error (rMSE), mentions weak error as an open problem,  and stability results. Both schemes (plus the proposed SSM one) hold  under the same conditions: Assumption \ref{Ass:Monotone Assumption} and a sufficiently high integrable initial condition (as in Theorem \ref{theo:SSTM convergence rate for MV-sde}).

\begin{table}[h!tb]
\centering
\begin{tabular}{c|c|c}
% \hline
Methods  &  rMSE & Stability \\ 
\hline
\hline 
Taming \cite{reis2018simulation}
&
    0.5
    &  Unknown    
    \\ 
   \hline
Adaptive \cite{reisinger2020adaptive} 
& 0.5      &   Unknown      \\ \hline
SSM     
& 0.5         &   Contraction Theorem \ref{theo:SSTM:stabilty}   \\ \hline
\end{tabular}
\caption{Information regarding the different methods. Convergence of Taming \cite{reis2018simulation} and Adaptive \cite{reisinger2020adaptive} scheme under the same conditions as Assumption \ref{Ass:Monotone Assumption}. The Root MSE (rMSE) is the metric presented in \eqref{eq:rMSEofSSM}. Weak error analysis has not been carried out but experimental work points in the direction of weak error order $1$ for the three methods.
}
\label{table with Schemes}
\end{table}

\subsection{The Taming method}
Taming \cite{reis2018simulation} approximates \eqref{Eq:MV-SDE Propagation} as follows (see also \cite[Section 4]{reisinger2020adaptive})
\begin{align}
	\label{Eq:Tamed MVSDE}
 \bar{X}_{n+1}^{i,N,M}
	=
	\bar{X}_{n}^{i,N,M}
	&+ 
	\dfrac{\widehat b\Big(t_{n} , \bar{X}_{n}^{i,N,M}, \bar{\mu}^{X,N}_{n} \Big) }
	{1+ M^{-\alpha} \Big\vert \widehat b\Big(t_{n} , \bar{X}_{n}^{i,N,M}, \bar{\mu}^{X,N}_{n} \Big) \Big\vert } h
	%\\
	%& \qquad \qquad \nonumber
	+ 
	\sigma\Big(t_{n} , \bar{X}_{n}^{i,N,M}, \bar{\mu}^{X,N}_{n} \Big)
	\Delta W_{n}^{i} ,\qquad i\in \llbracket 1,N\rrbracket.
\end{align}
where $\bar{\mu}^{X,N}_{n}(\dd x) = \frac{1}{N} \sum_{j=1}^N \delta_{\bar{X}_{n}^{j,N,M}}(\dd x)$,$~\Delta W_{n}^{i}=W_{t_{n+1}}^{i}-W_{t_{n}}^{i} $ with $\bar{X}_{0}^{i,N,M}= X_{0}^{i}$. The parameter $\alpha \in (0, 1]$ is a tuning parameter  where setting $\alpha=1/2$ delivers a rMSE convergence rate of order $1/2$ while setting $\alpha=1$ delivers a rMSE convergence rate of order $1$ (for a constant diffusion $\sigma$).

%%%%%%%%%%%%%%%%%%%%%%%%%%%%%%%%%%%%%%%%%%%%%%%%
%%%%%%%%%%%%%%%%%%%%%%%%%%%%%%%%%%%%%%%%%%%%%%%%%%%
%%%%%%%%%%%%%%%%%%%%%%%%%%%%%%%%%%%%%%%%%%%%%%%%%
% \newpage

\subsection{Adaptive time-stepping method}

Adaptive \cite{reisinger2020adaptive} approximates \eqref{Eq:MV-SDE Propagation} as follows for  
% With the notation above, set , $M h=T,~i\in\{1,\cdots,N\}$ and 
$t_n\in[k_n h,(k_n+1)h),k_n\in\mathbb{N}$ and
% , $~n\in\llbracket 0,M-1\rrbracket,~X\in \mathbb{R}^d$ reads:
\begin{align}
	\label{Eq:Adaptive MVSDE}
 \bar{X}_{t_{n+1}}^{i,N}
	=
	\bar{X}_{t_{n}}^{i,N}
	&+ 
	\widehat b\Big(t_{n} , \bar{X}_{t_{n}}^{i,N}, \bar{\mu}^{X,N}_{k_n h} \Big) 
	 h_n^i
	%\\
	%& \qquad \qquad \nonumber
	+ 
	\sigma\Big(t_{n} , \bar{X}_{t_{n}}^{i,N}, \bar{\mu}^{X,N}_{k_n h} \Big)
	\Delta W_{t_{n}}^{i} ,\qquad i\in \llbracket 1,N\rrbracket. 
\end{align}
where $\bar{\mu}^{X,N}_{k_n h}(\dd x) = \frac{1}{N} \sum_{j=1}^N \delta_{\bar{X}_{k_n h}^{j,N,M}}(\dd x)$, $t_{n+1}=t_n+h_n^i$, $~\Delta W_{t_{n}}^{i}=W_{t_{n+1}}^{i}-W_{t_{n}}^{i} $ with $\bar{X}_{0}^{i,N,M}= X_{0}^{i}$ and for a map $\mathbf{h}^{\delta}(x):\mathbb{R}^d\rightarrow [0,h]$ 
\[h_n^i=\min\Big\{\mathbf{h}^{\delta}(\bar{X}_{t_{n}}^{i,N}),(k_n+1)h-t_n  \Big\}.\]

The function $\mathbf{h}^{\delta}$ is specified at each example and is to be understood similarly to the taming technique. In essence, $\mathbf{h}^{\delta}$ is to be chosen such that $|\widehat b(x) \mathbf{h}^{\delta}(x)|$ is of linear growth.  
For Adaptive, one modifies the timestep $h$ in a dynamic fashion to control the growth of $\widehat b$ while taming modifies the drift $\widehat b$ to control the growth across the application of the scheme. The rMSE convergence rate of order $1/2$, see \cite{fang2020adaptive} or  \cite{reisinger2020adaptive}.

\bibliographystyle{abbrv}
% \bibliography{BibFileGoesHere}  

\end{document}